\documentclass[12pt,a4paper]{article}\usepackage[top=20mm,left=12mm]{geometry}\textwidth18.2cm\textheight25.3cm
\usepackage{amsmath,amssymb,amsthm,amsfonts}
\usepackage{lineno,url}\usepackage{float}
\floatstyle{plaintop}\restylefloat{table}
\usepackage[tableposition=top]{caption}
\usepackage{graphicx}
\usepackage{multirow,pstricks,fancyvrb,listings}
\usepackage{paralist, longtable}
\usepackage{setspace}
\newtheoremstyle{myplain}
  {-\baselineskip\topsep}   
  {\topsep}   
  {\itshape\setstretch{1.05}}  
  {0pt}       
  {\bfseries} 
  {.}         
  {5pt plus 1pt minus 1pt} 
  {}       
\newtheoremstyle{mydef}
  {-\baselineskip\topsep}   
  {\topsep}   
  {\setstretch{1.05}}  
  {0pt}       
  {\bfseries} 
  {.}         
  {5pt plus 1pt minus 1pt} 
  {}       
\theoremstyle{myplain}
\newtheorem{theorem}{Theorem}[section]
\newtheorem{remark}[theorem]{Remark}
\def\bexe{\begin{exercise}}\def\eexe{\eex\end{exercise}}
\def\bsol{\begin{solution}}\def\esol{\eex\end{solution}}
\def\bexa{\begin{example}}\def\eexa{\end{example}}
\def\brem{\begin{remark}}\def\erem{\end{remark}}
\def\bthm{\begin{theorem}}\def\ethm{\end{theorem}}
\def\blem{\begin{lemma}}\def\elem{\end{lemma}}
\def\bcor{\begin{corollary}}\def\ecor{\end{corollary}}
\def\bdefi{\begin{definition}}\def\edefi{\end{definition}}
\newcommand{\IDEA}{\textbf{Idea of the Proof.} }

\def\bmip{\begin{minipage}{\textwidth}}\def\emip{\end{minipage}}
\def\huga#1{\begin{gather} #1 \end{gather}}

\def\hual#1{\begin{align} #1 \end{align}}

\newcommand{\R}{{\mathbb R}}
\newcommand{\C}{{\mathbb C}}\newcommand{\N}{{\mathbb N}}
\newcommand{\Q}{{\mathbb Q}}\newcommand{\Z}{{\mathbb Z}}

\def\CA{{\cal A}}\def\CD{{\cal D}}  
\def\CG{{\cal G}}\def\CH{{\cal H}}
\def\CO{{\cal O}}
\def\CM{{\cal M}}

\def\setm{\setminus}
\def\spani{{\rm span}}

\def\vx{\vec x}

\def\uti{\tilde{u}}

\def\ga{\gamma}\def\om{\omega}
\def\noi{\noindent}\def\ds{\displaystyle}
\def\vt{\vartheta}\def\pa{{\partial}}\def\lam{\lambda}
\newcommand{\bi}{\begin{itemize}}\newcommand{\ei}{\end{itemize}}
\newcommand{\ben}{\begin{enumerate}}\newcommand{\een}{\end{enumerate}}
\newcommand{\bce}{\begin{center}}\newcommand{\ece}{\end{center}}
\newcommand{\bci}{\begin{compactitem}}\newcommand{\eci}{\end{compactitem}}
\newcommand{\bcen}{\begin{compactenum}}\newcommand{\ecen}{\end{compactenum}}
\newcommand{\bcena}{\begin{compactenum}[(a)]}
\newcommand{\reff}[1]{(\ref{#1})}

\newcommand{\spr}[1]{\left\langle #1 \right\rangle}
\newcommand{\hs}[1]{{\hspace{#1}}}\newcommand{\vs}[1]{{\vspace{#1}}}

\def\eps{\varepsilon}
\def\ra{\rightarrow}

\newcommand{\barr}{\begin{array}}\newcommand{\earr}{\end{array}}
\newcommand{\bpm}{\begin{pmatrix}}\newcommand{\epm}{\end{pmatrix}}
\newcommand{\bsm}{\left(\begin{smallmatrix}}
\newcommand{\esm}{\end{smallmatrix}\right)}
\newcommand{\ba}{\begin{array}}\newcommand{\ea}{\end{array}}
\def\dd{\, {\rm d}}\def\ri{{\rm i}}

\def\er{{\rm e}}
\def\re{{\rm Re}}
\def\om{\omega}\def\Om{\Omega}
\def\hot{{\rm h.o.t}}\def\ddt{\frac{\rm d}{{\rm d}t}}
\def\del{\delta}

\def\eex{\hfill\mbox{$\rfloor$}}

\def\al{\alpha}

\def\Ga{\Gamma}

\def\bd{\begin{displaymath}} \def\ed{\end{displaymath}}
\def\ba{\begin{array}} \def\ea{\end{array}}
  
\def\eps{\varepsilon}

\def\pdep{{\tt pde2path}}

\def\p{\tt p} 
\def\oop{{\tt OOPDE}}

\def\mlab{{\tt Matlab}}\def\ptool{{\tt pdetoolbox}}

\def\dhome{/hh/path/pde2path/demos}
\def\hdhome{./hopfdemos}
\def\Gbc{G_{\text{BC}}}

\definecolor{codegreen}{rgb}{0,0.6,0}
\definecolor{codegray}{rgb}{0.5,0.5,0.5}
\definecolor{codepurple}{rgb}{0.58,0,0.82}
\definecolor{backcolour}{rgb}{0.95,0.95,0.92}
 
\lstdefinestyle{mystyle}{
    backgroundcolor=\color{backcolour},   
    commentstyle=\color{codegreen},
    keywordstyle=\color{black},
    numberstyle=\small\color{codegray},
    stringstyle=\color{codepurple},
    basicstyle=\footnotesize\ttfamily,
    breakatwhitespace=false,         
    breaklines=true,                 
    captionpos=b,                    
    keepspaces=true,                 
    numbers=left,                    
    numbersep=5pt,                  
    showspaces=false,                
    showstringspaces=false,
    showtabs=false,                  
    tabsize=2, 
stepnumber=0, 
  xleftmargin=4mm,
}
 
\newlength{\tew}

\def\neig{n_\text{eig}}\def\ig{\includegraphics}
\def\rds{{\rm ds}}
\def\hoxi{\xi_{\text{H}}}
\def\ut{\tilde{u}}\def\ind{{\rm ind}}\def\emu{{\rm err}_{\ga_1}}

\def\fltol{{\rm tol}_{{\rm fl}}}
\def\heda{{\bf HD1}}\def\hedb{{\bf HD2}}
\def\fla{{\bf FA1}}\def\flb{{\bf FA2}}\def\rO{{\rm O}}\def\sO{{\rm SO}}

\usepackage{hyperref}
\numberwithin{equation}{section}
\renewcommand{\arraystretch}{1.15}\renewcommand{\baselinestretch}{1.1}
\def\taskip{\renewcommand{\arraystretch}{1}\renewcommand{\baselinestretch}{1}}
\def\teskip{\renewcommand{\arraystretch}{1.1}\renewcommand{\baselinestretch}{1.1}}
\def\hulst#1#2{\taskip\lstinputlisting[#1]{#2}\teskip}
\def\hutab#1{\taskip\begin{\table}#1\end{table}\teskip}

\begin{document}
\text{}\vspace{10mm}
\begin{center}\Large
User guide on Hopf bifurcation and time periodic orbits with \pdep\\[4mm]
\normalsize Hannes Uecker \\[2mm]
\small 
Institut f\"ur Mathematik, Universit\"at Oldenburg, D26111 Oldenburg, 
hannes.uecker@uni-oldenburg.de\\[2mm]
\normalsize
\today
\end{center}
\noi
\begin{abstract} 
  We explain the setup for using the \pdep\ libraries for Hopf
  bifurcation and continuation of branches of periodic orbits and give
  implementation details of the associated demo directories. See
  \cite{hotheo} for a description of the basic algorithms and the
  mathematical background of the examples.  Additionally we explain
  the treatment of Hopf bifurcations in systems with continuous
  symmetries, including the continuation of traveling waves and
  rotating waves in $\rO(2)$ equivariant systems as relative
  equilibria, the continuation of Hopf bifurcation points via extended
  systems, and some simple setups for the bifurcation from periodic
  orbits associated to critical Floquet multipliers going through
  $\pm 1$.
\end{abstract}

\setcounter{tocdepth}{2}
\tableofcontents 

\section{Introduction}\label{i-sec}
In \cite{hotheo} we describe the basic algorithms in \pdep\ to study 
Hopf bifurcations%
\footnote{i.e.: the bifurcation of (branches of) 
time periodic orbits (in short Hopf orbits) 
from steady states; accordingly, we shall call these branches Hopf branches;}
in PDEs of the form 
\huga{\label{tform} 
M_d\pa_t u=-G(u,\lam), \quad u=u(x,t),\ x\in\Om,\ t\in\R, 
}
where $M_d\in\R^{N\times N}$ is the dynamical mass matrix, and 
$G(u,\lam):=-\nabla\cdot(c\otimes\nabla u)+a u-b\otimes\nabla u-f$. 
Here $u=u(x,t)\in\R^N$, $x\in\Omega$ with $\Om\subset\R^d$ some bounded domain, 
$d=1,2,3$, $\lam\in\R^p$ is a parameter (vector), 
and the diffusion, advection and linear tensors $c,b,a$, and the nonlinearity 
$f$, can  depend on $x,u,\nabla u$, and parameters. 
The  boundary conditions (BC) for \reff{tform} are of 
``generalized Neumann'' form, i.e., 
\begin{align}\label{gnbc}
{\bf n}\cdot (c \otimes\nabla u) + q u = g,
\end{align}
where ${\bf n}$ is the outer normal and again $q\in \R^{N\times N}$
and $g\in \R^N$ may depend on $x$, $u$, $\nabla u$ and
parameters, and over rectangular domains there additionally is the 
possibility of periodic BC in one or more directions. See \cite{hotheo}, 
and the steady state tutorials at \cite{p2phome}, for more details 
on $c,b,\ldots,g$. 

\pdep\ spatially discretizes the PDE \reff{tform}, \reff{gnbc} via piecewise 
linear finite elements, leading to the high--dimensional ODE problem 
(with a slight misuse of notation) 
\huga{\label{tformd} 
M\dot{u}=-G(u,\lam), \quad u=u(t)\in\R^{n_u}, \quad G(u)=Ku-Mf(u),  
}
where $n_u=n_pN$ is the number of unknowns, with $n_p$ the number of 
grid points. In \reff{tformd}, $M$ is the mass matrix, $K$ is the 
stiffness matrix, which typically corresponds to the diffusion 
term $-\nabla\cdot(c\otimes\nabla u)$, and $Mf:\R^{n_u}\ra\R^{n_u}$ contains  the rest, which we often also 
call the 'nonlinearity'. However, \reff{tformd} is really a sort of 
symbolic notation to express the spatially discretized version of 
\reff{tform}, and $K$ in \reff{tformd} can also involve nonlinear terms 
and first order derivatives coming from $b\otimes\nabla u$ in \reff{tform}. 
See, e.g., \cite{actut,symtut} for more details. 

Here we first present implementation details for the four Hopf bifurcation test problems considered in \cite{hotheo}, thus giving a tutorial on how to treat 
Hopf bifurcations in \pdep. See  \cite{p2phome}  
for download of the package, including the demo directories, 
and various documentation and tutorials. In particular, 
since the Hopf examples are somewhat more involved than the steady case we 
recommend to new users of \pdep\ to first look into \cite{actut}, 
which starts with some simple steady state problems. 

The first Hopf demo problem (demo {\tt cgl}, subdir of {\tt hopfdemos}) is a cubic--quintic  complex Ginzburg--Landau (cGL) equation,  
which we consider over 1D, 2D, and 3D cuboids with homogeneous 
Neumann and Dirichlet BC. Next we consider a 
Brusselator system (demo {\tt brussel}) from \cite{yd02}, which shows interesting 
interactions between Turing branches and Turing--Hopf branches.  
As a non--dissipative example we treat the canonical system 
for an optimal control problem (see also \cite{octut}) of ``optimal pollution''  (demo {\tt pollution}). This is still of 
the form \reff{tform}, but is ill--posed as an initial value problem, 
since it features ``backward diffusion''. Nevertheless, we continue 
steady states and obtain Hopf bifurcations 
and branches of periodic orbits. 

In \S\ref{hosymsec} 
we give three tutorial examples for Hopf bifurcations in systems with 
continuous symmetries, namely a reaction-diffusion problem with mass conservation (demo {\tt mass-cons}), 
a Kuramoto-Sivashinsky equation with translational and boost invariance 
(demos {\tt kspbc2} and {\tt kspbc4}), 
and a FHN type system with (breathing) pulses featuring an approximate 
translational symmetry (demo {\tt symtut/breathe}). 
Such symmetries were not considered systematically in \cite{hotheo}.  
They require phase conditions, first for the reliable 
continuation of steady states and detection of (Hopf or steady) bifurcations, 
 which typically lead to the coupling of 
\reff{tform} with algebraic equations. To compute Hopf branches
we then also need to set up suitable phase conditions 
for the Hopf orbits. See also \cite{symtut}, where preliminary results 
for the  breathing pulses are discussed, and one more example of 
Hopf orbits for systems with symmetries is considered, namely modulated 
traveling fronts in a combustion model.  
Moreover, additional to \cite{hotheo,symtut} we explain  some further routines such as continuation of Hopf bifurcation points, and some 
simple setups for branch switching {\em from} Hopf orbits, i.e., 
for Hopf pitchforks (critical multiplier $1$) and period doubling 
(critical multiplier $-1$). See also \cite[\S8]{pftut} for 
further examples of period--doubling bifurcations in a classical 
two--component Brusselator system, following \cite{YZE04}. 

In \S\ref{o2sec} we consider Hopf bifurcation in $\rO(2)$ 
equivariant systems, which generically leads to 
coexistence of standing waves (SWs) 
and traveling waves (TWs). We start with the cGL equation over an interval 
with periodic BC (pBC) (demo {\tt cglpbc}) and use an appropriate additional phase 
condition to continue TWs, periodic in the lab frame, but steady in the co--moving frame) as relative 
equilibria, from which we obtain modulated TWs via Hopf bifurcation in 
the co--moving frame. A similar approach for the cGL equation in a disk 
with Neumann BC (demo {\tt cgldisk}) yields spiral waves as rotating waves (RWs), and meandering 
spirals as modulated RWs. Then we review and extend an example 
from \cite[\S3.2]{hotheo}, namely a reaction diffusion system  
in a disk and with Robin BC (demo {\tt gksspirals}), following 
\cite{GKS00}. In all these examples, the 'interesting' Hopf 
bifurcation points have multiplicity at least two, leading to SWs vs TWs (or RWs). 
Our default branch switching so far only deals with 
simple Hopf bifurcation points systematically, but we use and
 explain an ad hoc modification for Hopf points of higher multiplicity, 
allowing to select, e.g., TW vs SW branches.  At the end of \S\ref{o2sec} 
we also explain a setup to compute periodic orbit branches with 
fixed period $T$ (freeing a second parameter), and for non--autonomous 
systems, i.e., with explicit $t$--dependence of $G$.

The user interfaces for Hopf bifurcations reuse the standard 
\pdep\ setup and no new user functions are necessary, except for the case 
of symmetries, which requires one additional user function. 
For the basic ideas of continuation and bifurcation 
in steady problems we refer to \cite{p2p} (and the references therein), 
for \pdep\ installation hints and review of data structures to \cite{qsrc}, 
for a general soft introduction to \cite{actut}, and for the basic algorithms implemented in the {\tt hopf} library of \pdep\ to \cite{hotheo}. 
Thus, here we concentrate on how to use these routines, and on recent 
additions.  

The basic setup of all demos is similar. Each demo directory contains: 
\bci 
\item 
Function files 
named {\tt *init.m} for setting up the geometry and the basic \pdep\ 
data, where  {\tt *} stands for the problem, e.g., {\tt cgl} (and later 
{\tt brussel, \ldots}). 
\item Main script files, such as {\tt cmds*d.m} 
where * stands for the space dimension. 
\item Function files {\tt sG.m} and {\tt sGjac.m} 
for setting up the rhs of the equation and its Jacobian. 
Most examples are 2nd order semilinear, i.e., of the form 
$\pa_t u=-G(u)=D\Delta u+f(u)$ with diffusion matrix $D\in\R^{N\times N}$, and 
we put the 'nonlinearity' $f$ (i.e., everything except diffusion) into a function  {\tt nodalf.m}, which is then 
called in {\tt sG.m}, but also in mesh-adaption and in 
normal form computations. An exception is the KS equation, 
see \S\ref{kssec}. 

\item  a function {\tt oosetfemops.m} for setting up the system matrices.  
\item 
A few auxiliary functions, for instance small 
modifications of the basic plotting routine {\tt hoplot.m} from 
the {\tt hopf} library, which we found convenient to have problem adjusted 
plots. 
\item Some auxiliary scripts {\tt auxcmds.m}, which 
contain commands, for instance for creating movies or 
for mesh--refinement, which are not genuinely related to the Hopf computations,   
but which we find 
convenient for illustrating either some mathematical aspects of the models, 
or some further \pdep\ capabilities, and which we hope the user will 
find useful as well. 
\item For the demo {\tt pollution} (an optimal control problem) we additionally have the functions {\tt polljcf.m}, which 
implements the objective function, and {\tt pollbra.m}, which combines 
output from the standard Hopf output {\tt hobra.m} and the standard OC 
output {\tt ocbra.m}. 
\eci 

In all examples, the meshes are chosen rather coarse, to quickly 
get familiar with the software. 
We did check for all examples that these coarse meshes give reliable 
results by running the same simulations on finer meshes, without 
qualitative changes. We give hints about the timing and 
indications of convergence, but we refrain from a genuine 
convergence analysis. 
In some cases (demos {\tt cgl} in 3D and {\tt brussel} in 2D, 
and {\tt cgldisk, gksspirals}) we additionally 
switch off the on the fly computation of Floquet multipliers
and instead compute the multipliers a posteriori 
at selected points on branches. Computing the multipliers 
simultaneously is possible as well, but may be slow.  
 Nevertheless, even with the coarse meshes  
some commands, e.g., the continuation 
of Hopf branches in 3+1D (with about 120000 total degrees of freedom), may take several minutes. 
All computational times given in the following 
are from a 16GB i7 laptop with Linux Mint 18 and \mlab\ 2016b. 

In \S\ref{exsec} to \S\ref{o2sec} we explain the 
implementations for the four demos from \cite{hotheo} and the extensions, and thus 
in passing also the main data-structures and functions from the {\tt hopf} 
library. In particular, in \S\ref{brusec} we extend some results from \cite{hotheo} by Hopf point continuation, and in \S\ref{hosymsec} explain the setup for the systems 
with symmetries and hence constraints. For the unconstrained theory, 
and background on the first three example problems, we refer to \cite{hotheo}, 
but for easier reference and to explain the setup with constraints,  we also 
give the pertinent formulas in Appendix A. This also helps understanding a  
number of new functions, for instance for continuation of TWs and RWs, and for bifurcation from periodic orbits, and Appendix B contains an overview of the functions in the {\tt hopf} library and of the relevant data 
structures for quick reference. 
For comments, questions, and bugs, please mail to 
{\tt hannes.uecker@uni-oldenburg.de}. 

\vs{3mm}
\noindent 
{\bf Acknowledgment.}  Many thanks to Francesca Mazzia 
for providing TOM \cite{MT04}, which was 
essential help for setting up the {\tt hopf} library; to Uwe Pr\"ufert 
for providing \oop; 
to Tomas Dohnal, Jens Rademacher and Daniel Wetzel for some 
testing of the Hopf examples; to Daniel Kressner for {\tt pqzschur}; 
and to Dieter Grass for the cooperation 
on distributed optimal control problems, which was one of my main 
motivations to implement 
the {\tt hopf} library. 

\section{The cGL equation as an introductory example: Demo {\tt cgl}}\label{exsec}\label{cglsec}
\def\dname{cgl}
We consider the cubic-quintic complex Ginzburg--Landau equation 
\huga{\label{cAC0} 
\pa_t u=\Delta u+(r+\ri\nu)u-(c_3+\ri \mu)|u|^2 u-c_5|u|^4u, \quad 
u=u(t,x)\in\C, 
}
with real parameters $r,\nu,c_3,\mu,c_5$. 
In real variables $u_1,u_2$ with $u=u_1+\ri u_2$, 
\reff{cAC0} can be written as the 
2--component system  
\huga{\label{cAC} 
\pa_t \bpm u_1\\ u_2\epm =\bpm \Delta+r&-\nu\\\nu&\Delta+r\epm
\bpm u_1\\ u_2\epm-(u_1^2+u_2^2)\bpm c_3 u_1-\mu u_2\\ 
\mu u_1+c_3 u_2\epm-c_5(u_1^2+u_2^2)^2\bpm u_1\\ u_2\epm. 
}
We set 
$\text{$c_3=-1, c_5=1, \nu=1, \mu=0.1$, }$ 
and use $r$ as the main bifurcation parameter.  Considering 
\reff{cAC} on, e.g., a (generalized) rectangle 
\huga{\label{omdef}
\Om=(-l_1\pi,l_1\pi)\times\cdots\times(-l_d\pi,l_d\pi) 
} 
with homogeneous Dirichlet BC or Neumann BC, or with periodic BC, 
we can explicitly calculate all 
Hopf bifurcation points from the trivial branch $u=0$, 
located at $r=|k|^2:=k_1^2+\ldots+k_d^2$, with wave vector $k\in \Z/(2l_1)\times 
\ldots\times \Z/(2l_d)$. 

In particular from the bifurcation point of view, an important feature 
of the cGL equation \reff{cAC0} are its various symmetries, as also 
discussed in \cite{symtut}. For homogeneous Neumann and Dirichlet BCs, 
\reff{cAC0} has the gauge (or rotational) symmetry $u\mapsto \er^{\ri \phi} u$, 
i.e., \reff{cAC0} is equivariant wrt the action of the special 
orthogonal group $\sO(2)$. Periodic boundary conditions imply a translation invariance as an in general additional $\sO(2)$ equivariance, as on the real line.  However, for wavetrains (or traveling waves) $u(t,x)=R\exp(\ri(kx-\omega t))$, where $R>0$ and $\omega,k\in\R$ are the amplitude, frequency and wave number, respectively,  the rotation and translation have the same group orbits. Therefore, in contrast 
to the steady case \cite{symtut}, for our purposes here 
the gauge symmetry will not play a role. Additional, \reff{cAC0} has 
the reflection symmetry $x\mapsto -x$ (1D, and analogously for suitable BC 
over higher dimensional boxes with suitable BC). Thus, in summary, for pBC 
(and also for the cGL in a disk with e.g., Neumann BC where the role 
of spatial translation will be played by spatial rotation, 
the pertinent symmetry group is $\rO(2)$. Consequently, many Hopf bifurcation 
points from the trivial solution will be at least double, and in 
\S\ref{o2sec} we discuss 
this case and the associated questions of the bifurcation of standing vs 
traveling waves, and their numerical treatment in \pdep. 

Here we shall first focus on \reff{cAC0} over boxes with BC that break the 
translational invariance, such that only discrete symmetries remain, 
since, as noted above, for Hopf bifurcation the gauge invariance is equivalent to time shifts, and hence is automatically factored out by the phase condition 
for time shifts. Moreover, we shall choose boxes such that all HBPs are 
simple. 

\subsection{General setup}
The cGL demo directory consists, as noted above, 
of some function files to set up and 
describe \reff{cAC}, some script files to run the computations, 
and a few auxiliary functions and scripts to explain additional 
features, or, e.g., to produce customized plots. An overview 
is given in Table \ref{tab1}, and for this 'first' demo we discuss 
the main files in some detail, while in later demos we will mainly 
focus on differences to this basic template. 

\taskip
\begin{table}[ht]\caption{{\small Scripts and functions in {\tt hopfdemos/cgl}. 
Treating the 1D, 2D and 3D cases in one directory, the only dimension dependent 
files are the scripts, and the function cGLinit. The 2nd part of the table contains 
auxiliary functions and scripts which are not needed for the Hopf 
computations, but which illustrate additional \pdep\ features. }
\label{tab1}}
\bce 
\vs{-5mm}
{\small 
\begin{tabular}{l|p{0.7\tew}}
script/function&purpose,remarks\\
\hline
cmds*d&main scripts; for $*=1,2,3$, respectively, which are all quite similar, i.e., mainly differ in settings for output file names. 
Thus, only cmds1d is discussed in some detail below.\\%
p=cGLinit(p,lx,nx,par,ndim)&init function, setting up parameters and function 
handles, and, as the only space dimension $d$=ndim dependent points, 
the domain.  \\
p=oosetfemops(p)&set FEM matrices (stiffness K and mass M)\\
r=sG(p,u)&encodes $G$ from \ref{cAC} (including the BC)\\
f=nodalf(p,u)&the 'nonlinearity' in \reff{cAC}, i.e., everything except 
$D\Delta u$.\\
Gu=sGjac(p,u)&Jacobian $\pa_u G(u)$ of $G$. \\
\hline
auxcmds1&script, auxiliary commands, illustrating stability checks by time-integration, and a posteriori computation of Floquet multipliers\\
auxcmds2&script, auxiliary commands, illustrating (adaptive) mesh-refinement 
by either switching to natural parametrization, or via {\tt hopftref} \\
cmds1dconv&script showing some convergence of periods in 1D for 
finer $t$ discretization\\
plotana1&plot analytical Hopf-branch for cGL for comparison with numerics, used in {\tt cmds1d.m}, calls mvu=anafloq(rvek,s) \\
homov2d, homov3d&auxiliary functions to generate movies of Hopf orbits
\end{tabular}
}
\ece 
\end{table}\teskip

As main functions we have 
\bci
\item {\tt cGLinit.m}, which (depending on the spatial dimension) 
sets up the domain, mesh, boundary conditions, 
and sets the relevant function handles {\tt p.fuha.sG=@sG} and {\tt p.fuha.sGjac=@sGjac} to encode the rhs of \reff{cAC}; 
\item {\tt sG.m} and {\tt sGjac}, which encode \reff{cAC} and the associated 
Jacobian of $G$; 
\eci 
Then we have three script files, {\tt cmds*d.m}, where 
*=1,2,3 stands for the spatial dimension. These are all 
very similar, i.e., only differ in file names for output and some plotting 
commands, but the basic procedure is always the same: 
\bcen \item call {\tt cGLinit}, then {\tt cont} to find the 
HBPs from the trivial branch $u\equiv 0$, $r\in\R$; 
\item compute branches of periodic orbits (including Floquet multipliers)  
by calling {\tt hoswibra} and {\tt cont} again, 
then plotting. 
\ecen

Listings \ref{l2}-\ref{l5} discuss the dimension independent (function) files. 
We use the \oop\ setup, and thus we refer to \cite{actut} for a general 
introduction concerning the meaning of the stiffness matrix $K$, the 
mass matrix $M$ and the BC matrices $Q$ (and $\Gbc$, not used here), 
and the initialization methods {\tt stanpdeo*D}, $*=1,2,3$, setting 
up an \oop\ object which contains the geometry and FEM space, and the methods 
to assemble the system matrices. 

\hulst{caption={{\small {\tt \dname/cGLinit.m}, which collects some typical initialization commands. {\tt p=stanparam(p)} in line 2 sets the \pdep\ controls, switches and numerical constants to standard 
values; these can always be overwritten afterwards, and some typically are. In line 3 set the function handles to the rhs and its Jacobian, and similarly in 
line 4 
we (re)set the output function handle to {\tt hobra}, which can be used as 
 a standard output when Hopf bifurcations are expected. In lines 5-14, depending 
on the spatial dimension, we create a 1D, 2D or 3D  \oop\ objects,  
essentially consisting of the domain, the FEM setup and the boundary condition. 
In 1D, this is the interval (-lx,lx) with mesh width lx/nx, and 
homogeneous Neumann BC. In line 15 we finish this by preparing the 
associated BC matrices, and afterwards we put this {\tt PDE}--object, 
the number of grid points, and the associated norm weight $\xi$ into {\tt p}. 
Calling {\tt setfemops} in line 17 then immediately refers to {\tt oosetfemops}, 
see Listing \ref{l2}. 
In line 18 
we initialize the solution vector (here with the explicitly known trivial solution 
$u=0$ and append the parameters. 
In the remainder of cGLinit we set some additional controls, mostly explained by the comments. We only remark that p.sw.bifcheck=2 in line 23 tells \pdep\ 
to use algorithm HD2 \mbox{\cite[\S2.1]{hotheo}} to detect bifurcations, 
by computing $\neig$=p.nc.neig=20 eigenvalues near $0$. This is a suitable 
choice since $\pa_u G(0)$ has no real eigenvalues. p.nc.mu1 in line 24 
refers to $\mu_1$ from \mbox{\cite[Remark 2.2]{hotheo}}. 
Finally, {\tt p.vol} in lines 6,8 and 10 is used in the norm \reff{defnorm}, and {\tt p.x0i} is a point index 
for plotting the time-series $t\mapsto u(t,x_{\tt p.x0i})$. }},
label=l1,language=matlab,stepnumber=5, firstnumber=1}{\hdhome/cgl/cGLinit.m}

\hulst{caption={{\small {\tt \dname/sG.m}. Given K and M from 
oosetfemops, we only need to compute the 'nonlinearity' $f$, which we 
outsource to {\tt nodalf}, see Listing \ref{l4}, and then compute 
the rhs $G$ (or residual $r$) as $G(u)=Ku-Mf$. Note that the BC are already 
included in p.mat.K via line 7 of {\tt oosetfemops}. }},
label=l3,language=matlab,stepnumber=5, firstnumber=1}{\hdhome/cgl/sG.m}

\hulst{caption={{\small {\tt \dname/nodalf.m}. The 'nonlinearity' 
(which includes linear terms, i.e., everything except the diffusion terms) 
$f$ from \reff{cAC}. We extract the two components $u_1$ and $u_2$, 
and the parameters from u, introduce an auxiliary variable ua$=|u|^2$, 
and write down the components of $f$ in a standard \mlab\ way. }},
label=l4,language=matlab,stepnumber=5, firstnumber=1}{\hdhome/cgl/nodalf.m}

\hulst{caption={{\small {\tt \dname/sGjac.m}. Similar to 
{\tt sG}, the main problem specific part is $\pa_u f$, put 
into {\tt njac}, the implementation of which follows immediately from {\tt nodalf}. Gu in line 5 is then rather generic. }},
label=l5,language=matlab,stepnumber=5, firstnumber=1}{\hdhome/cgl/sGjac.m}

\hulst{caption={{\small {\tt \dname/oosetfemops.m}. This sets 
the stiffness matrix K, the mass matrix M, and the BC matrix Q for 
\reff{cAC}; see \mbox{\cite[\S1]{actut}} for the meaning of these matrices. 
K$\in\R^{n_p\times n_p}$ in line 2 is the 
'scalar' (i.e., one component) Neumann Laplacian, while M is the scalar 
mass matrix. Similarly, Q$\in\R^{n_p\times n_p}$ in line 3 is a BC matrix 
for one component, and its content depends on the BC set in lines 7, 9 and 11  
of {\tt cGLinit}. In line 6 of oosetfemops we create a zero matrix for 
convenience, and in lines 5,6 we then set up the FEM matrices for the {\em system}  
\reff{cAC}. Here both diagonal blocks of p.mat.K are equal, because so are 
the diffusion constants for both components in \reff{cAC} and the 
BC we consider. However, this setup is quite flexible to implement 
also more complicated differential operators, including off-diagonal blocks 
('cross diffusion'), first order differential operators, and different 
BC in different components. 
}},
label=l2,language=matlab,stepnumber=5, firstnumber=1}{\hdhome/cgl/oosetfemops.m}


\subsection{1D}
\noi
In 1D we use Neumann BC, and $n_x=31$ spatial, and (without mesh-refinement) 
$m=21$ temporal discretization points. 
Listings \ref{l7} and {\ref{l7b} shows the main script file {\tt cmds1d.m} for 1D 
computations (with some omissions wrt to plotting), while Fig.~\ref{f1} 
shows some output generated by running {\tt cmds1d}. The 
 norm in (a) is 
\huga{\label{defnorm}
\|u\|_*:=\|u\|_{L^2(\Om\times (0,T), \R^N)}/\sqrt{T|\Om|}, 
}
which is our default for plotting of Hopf branches. 
During the continuation the default plotting routine {\tt hoplot} also plots the time--series 
$t\mapsto u_1(x_0,t), u_2(x_0,t)$ for some mesh point $x_0$, selected 
by the index {\tt p.hopf.x0i}, which is set in {\tt cGLinit} (see 
also Fig.~\ref{f1}(b))). 
The simulations run in less than 10 seconds 
per branch, but the rather coarse meshes 
lead to some inaccuracies. For instance, the first three HBPs, 
which analytically are at $r=0, 1/4, 1$, are obtained at 
$r=6*10^{-5}, 0.2503, 1.0033$, and (b) also shows some 
visible errors in the period $T$. 
However, these numerical errors quickly decay if we increase $n_x$ and $m$, 
and runtimes stay small. 

\begin{figure}[H]
\bce{\small 
\begin{tabular}{lp{0.73\textwidth}}
(a) BD, norm $\|u(\cdot,\cdot;r)\|_*$&(b) Example plots\\
\ig[width=0.24\textwidth, height=49mm]{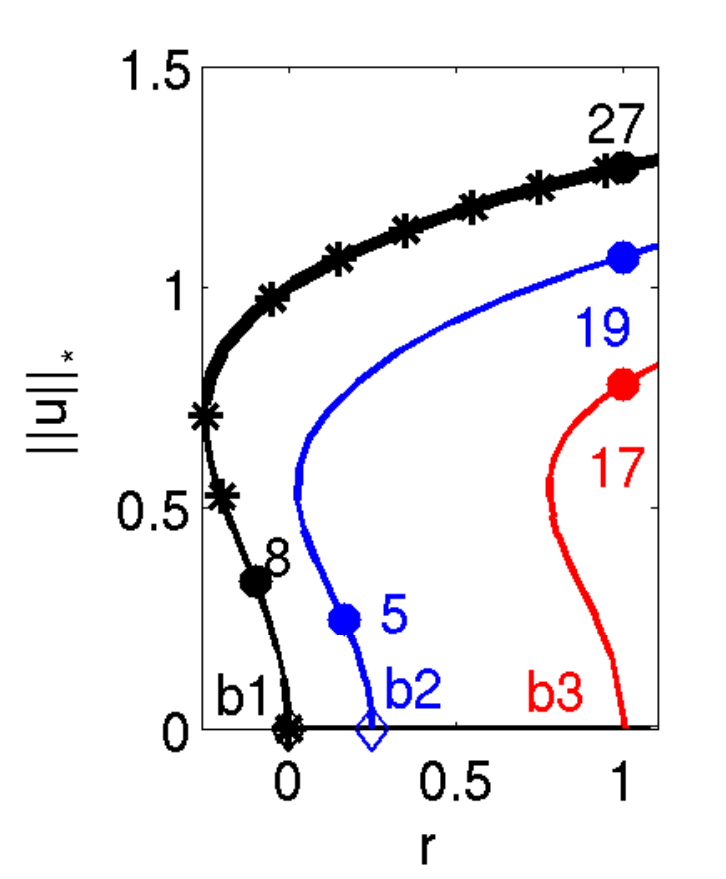}&
\hs{-4mm}
\ig[width=0.25\textwidth,height=49mm]{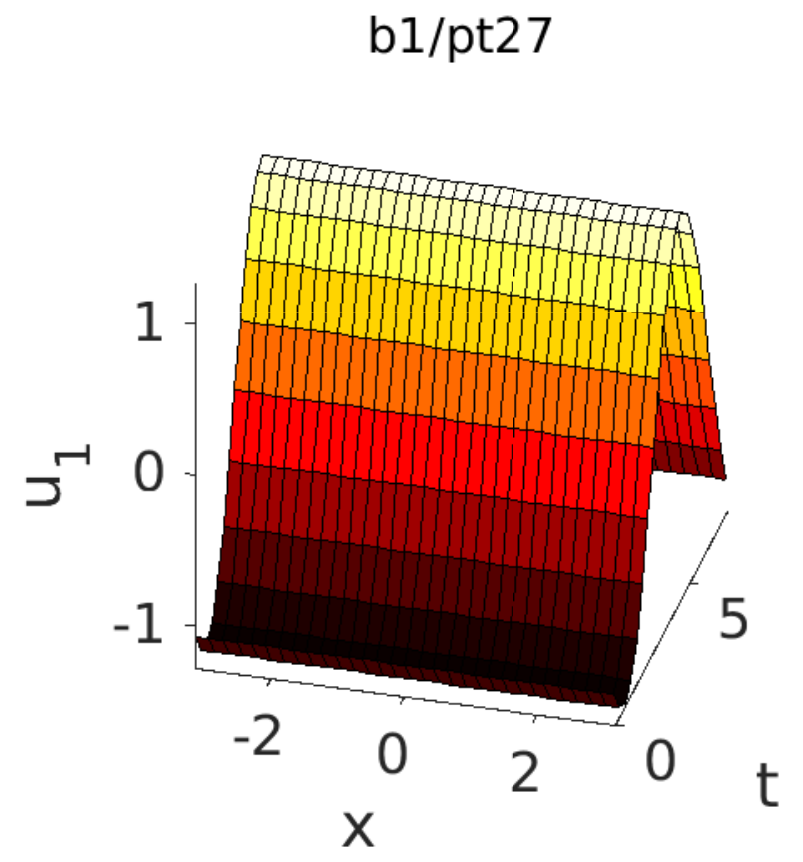}
\ig[width=0.21\textwidth]{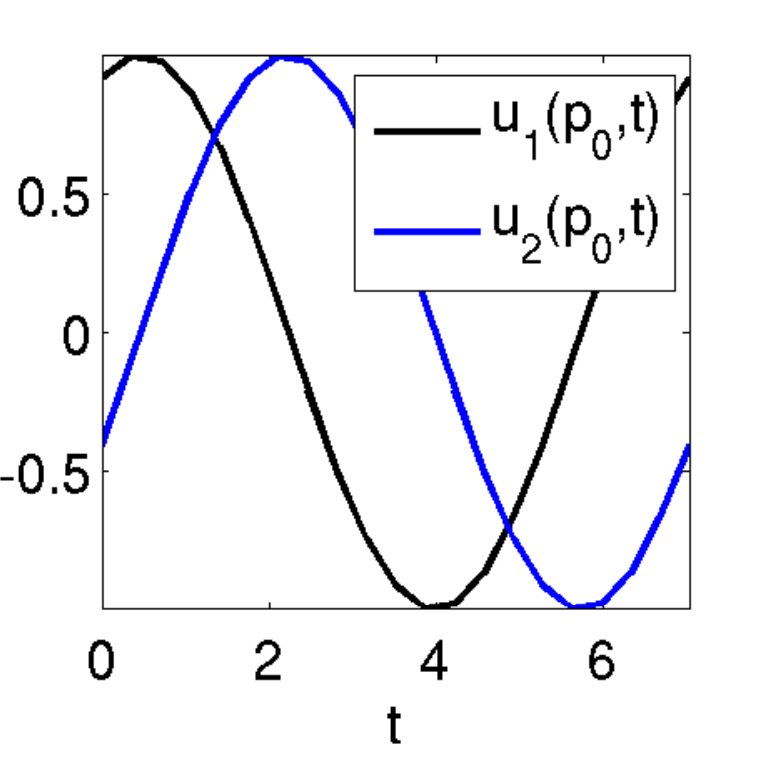}
\ig[width=0.25\textwidth,height=49mm ]{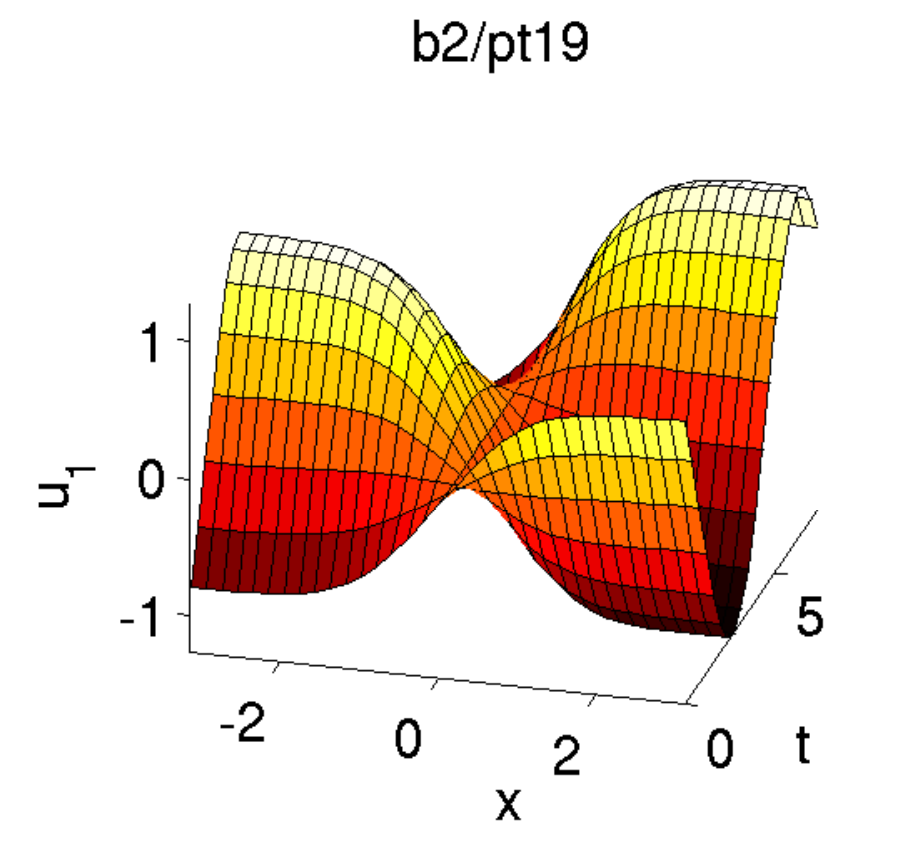}
\end{tabular}
\begin{tabular}{p{0.55\textwidth}p{0.4\textwidth}}
(c) Multipliers at b1/pt8 ($\ind=1$) and b2/pt5 ($\ind=3$) (left),  
and at b1/pt27 ($\ind=0$) (right)
&{\small (d) left: BD, period $T(r)$. Right: numerical periods (for $m=20,40,60$) and analytical period (black dots) on the 1st branch} \\[1mm]
\raisebox{15mm}{\begin{tabular}{l}
\ig[height=30mm]{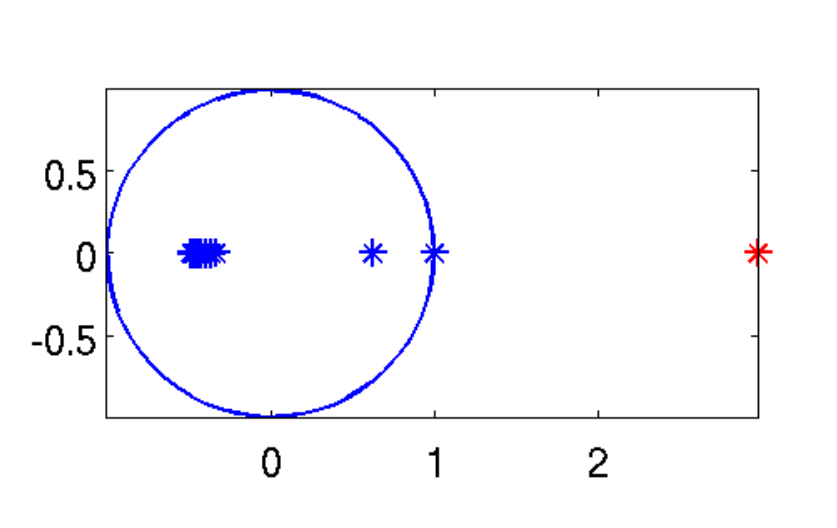}\\[-7mm]
\raisebox{10mm}{\ig[height=25mm]{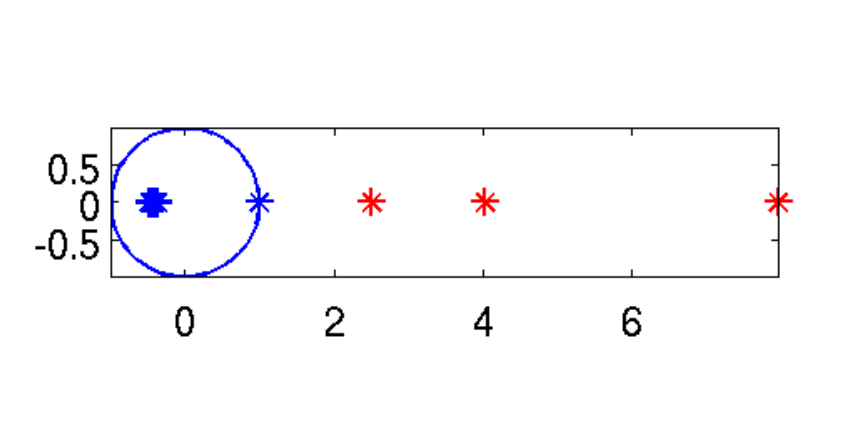}}
\end{tabular}}
\raisebox{10mm}{\ig[height=30mm]{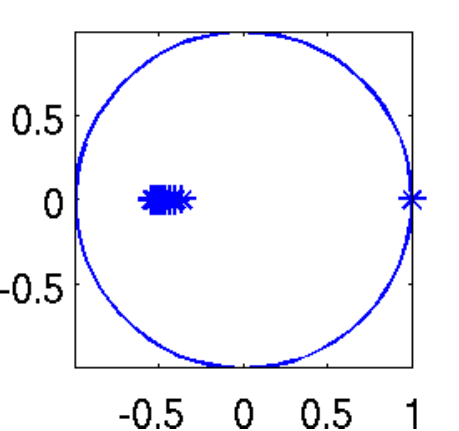}}
&\ig[width=35mm,height=49mm]{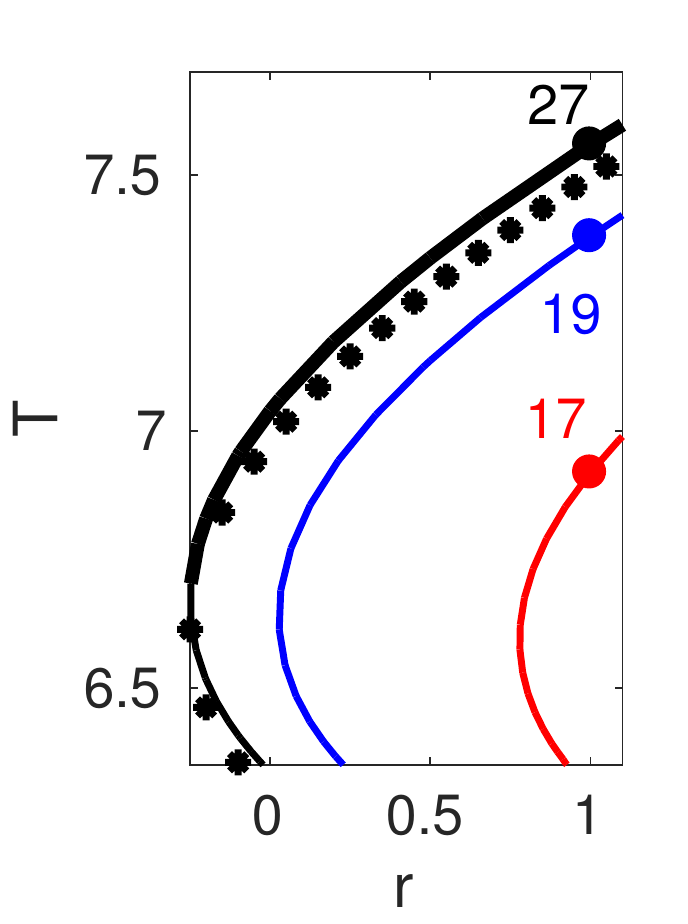}\ig[width=35mm,height=49mm]{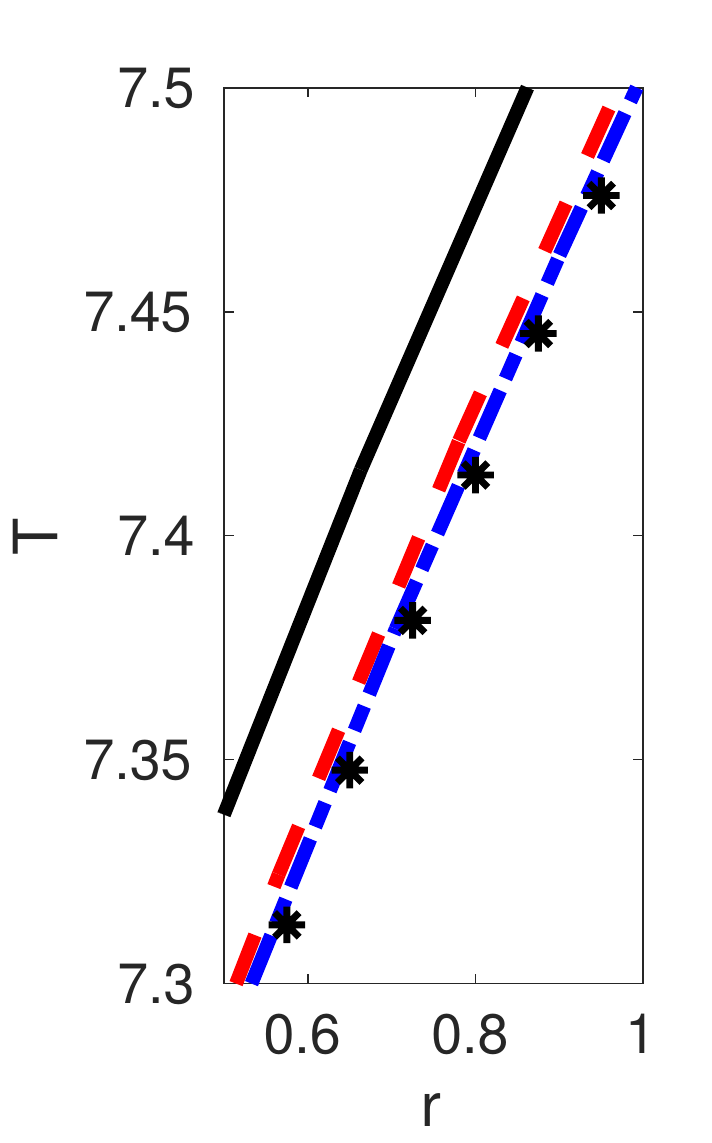}
\end{tabular}
}
\ece 

\vs{-18mm}
   \caption{{\small Selected outputs from {\tt cmds1d.m}, i.e., numerical bifurcation diagrams, example plots 
and (leading 20) Floquet multipliers 
for \reff{cAC} on the domain 
$\Om=(-\pi,\pi)$ with Neumann BC, 30 grid--points in $x$.  
Parameters $(\nu,\mu,c_3,c_5)=(1,0.1,-1,1)$, 
hence bifurcations at (restricting to the first three branches) 
$r=0$ ($k=0$, spatially homogeneous branch, black), 
$r=1/4$ ($k=1/2$, blue) and $r=1$ ($k=1$, red). 
The thick part of the black line in (a) 
indicates the only stable periodic solutions. 
 The black dots in (a) and (d) are from the (spatially homogeneous) analytical solution, see \cite{hotheo}.
For $m=20$ there is a visible error in $T$. 
The right panel of (d) shows the numerical $T$ for different $m$ ($m=20$ black, 
$m=40$ red-dashed, $m=60$ blue-dotted), which 
illustrates the convergence of the numerical solution towards the analytical solution \reff{cACs}. Similarly, the periods also converge on the other 
branches (see {\tt cmds1dconv.m}). The second plot in (b) 
shows a time series at the point {\tt p.hopf.x0i} from {\tt b1/pt27}. 
See also Fig.~\ref{f1b}(b) for a plot of {\tt b3/pt17}. 
  \label{f1}}}
\end{figure}

As usual we recommend 
to run {\tt cmds1d} cell-by-cell to see the effect(s) of the individual cells. 

\hulst{caption={{\small {\tt \dname/cmds1d.m} (first four cells). In cell 1 we initialize the problem and continue the trivial 
branch (with standard settings) to find the HBPs. In cell 2 we then 
compute the first 2 bifurcating Hopf branches in the arclength setting. 
See Appendix \ref{appd} for comments on {\tt hoswibra}, which 
sets all the data structures for periodic orbit continuation and 
of course an initial guess, and thus is the main routine here. 
In line 15 we switch on the Floquet computation with {\tt floq}, 
see \S\ref{flsec}.  In line 16
we set parameters for the (optional) bordered elimination linear system solver 
{\tt lssbel}, see \mbox{\cite[Remark 2.3]{hotheo}}, and in line 17 for the
 preconditioned ilupack solver 
\cite{ilupack} {\tt lssAMG} as an inner solver for {\tt lssbel}. 
This is optional, and controlled by the switch {\tt AMG} in line 18. 
See lines 19,20 for the convenience functions to switch on these solvers and 
to set parameters. 
For the present 1D problem, both {\tt lssAMG} or just {\tt lss} are roughly equally fast, 
but for larger scale problems {\tt lssAMG} is significantly faster. 
In any case, without {\tt ilupack}, {\tt lssbel} gives a significant speedup over {\tt lss} for bordered systems, see also \cite{lsstut} for a tutorial on these solvers. In cell 3 
we do the initial steps for the third Hopf branch  in natural parametrization, 
which gives a refinement of the $t$-mesh by TOM from $m=21$ to $m=41$ 
(here uniform due to the harmonic nature of the time-dependence). 
We then switch to arclength and proceed as before. 
}},
label=l7,language=matlab,stepnumber=5, firstnumber=1, linerange=2-32}{\hdhome/cgl/cmds1d.m}

\hulst{caption={{\small {\tt \dname/cmds1d.m} continued, 
to illustrate (with some omissions) the plot of bifurcation diagrams and solutions. 
Since in {\tt cGLinit} we set p.fuha.outfu=@hobra, i.e., to the standard 
Hopf branch output, and since we have 5 parameters in the problem, 
the period $T$ is at (user-)component 6 of the branch, then follow min and max, and  component 9  contains 
the $L^2$ norm; see also \cite{ooschnaktut} 
for details on the organization of the branch data and on {\tt plotbra}. 
 }},
label=l7b,language=matlab,stepnumber=5, firstnumber=1, linerange=34-37}{\hdhome/cgl/cGL1dcmds_red.m}

Switching to continuation in another parameter works just as for 
stationary problems by calling {\tt p=hoswiparf(\ldots)}.  
See Cells 1 and 2 of {\tt cgl/auxcmds1.m} for an 
example, and Fig.~\ref{f1b}(a) for illustration. Cells 3 and 4 
 of {\tt auxcmds1.m} 
then contain examples for mesh-refinement in $t$, for which there 
are essentially two options. The first is to use {\tt p.sw.para=3} and the mesh-adaption of TOM, the second is {\tt hopftref}, see Listing \ref{l9}.

  \hulst{language=matlab,stepnumber=5, firstnumber=1, linerange=1-5}{\hdhome/cgl/auxcmds1.m}

\hulst{caption={{\small {\tt \dname/auxcmds1.m}. Cells 1 and 2 
illustrate switching to another continuation parameter, while cells 3 and 4 
give simple examples of mesh-adaption in $t$. In cell 3 we use the  error estimator build into TOM. In cell 4 we use {\tt hopftref}, which is a purely ad hoc 
refinement, and which requires a time $t^*$  where to refine from the user. 
In some cases, the convenience function {\tt hogradinf(p)}, 
which inter alia returns the time $t^*$ where $\|\pa_t u(\cdot,t)\|_{\infty}$ 
is maximal, is useful, 
though not in this problem since the solutions considered here are rather 
time harmonic. Note that in contrast to cell 3, or to the routine 
{\tt meshada} for spatial mesh refinement, 
neither {\tt hogradinf(p)} nor {\tt hopftref} 
deal with error estimates in any sense.}},
label=l9,language=matlab,stepnumber=5, linerange=8-20}{\hdhome/cgl/auxcmds1.m}

\begin{figure}[ht]
\bce{\small
\begin{tabular}{lll}
(a) BD $T(c_5)$&(b) $u_1$ at b3/pt17& (c) $u_1$ from (b) after calling 
{\tt hopftref}.\\
\ig[width=0.3\tew,height=40mm]{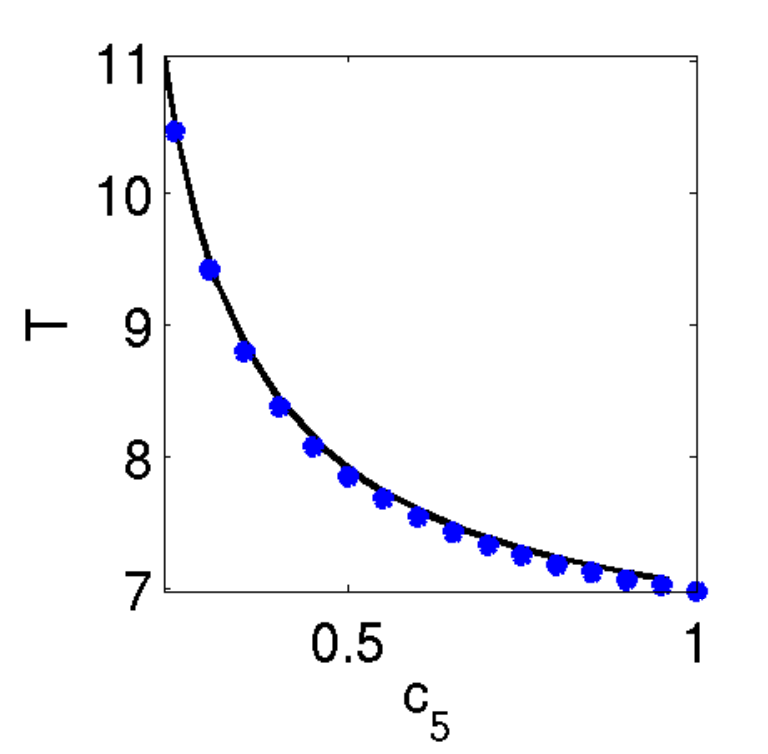}&
\ig[width=0.25\tew]{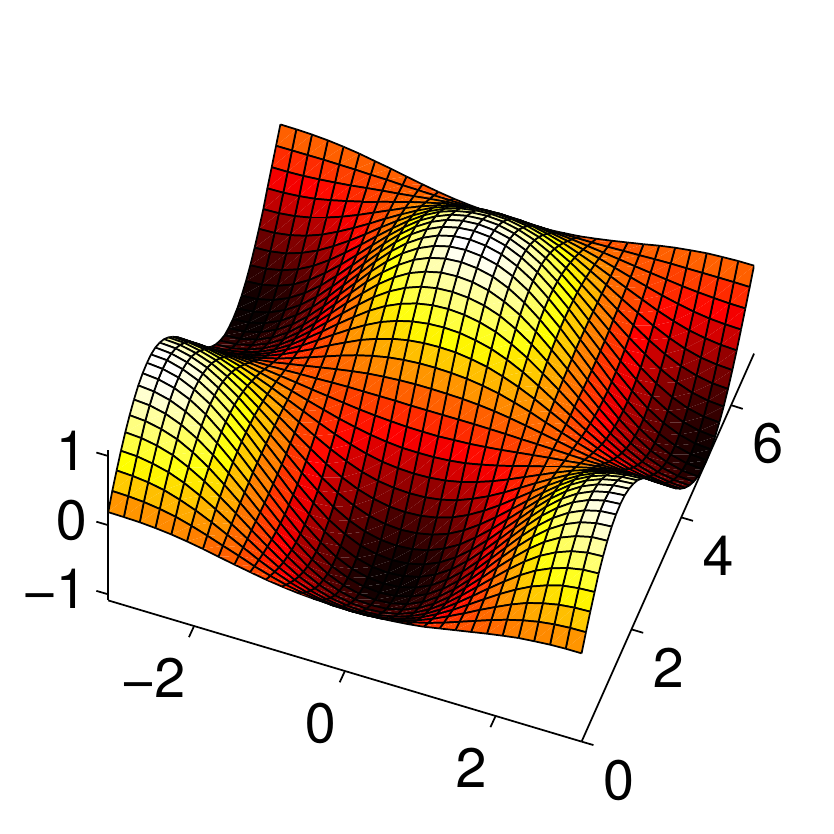}&
\ig[width=0.25\tew]{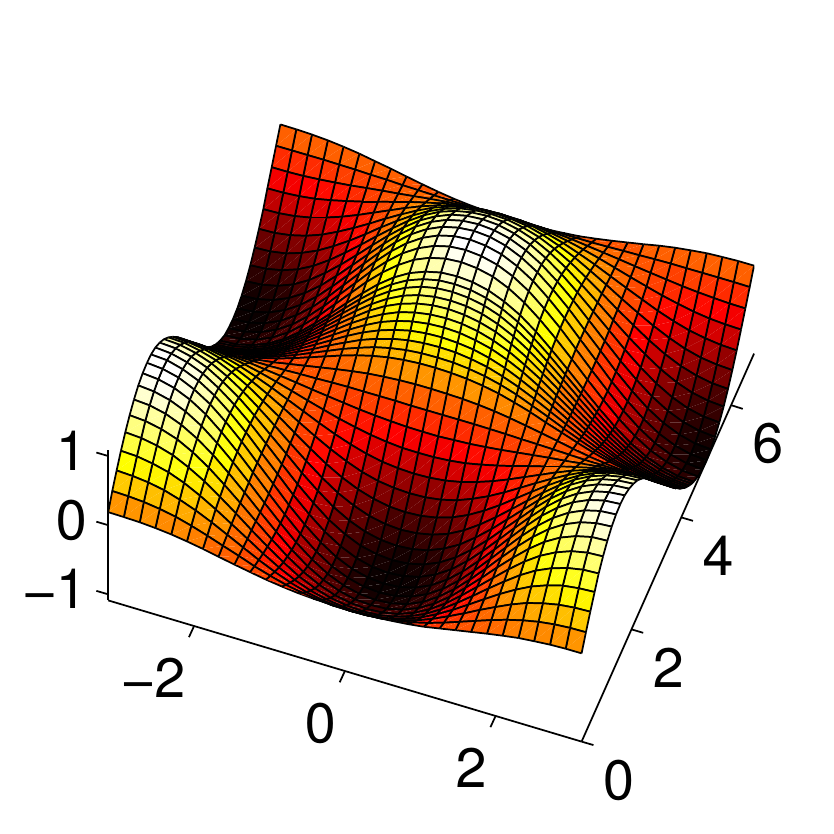}
\end{tabular}}
\ece 

\vs{-5mm}
   \caption{{\small Example outputs from {\tt auxcmds1.m}. (a) Continuing the solution b1/pt28 from Fig.~\ref{f1}(a,b)  
in $c_5$, with comparison to the analytical formula \cite[\S3.1]{hotheo}. 
(b), (c) Solution at {\tt b3/pt17} before and after mesh-refinement in $t$ 
via {\tt hopftref}, here near $t^*=4$. 
  \label{f1b}}}
\end{figure}

\subsection{Remarks on Floquet multipliers and time integration} 
\label{flsec}
For the Floquet multipliers $\ga_j$, $j=1,\ldots,n_u$ ($n_u=Nn_p$ with 
$n_p$ the number of spatial discretization points, see \reff{tformd}) we recall from \cite[\S2.4]{hotheo} 
that we have two algorithms for their computation:%
\footnote{in the software we typically call the Floquet multipliers $\mu$ instead of $\ga$}
\bci
\item \fla\ (encoded 
in the function {\tt floq}) computes 
$0\le {\tt p.hopf.nfloq}\le n_u$ multipliers as 
eigenvalues of the monodromy matrix $\CM$. 
\item
\flb (encoded in {\tt floqps}) 
uses a periodic Schur decomposition of the matrices building $\CM$ 
to compute all $n_u$  multipliers.
\eci 
\flb\ is generally much more accurate and robust, but may be slow.%
\footnote{ For {\tt floqps} one needs to {\tt mex} percomplex.f(F) in the directory 
{\tt pqzschur}, see the README file there.} See also line 15 in Listing 
\ref{l7}. 
There always is the trivial Floquet multiplier $\ga_1=1$ associated to  translational in $t$, and we use $\emu:=|\ga_1-1|$ with the numerical 
$\ga_1$ as a measure for the accuracy of the multiplier computation. 
Furthermore we define the index of a periodic orbit $u_H$ as 
\huga{\label{inddef} 
\ind(u_H)=\text{number 
of multipliers $\ga$ with $|\ga|>1$ (numerically: $|\ga|>1+$p.hopf.fltol)},  
}
such that $\ind(u_H)>0$ indicates instability.

On {\tt b1} in Fig.~\ref{f1}, initially there is 
one unstable multiplier $\ga_2$, i.e., $\ind(u_H)=1$,  
which passes through 1 to enter the unit 
circle at the fold. On b2 we start with $\ind(u_H)=3$, and 
have $\ind(u_H)=2$ after the fold. Near $r=0.45$ another multiplier 
moves through 1 into 
the unit circle, such that afterwards we have $\ind(u_H)=1$, with, 
for instance $\ga_2\approx 167$ at $r=1$. Thus, we may expect 
a  bifurcation near $r=0.45$, 
and similarly we can identify 
a number of possible bifurcation on b3 and other branches. 
The trivial multiplier $\ga_1$ is $10^{-12}$ 
close to $1$ in all these computations, using {\tt floq}. 

In {\tt cgl/auxcmds2.m} we revisit these multiplier computations, and 
complement them with time-integration. For the latter, the idea 
is to start time integration from some point on the periodic orbit, 
e.g.~$u_0(\cdot)=u_H(\cdot,0)$, and to monitor, inter-alia, 
$e(t):=\|u(t,\cdot)-u_0(\cdot)\|$, 
where by default $\|\cdot\|=\|\cdot\|_\infty$. 
Without approximation error for the computation of $u_H$ (including 
the period $T$) and of $t\mapsto u(\cdot,t)$ we would have $e(nT)=0$. In general, even if 
$u_H$ is stable we cannot expect that, in particular due to 
errors in $T$ which will accumulate with $n$,  but nevertheless 
we usually can detect instability of $u_H$ if at some $t$ there is 
a qualitative change in the time--series of $e(t)$.\footnote{\label{fntint}
The time 
integration {\tt hotintxs}  takes inter alia the number {\tt npp} of 
time steps per period $T$ as argument. 
Time integration is much faster than the BVP solver used 
to compute the periodic orbits, and thus {\tt npp} can be chosen 
significantly larger than the number $m$ of time-discretization 
points in the BVP solver. Thus, choosing ${\tt npp}=5m$ or 
${\tt npp}=10m$ appears a reasonable practice.} 
In Fig.~\ref{fstab1}(a), where we use the smaller amplitude periodic 
solution at $r=0$ for the IC, this happens right from the start. 
Panel (b) illustrates 
the stability of the larger amplitude periodic solution at $r=0$, 
while in (c) the instability of the solution on {\tt h2} at $r=1$ 
manifests around $t=30$, with subsequent convergence 
to the (stable) spatially homogeneous periodic orbit

\hulst{caption={{\small {\tt \dname/auxcmds2.m}. Cells 1-3 deal 
with Floquet computations as indicated in the comments. Cells 4-6 deal 
with time integration. In line 10 we set up $u(\cdot,t_0)$ from the Hopf-orbit 
in {\tt 1db1/pt8} as an initial condition, and set some parameters. This is used in Cell 5 for time integration via {\tt hotintxs} (see source for documentation), and Cell 6 plots the results, see Fig.~\ref{fstab1}.}},
label=l8,language=matlab,stepnumber=5, firstnumber=1, lastline=44}{\hdhome/cgl/auxcmds2.m}

\begin{figure}[ht]
{\small 
\bce
\begin{tabular}{lll}
(a) Time series and solution for IC b1/pt8 &(b) IC b1/pt27&(c) IC b2/pt19\\
\hs{-4mm}\ig[width=0.2\tew]{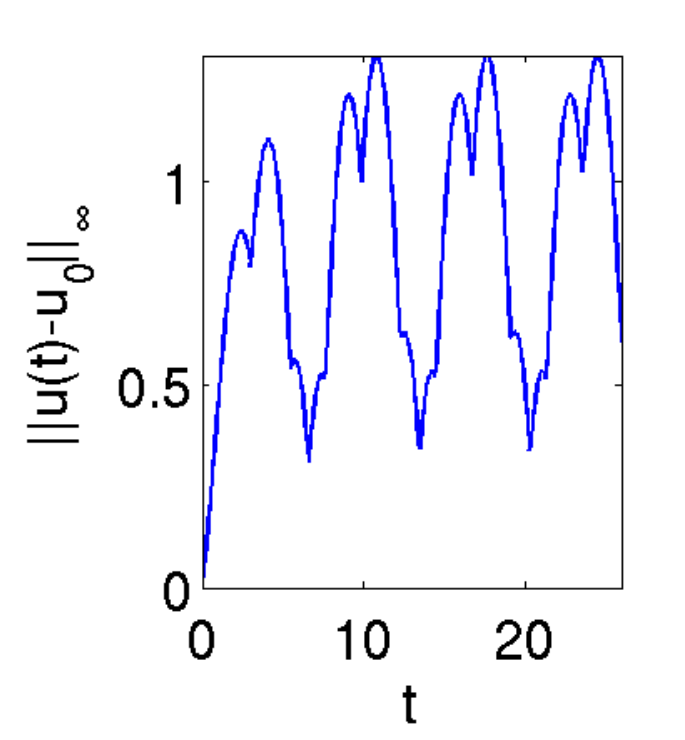}\ig[width=0.2\tew]{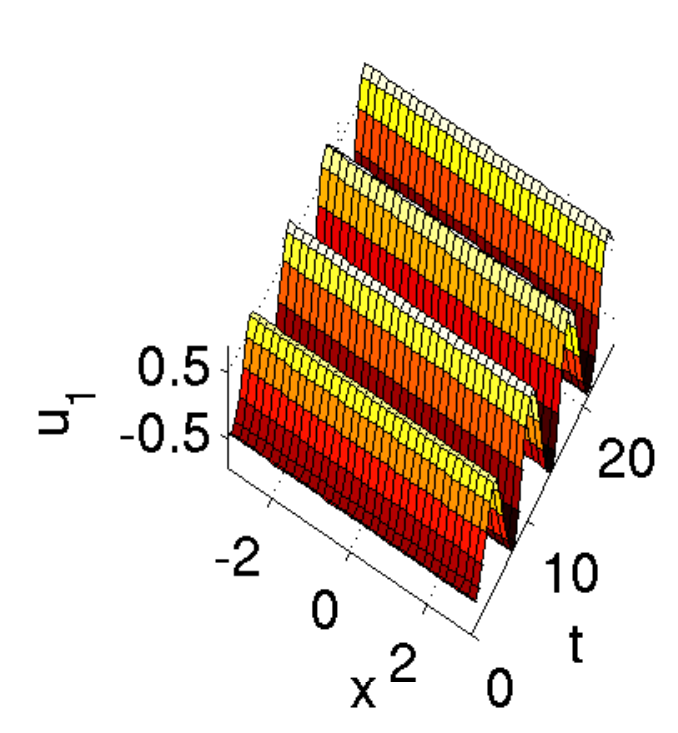}&
\hs{-4mm}\ig[width=0.2\tew]{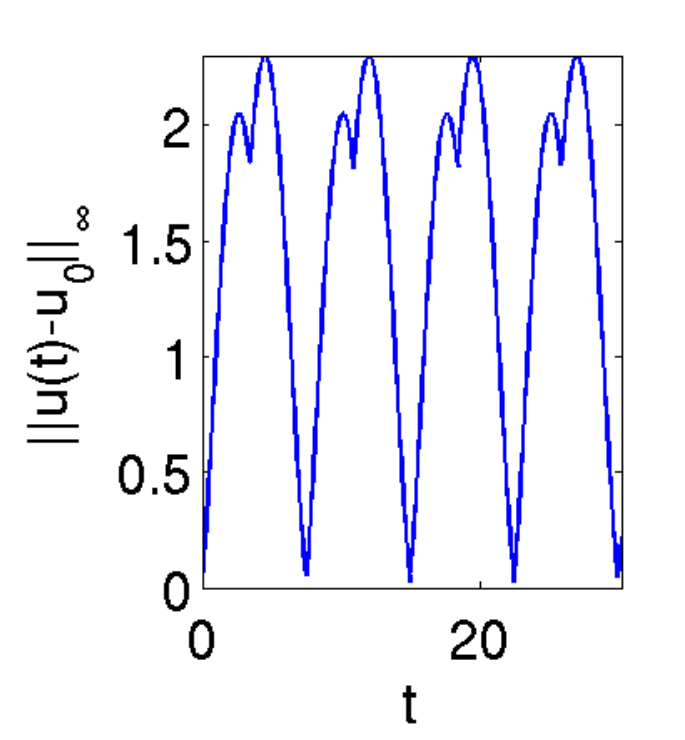}&
\hs{-4mm}\ig[width=0.2\tew]{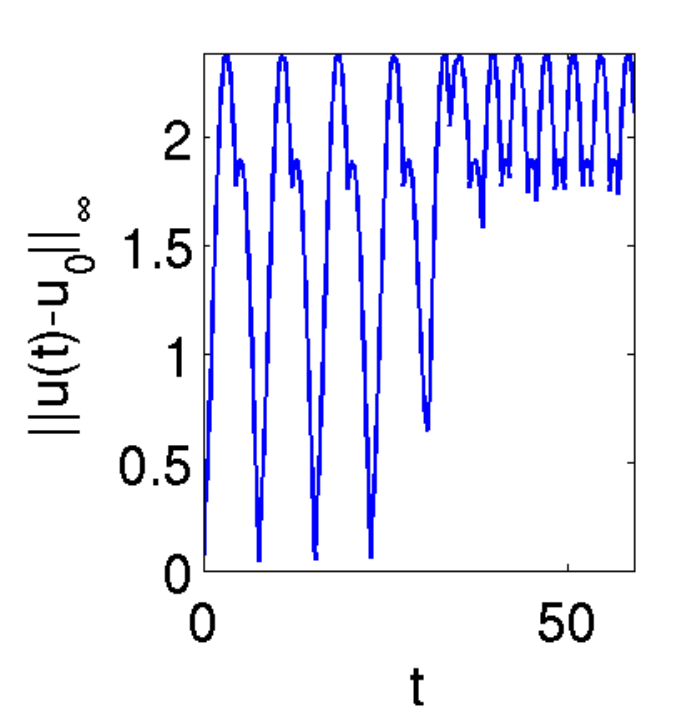}\hs{-3mm}\ig[width=0.21\tew]{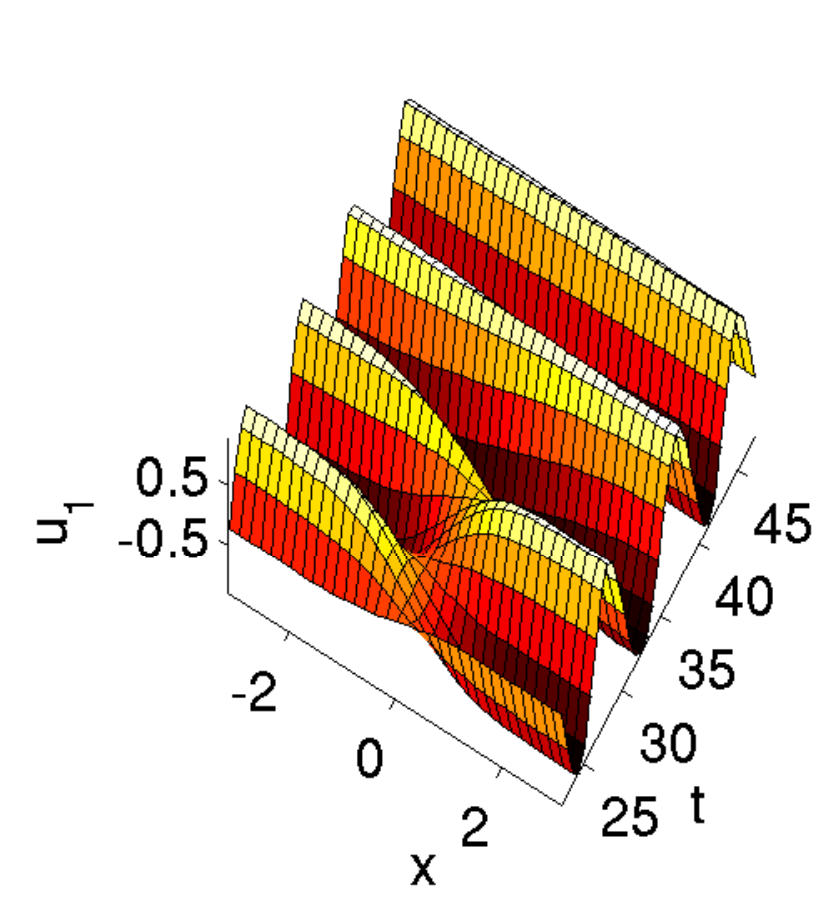}
\end{tabular}
\ece
}
\vs{-5mm}
   \caption{{\small Selected output from {\tt auxcmds2.m}, i.e., stability experiments for \reff{cAC} in 1D. 
(a) IC h1/pt8, 
time series of $\|u(\cdot,t)-u_0\|_\infty$ 
and $u_1(x,t)$, showing the convergence to the larger amplitude 
solution at the same $r$. (b) IC h1/pt27 from Fig.~\ref{f1}, where we plot 
 $\|u(\cdot,t)-u_0\|_\infty$ for $t\in [0,4T]$, which shows stability of 
the periodic orbit, and a good agreement for the temporal period 
under time integration. (c) instability of b2/pt19 from Fig.~\ref{f1}, 
and again convergence to the solution on the b1 branch. 
Note that the time--stepping is much 
finer than the appearance of the solution plots, but we only save 
the solution (and hence plot) every 100th step, cf.~footnote \ref{fntint}. 
  \label{fstab1}}}
\end{figure}

\subsection{2D}  
In 2D we choose homogeneous Dirichlet BC for $u_1,u_2$, see lines 8,9 
in {\tt cGLinit}, and {\tt oosetfemops.m}. 
Then the first two HBPs are at 
$r_1=5/4$ ($k=(1/2,1)$, and $r_2=2$ ($k=(1,1)$). 
The script file {\tt cmds2d.m} follows the same principles as 
{\tt cmds1d.m}, and includes some time integration as well, 
and in the last cell an example 
for creating a movie of a periodic orbits. 

Figure \ref{f2} shows some results from {\tt cmds2d.m}, obtained on a coarse 
mesh of $41\times 21$ points, hence $n_u=1722$ spatial unknowns, 
yielding the numerical values $r_1=1.2526$ 
and $r_2=2.01$. With $m=20$ temporal discretization points, the 
computation of each Hopf branch then takes about a minute. 
Again, the numerical HBPs converge to the exact values when decreasing 
the mesh width, but at the prize of longer computations for the Hopf branches. 
For the Floquet multipliers we obtain a similar picture as in 1D. 
The first branch has $\ind(u_H)=1$ up to the fold, and $\ind(u_H)=0$ 
afterwards, and on b2 $\ind(u_H)$ decreases from 3 to 2 at the fold 
and to 1 near $r=7.2$. Panel 
(c) illustrates the 2D analogue of Fig.~\ref{fstab1}(c), i.e., 
the instability of the second Hopf branch and stability of the first.

\begin{figure}[ht]
{\small 
\bce
\begin{tabular}{lll}
(a) BD &(b) solution snapshots&(c) Instability of b2/pt10, conv.~to b1 \\
\hs{-2mm}\ig[width=0.17\tew]{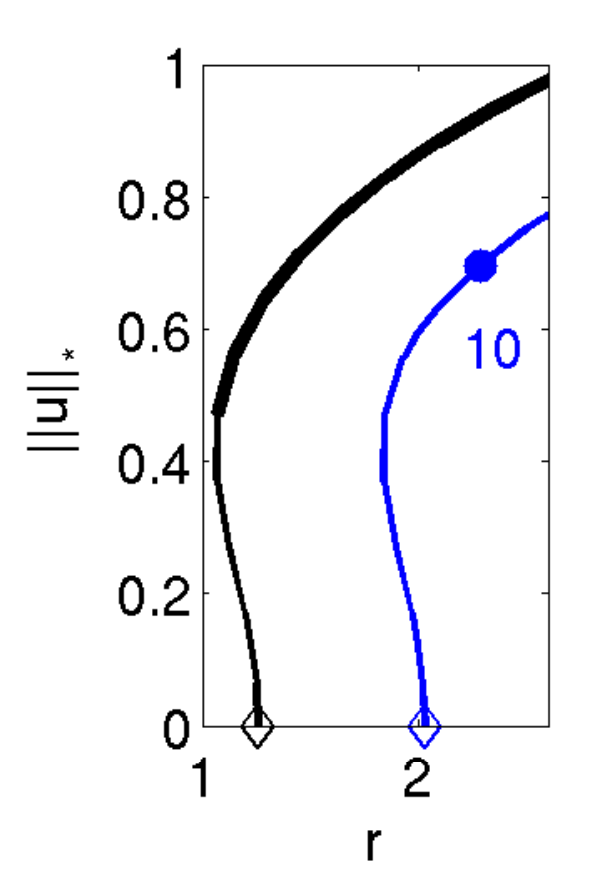}&\hs{-4mm}
\ig[width=0.25\tew]{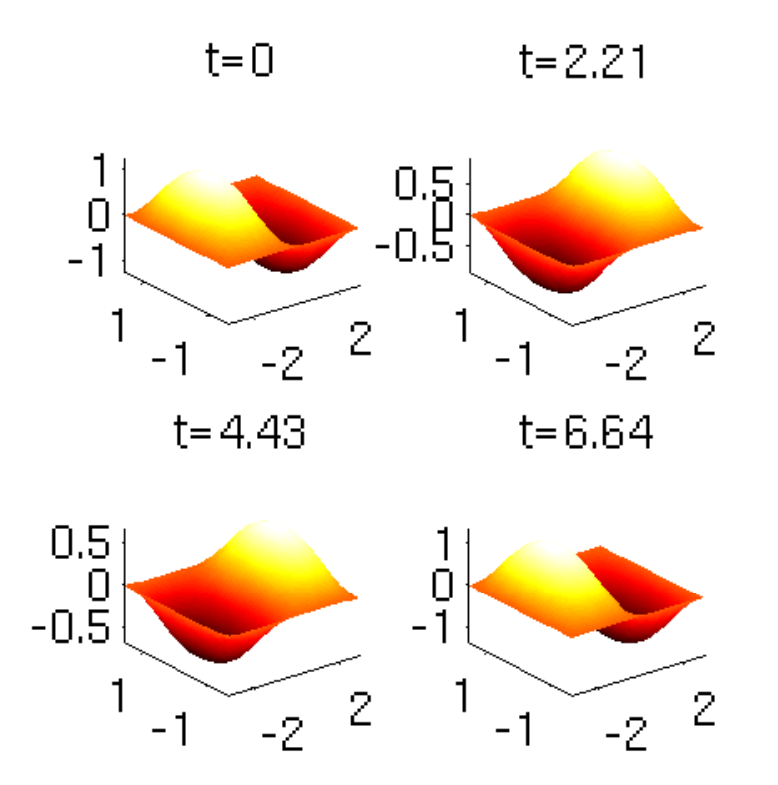}&
\hs{-4mm}\ig[width=0.17\tew]{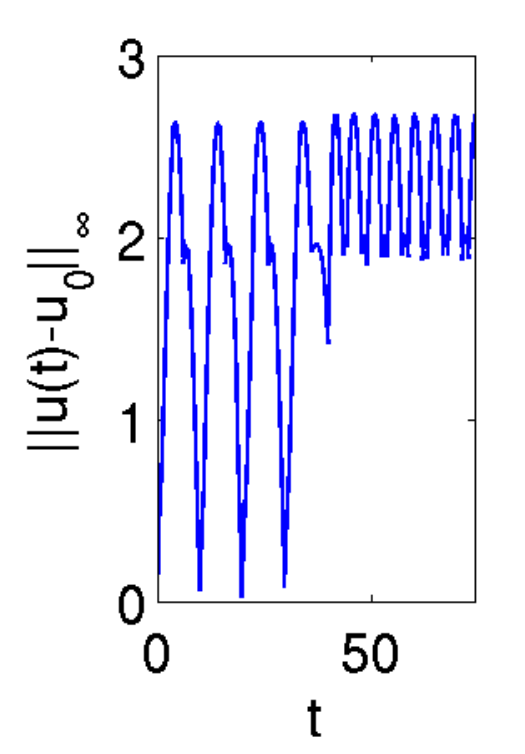}
\ig[width=0.4\tew]{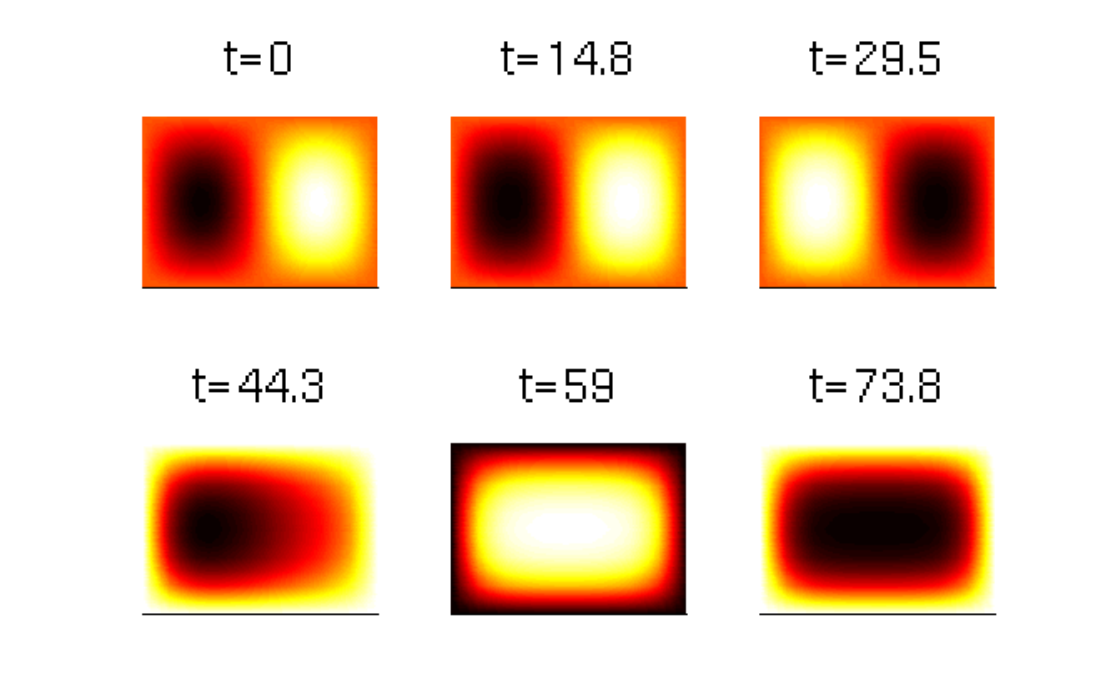}
\end{tabular}
\ece
}

\vs{-5mm}
   \caption{{\small Example plots from {\tt cmds2d.m}. (a) Bifurcation diagrams of the first 2 Hopf branches for 
 \reff{cAC} in 2D. (b) Solution snapshot from b2/pt10, at 
$t=0, \frac 3 {10} T, \frac 6 {10} T, \frac 9 {10} T$. (c) 
Time integration starting from (b) ($t=0$), with 
convergence to the first Hopf branch.  }
  \label{f2}}
\end{figure}

\subsection{3D} 
To illustrate that exactly the same setup also works in 3D, 
in {\tt cmds3d.m} and Fig.~\ref{f4} we consider \reff{cAC} 
over $\Om=(-\pi,\pi)\times (-\pi/2,\pi/2)\times(-\pi/4,\pi/4)$. Here we 
use a {\em very} coarse tetrahedral mesh of $n_p{=}2912$ points, 
thus $5824$ DoF in space. Analytically, the first 2 HBPs are 
$r_1{=}21/4$ ($k=(1/2,1,2)$) and $r_2{=}6$ ($k=(1,1,2)$, but with the 
coarse mesh we numerically obtain $r_0{=}5.47$ and $r_1{=}6.29$. 
Again, this can be greatly improved by, e.g., halving the spatial 
mesh width, but then the Hopf branches become very expensive.  
Using $m=20$, the computation of the branches (with 15 
continuation steps each) in Fig.~\ref{f4} takes about 10 minutes, 
and a call of  {\tt floqap} to 
a posteriori compute the Floquet multipliers about 50 seconds. Again, 
on b1, $\ind(u_H){=}1$ up to fold and $\ind(u_H){=}0$ afterwards, 
while on b2 $\ind(u_H)$ decreases from 3 to 2 at the fold 
and to 1 at the end of the branch, 
and time integration from an IC from b2 
yields convergence to a periodic solution from b1. 

The script {\tt cmds2d.m} follows the same principles as the 1D and 
2D scripts. In 3D, the ``slice plot'' 
in Fig.~\ref{f4}(b), indicated by {\tt p.plot.pstyle=1} should be used 
as a default setting, while the isolevels in (c) (via {\tt p.plot.pstyle=2}) 
often require some fine tuning. Additionally we provide a ``face plot'' option 
{\tt p.plot.pstyle=3}, which however is useless for Dirichlet BC. 

\begin{figure}[ht]
{\small 
\bce
\begin{tabular}{l}
(a) BD, $\|u\|_*$ and $T$\hs{5mm} (b) Example slice plot\\
\ig[width=0.22\tew, height=38mm]{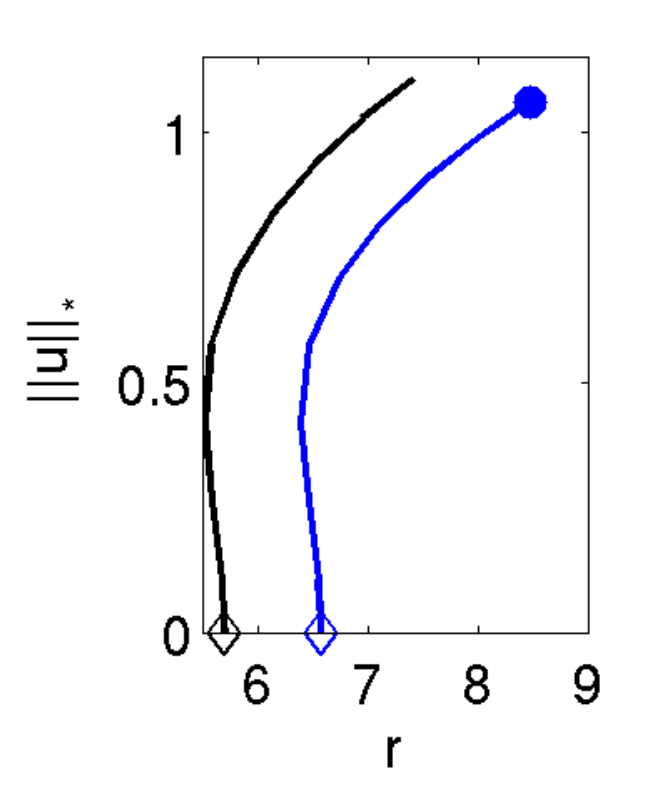}
\hs{-0mm}\ig[width=0.55\tew]{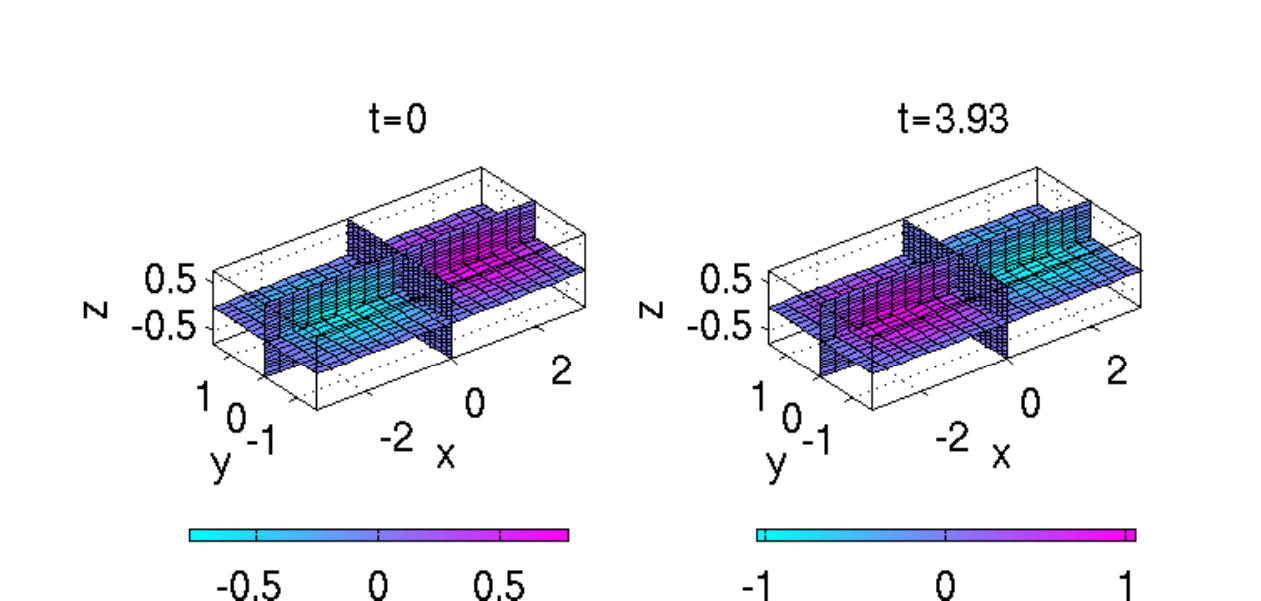}\\
(c) Example isoplot  \\
\hs{-10mm}\ig[width=0.6\tew]{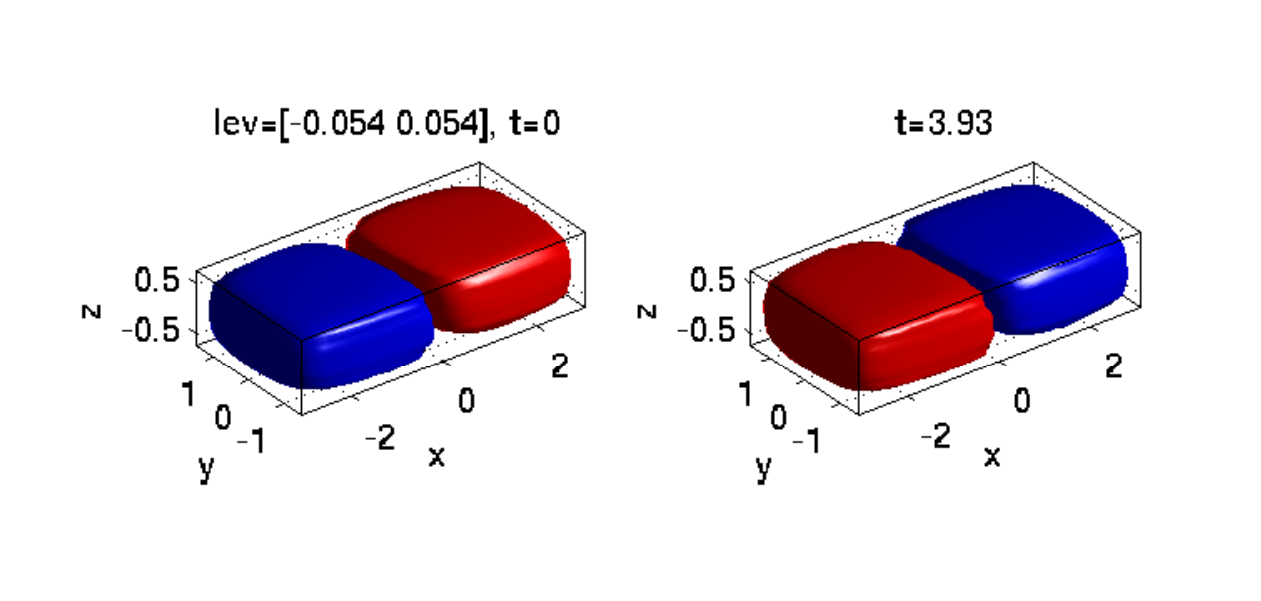}
\end{tabular}
\ece
}
\vs{-14mm}
   \caption{{\small Example plots from {\tt cmds3d.m}. (a) Bifurcation diagram of first 2 Hopf branches for 
 \reff{cAC} in 3D. (b,c) Solution snapshots at $t=0$ and $t=T/2$ for the 
blue dot in (a); 
slice-plot in (b), and isolevel plot in (c) with levels 
$0.525m_1+0.475m_2$ and $0.475m_1+0.525m_2$, 
where $m_1=\min_{x,t} u_1(x,t)$ and $m_2=\max_{x,t} u_1(x,t)$. 
  \label{f4}}}
\end{figure}

\section{An extended Brusselator: Demo {\tt brussel}}\label{brusec} \def\dname{brussel}\def\dhome{./hopfdemos/brussel}
In \cite[\S3.2]{hotheo} we consider an example with an interesting interplay between 
stationary patterns and Hopf bifurcations, and 
where there are typically many eigenvalues with small real parts, 
such detecting 
HBPs with {\tt bifcheck=2} without first using 
{\tt initeig} for setting a guess for a shift $\om_1$ 
is problematic. The model, following \cite{yd02} is an 'extended Brusselator', 
namely the three component reaction--diffusion system 
\huga{\label{ebru}
\pa_t u=D_u\Delta u+f(u,v)-cu+dw,\quad  
\pa_t v=D_v\Delta v+g(u,v), \quad \pa_t w=D_w\Delta w+cu-dw,
}
where $f(u,v)=a-(1+b)u+u^2v$, $g(u,v)=bu-u^2v$, with kinetic parameters 
$a,b,c,d$ and diffusion constants $D_u,D_v,D_w$. 
We consider \reff{ebru} on rectangular domains in 1D and 2D, with 
homogeneous Neumann 
BC for all three components. The system has the trivial spatially homogeneous 
steady state 
$$U_s=(u,v,w):=(a,b/a,ac/d),
$$ 
and in suitable parameter regimes it shows co-dimension 2 points between 
Hopf, Turing--Hopf (aka wave), and (stationary) Turing bifurcations 
from $U_s$. A discussion of these instabilities of $U_s$ in the $a-b$ plane 
is given in \cite{yd02} for fixed parameters 
\huga{\label{brupar}
(c,d,D_u,D_v,D_w)=(1,1,0.01,0.1,1).  
} 
In our simulations we additionally fix $a=0.95$, and take 
$b$ as the primary bifurcation parameter. 

For the quite rich bifurcation results, which include primary spatially homogeneous and patterned 
Hopf bifurcations from $U\equiv U_s$, and Turing bifurcations from $U_s$  followed 
by secondary Hopf bifurcations, we refer to \cite[\S3.2]{hotheo}. 
Regarding the implementation, Table \ref{brutab} lists the scripts and functions in {\tt brussel}. Except for the additional 
component ($N=3$ instead of $N=2$) this is quite similar to {\tt cgl}, 
with one crucial difference, in particular in 2D, on which we focus in 
\S\ref{br2dsec}. First, however, we shall focus on 1D and additional to \cite[\S3.2]{hotheo} compute bifurcation lines in the $a$--$b$ parameter 
plane by branch point continuation and Hopf point continuation, 
and compute secondary bifurcations {\em from} Hopf orbits.

\taskip
\begin{table}[ht]\caption{Scripts and functions in {\tt hopfdemos/brussel}. 
\label{brutab}}
{\small 
\begin{tabular}{l|p{0.77\tew}}
script/function&purpose,remarks\\
\hline
bru1dcmds&basic BDs for 1D, including some time integration. \\
bru1dcmds\_b&extension of bru1dcmds, dealing with bifurcations from the 
primary Hopf orbit\\
bru2dcmds&script for 2D, including preparatory step {\tt initeig} for guessing 
$\ri\om$ for Hopf bifurcations, and some time integration\\
cmdsHPc&script for Hopf and branch point continuation to compute Fig.~\ref{brfig1}(a)\\
auxcmds1&1D auxiliaries, illustrating spatial mesh refinement on 
Turing branches\\ 
auxcmds2&2D auxiliaries: 
illustration of problems with many small real eigenvalues\\
\hline
e2rsbru&elements to refine selector, interface to \oop's equivalent of 
{\tt pdejmps}\\ 
evalplot&script for plotting eigenvalues for linearization around 
spat.~homogeneous solution, see \cite[Fig.7(b)]{hotheo}. \\
bruinit&initialization as usual\\
oosetfemops&the FEM operator for \reff{ebru}, \oop\ setting\\
sG, sGjac, nodalf&rhs, Jacobian, and nonlinearity, as usual\\
bpjac, hpjac&computing (directional) second derivatives for BP and HP continuation
\end{tabular}
}
\end{table}\teskip

\subsection{1D}
\label{br1dsec}
Listing \ref{brul1} shows the startup in 1D. We fix $a=0.95$ and choose $b$ as 
the continuation parameter, starting at $b=2.75$, over the domain $\Om=(-l_x,l_x)$, 
$l_x=\pi/k_{TH}$, where $k_{TH}=1.4$ is chosen to have the first 
bifurcation from $U_s$ to a Turing-Hopf (or wave) branch. This follows from, 
e.g., \cite{yd02}, see also \cite{hotheo}, 
but also from the bifurcation lines and 
spectral plots in Fig.~\ref{brfig1}(a,b), which we explain below. 
The 'standard' files such as {\tt bruinit.m, oosetfemops.m, sG.m}, and {\tt  sGjac.m} 
are really standard and thus we refer to their sources. Moreover, 
regarding {\tt initeig} in line 5 of {\tt bru1dcmds.m} we refer 
to \S\ref{br2dsec}, where this becomes crucial. The continuation 
in line 6 of {\tt bru1dcmds.m} then yields the first three bifurcations 
as predicted from Fig.~\ref{brfig1}(a,b), and subsequently further 
steady and Hopf bifurcations. See Fig.~\ref{brfig1}(c) for the BD of 
these branches, and \cite[\S3.4]{hotheo} for further discussion. 
Here we first want to explain how Fig.~\ref{brfig1}(a) 
can be computed by Hopf point continuation (HPC) and branch point 
continuation BPC, see Listing \ref{brul2}.  

\hulst{caption={{\small {\tt \dname/bru1dcmds.m} (first 6 lines). 
The {\tt initeig} in line 5 is not strictly necessary in 1D, but useful 
for speed. The remainder of {\tt bru1dcmds} computes a number 
of steady and Hopf bifurcations from {\tt hom1d}, and some secondary 
Hopf bifurcations from Turing branches, and we refer to \cite{hotheo} 
for the associated BDs and solution plots.}},label=brul1,language=matlab,stepnumber=5, linerange=2-7,firstnumber=1}{\dhome/bru1dcmds.m}

\hulst{caption={{\small {\tt \dname/bruHPCcmds.m} (first 15 lines). 
The ideas of HP continuation and of BP continuation are explained in 
\S\ref{hpcsec} and \cite{pftut}, respectively. }},label=brul2,language=matlab,stepnumber=5, linerange=1-15,firstnumber=1}{\dhome/bruHPCcmds.m}

\subsubsection{Hopf point continuation}\label{hpcsec}
Similar to fold continuation and branch point continuation (\cite{p2p2} and \cite[\S3.4]{pftut}), Hopf point continuation (HPC) can be done 
via suitable extended systems. Here we use \cite[\S4.3.2]{govaerts}
\huga{\label{HPC1}
H(U):=\bpm G\\G_u\phi_r+\om M\phi_i\\G_u\phi_i-\om M\phi_r\\
c^T\phi_r-1\\
c^T\phi_i\epm
=0\in\R^{3n_u+2}, \quad U=(u,\phi_r,\phi_i,\om,\lam), 
}
where $\ri\om\in\R$ is the desired eigenvalue of $G_u$, 
$\phi=\phi_r+\ri\phi_i\in\C^{n_u}$ an associated eigenvector, and 
$c_r\in\R^{n_u}$ is a normalization vector. We thus have $3n_u+2$ equations 
for the $3n_u+2$ real unknowns $U$, and 
in \cite[Proposition 4.3.3]{govaerts} it is shown that \reff{HPC1} 
is regular at a simple Hopf bifurcation point. 

Thus, \reff{HPC1} can be used for localization of (simple) Hopf points 
(implemented in the \pdep\ function {\tt hploc}) if a 
sufficiently good initial guess $U$ is given, and, moreover, freeing 
a second parameter $w$ we can use the extended system 
\huga{\label{HPC2}
\CH(U,w)=\bpm H(U,w)\\p(U,w,{\rm ds})\epm=\bpm 0\\0\epm \in\R^{(3n_u+2)+1}
} 
for HPC, where $p(U,w,{\rm ds})$ is the standard arclength condition, 
with suitable weights for the parameters $\lam,w$.  This requires the Jacobian 
\huga{\label{HPC3}
\pa_{U}H(U)=\bpm
G_u&0&0&0&G_\lam\\
\pa_u(G_u\phi_r)&G_u&\om M&M\phi_i&\pa_\lam(G_u\phi_r)\\
\pa_u(G_u\phi_i)&-\om M&G_u&-M\phi_r&\pa_\lam(G_u\phi_i)\\
0&c^T&0&0&0\\
0&0&c^T&0&0\epm\in\R^{(3n_u+2)\times(3n_u+2) }. 
}
Here, $G_u$ is already available, the $\pa_\lam \cdot$ expressions 
are cheap from finite differences, as well as the $w$ derivatives needed 
in the arclength continuation, and expressions such as $\om M$ are 
easy. Thus the main task is to compute the directional 2nd derivatives 
\huga{\label{HPC4}
\bpm \pa_u(G_u\phi_r)\\\pa_u(G_u\phi_r)\epm\in\R^{2n_u\times n_u}. 
}
This can be done numerically, but this may be expensive, and for semilinear 
problems $G(u)=Ku-Mf(u)$ it is typically easy to write a function 
{\tt hpjac} (or with some other problem specific name {\tt name}) which returns \reff{HPC4}, and which should be 
registered as {\tt p.fuha.spjac=@hpjac} (or {\tt p.fuha.spjac=@name}). 
 For instance, for $N=2$ components 
and $\phi_r=(\phi_1,\phi_2)$ 
we have 
\huga{\label{HPC4b}
\pa_u(G_u\phi_r)
=-M\pa_u\bpm f_{1,u_1}\phi_1{+}f_{1,u_2}\phi_2\\
f_{2,u_1}\phi_1{+}f_{2,u_2}\phi_2\epm 
=-M\bpm f_{1,u_1u_1}\phi_1{+}f_{1,u_2u_1}\phi_2&
f_{1,u_1u_2}\phi_1{+}f_{1,u_2u_2}\phi_2\\
f_{2,u_1u_1}\phi_1{+}f_{2,u_2u_1}\phi_2&
f_{2,u_1u_2}\phi_1{+}f_{2,u_2u_2}\phi_2\epm, 
}
and we obtain the same expression for $\pa_u(G_u\phi_i)$ with 
 $\phi_i=(\phi_1,\phi_2)$. Accordingly, {\tt brussel/hpjac.m} 
returns $\bpm \pa_u(G_u\phi_r)\\\pa_u(G_u\phi_r)\epm$ for the three 
component semilinear system \reff{ebru}. For general testing 
we also provide the function {\tt hpjaccheck}, which checks 
{\tt p.fuha.hpjac} against finite differences. 

To initialize HPC, the user can call {\tt p=hpcontini('hom1d','hpt1',1,'hpc1')}, 
see line 2 of Listing \ref{brul2},
where the third argument 
gives the new free parameter $w$. Here $w=a$, which is at position 1 in the 
parameter vector. This triples {\tt p.nu} and sets a number of 
further switches, for instance for automatically taking care of the structure 
of $\pa_UH(U)$ in \reff{HPC3}. For convenience {\tt hpcontini} also 
directly sets {\tt p.fuha.spjac=@hpjac}, which of course  
the user can reset afterwards. 
Then calling {\tt cont} will continue \reff{HPC2} in $w$, 
and thus we produce the 'wave' and 'Hopf' lines in Fig.~\ref{brfig1}(a). 
  Use {\tt hpcontexit}  to return 
to 'normal' continuation (in the original primary parameter). 

Similarly, BPC is based on the extended system 
\cite[\S3.3.2]{mei2000}
\huga{\label{bpe}
H(U)=\bpm G(u,\lam)+\mu M\psi\\
G_u^T(u,w)\psi\\
\|\psi\|_2^2-1\\
\spr{\psi,G_{\lam}(u,w)}\epm 
=0\in\R^{2n_u+2}, \quad U=(u,\psi,w),  
}
where $(u,\lam)$ is a (simple) BP (for the continuation in $\lam$), 
$\psi$ is an adjoint kernel vector, 
$w=(\lam,\mu)$ with $w_1=\lam$ the primary active parameter and 
$w_2=\mu$ as additional active parameter. The BPC requires 
the Jacobian $\pa_U H$ of which $\pa_u(G_u^T \psi)$ is potentially 
difficult to implement. However, again for semilinear problems 
$\pa_u(G_u^T \psi)$ has a similar structure as \reff{HPC4b}, see \cite[\S3.4]{pftut} for further details. In particular, for \reff{ebru} 
it can be implemented rather easily, see {\tt bpjac.m}. 
The actual BPC is then initialized by calling {\tt bpcontini}, see 
line 11 of Listing \ref{brul2}, and the BPC produces the 'Turing line' 
in Fig.~\ref{brfig1}(a). Using {\tt bpcontexit} returns 
to 'normal' continuation. 

\begin{figure}[ht]
\bce 
\begin{tabular}{p{0.33\textwidth}p{0.21\textwidth}p{0.3\textwidth}} 
{\small (a) Bifurcation lines in the $a$--$b$ plane}&{\small (b) spectral plots}
&{\small (c) BD (from \cite[Fig.7]{hotheo}) } \\
\ig[width=0.27\textwidth,height=55mm]{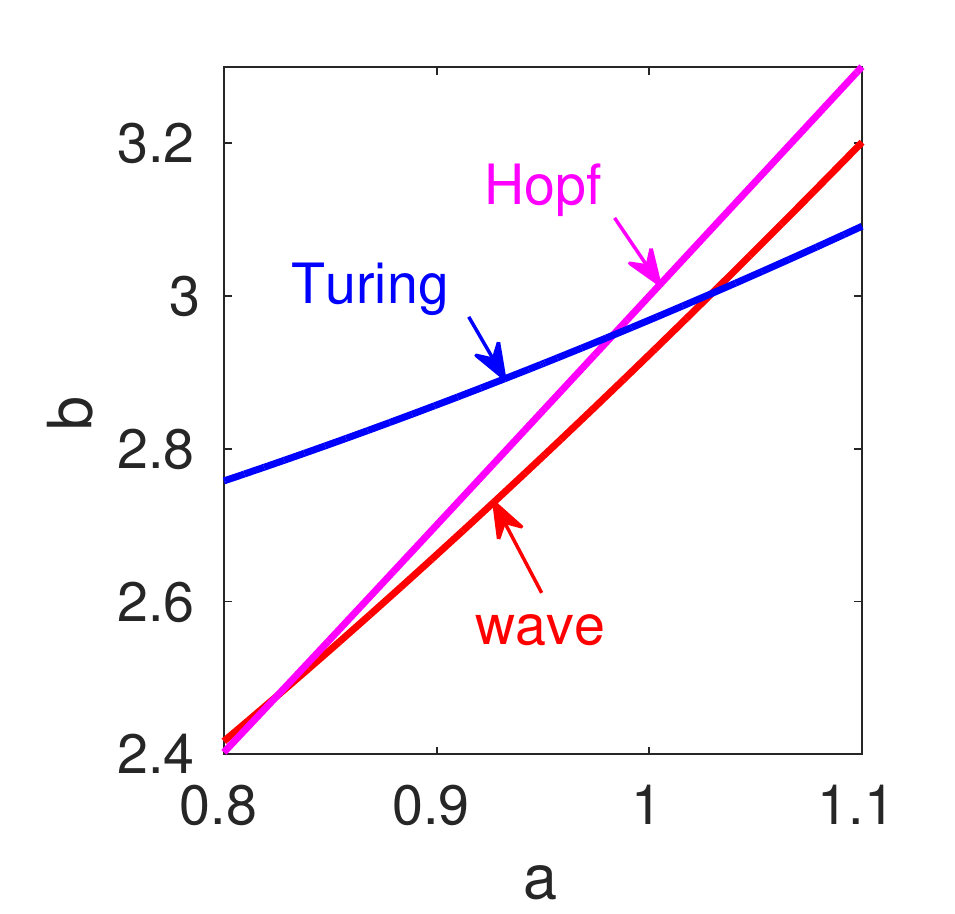}&
\raisebox{28mm}{\begin{tabular}{l}
\hs{-7mm}\ig[width=0.22\textwidth,height=25mm]{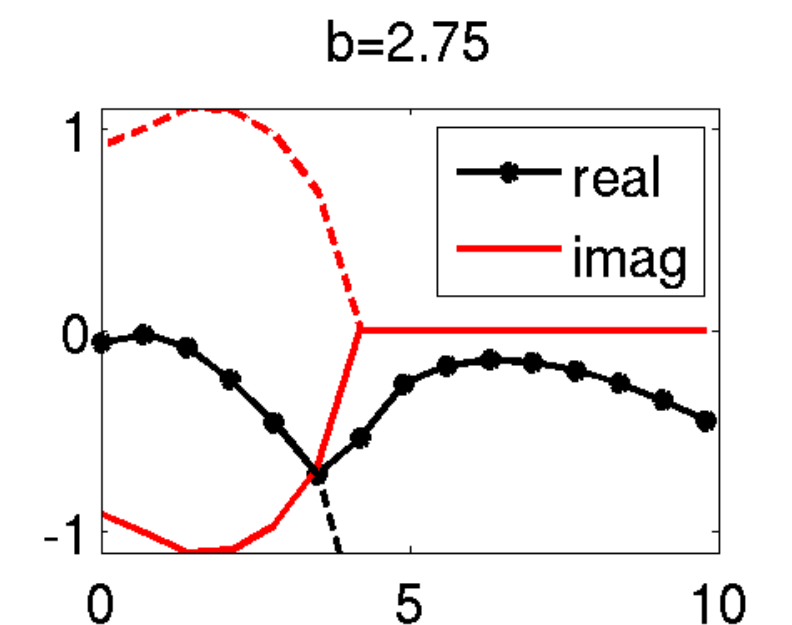}\\
\hs{-7mm}\ig[width=0.22\textwidth, height=25mm]{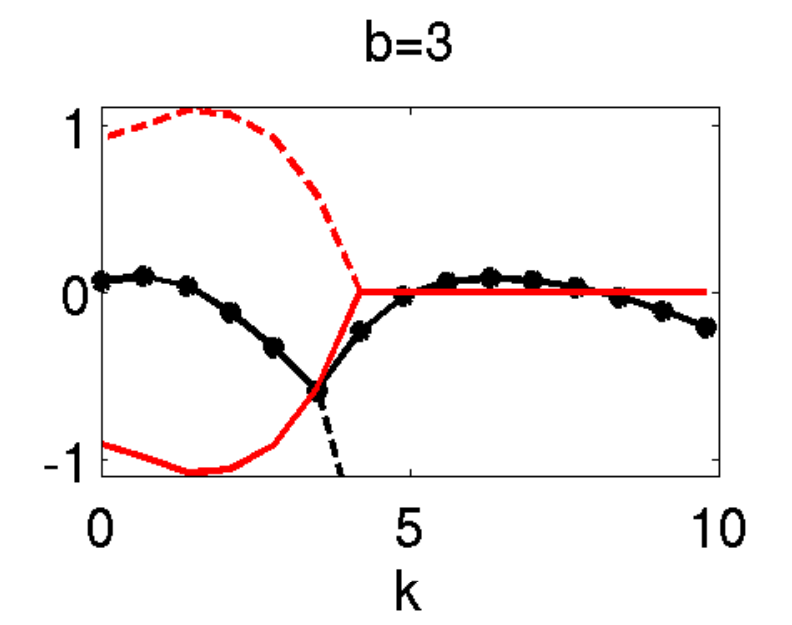}
\end{tabular}}
&\ig[width=0.25\textwidth,height=55mm]{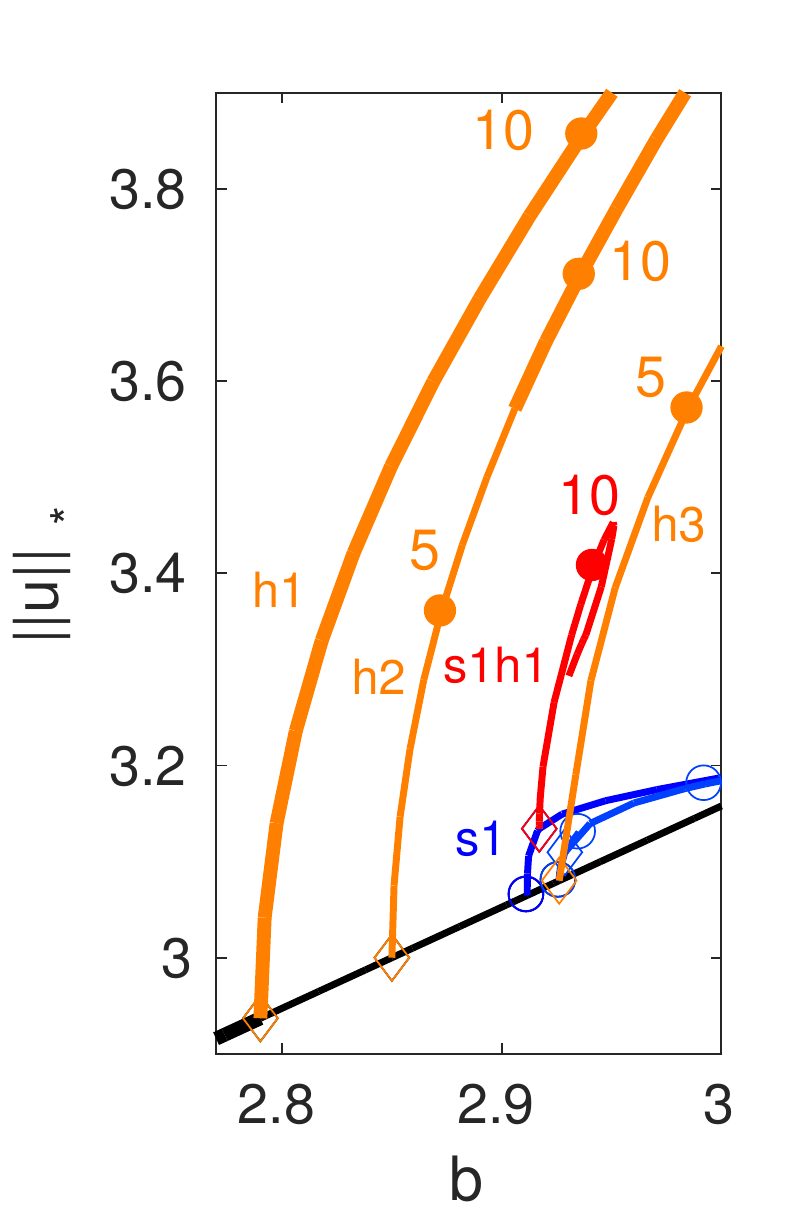}
\end{tabular}
\ece
\vs{-5mm}
   \caption{{\small Results for \reff{ebru}, from {\tt cmdsHPc.m} (a), 
{\tt evalplot.m} (b), and 
{\tt bru1dcmds.m} (c). (a) Bifurcation lines as obtained from 
branch point (Turing--line) and Hopf point (Hopf and wave line) continuation. 
Compare \cite[Fig.7a]{hotheo} or \cite{yd02}. 
(b) Spectrum of the linearization of \reff{ebru} around $U_s$, $a=0.95$ fixed. 
The dots on the real part show the admissible wave numbers on the 
subsequently used 1D domain. 
(c) Bifurcation diagram of Hopf and Turing branches; 
see also \cite[Fig.7]{hotheo} for the associated solution plots. 
  \label{brfig1}}}
\end{figure}

\subsubsection{Bifurcations from the first Hopf branch}\label{pobifsec}
As a second extension of what is presented for \reff{ebru} in \cite[\S3.3]{hotheo} 
we give some results on bifurcations from the first Hopf branch 
{\tt h1} in Fig.~\ref{brfig1}(c), associated  to 
Floquet multipliers going through $1$.%
 (Of course, a critical multiplier also goes through 1 for 
a periodic orbit fold as, e.g., for the cGL in \S\ref{cglsec}, but here we are interested in genuine bifurcations.)
The used method is described in Appendix \ref{appb}, together with 
the case of period doubling bifurcations associated  to 
Floquet multipliers going through $-1$.
\brem\label{bsprem}{\rm  Our methods are somewhat 
preliminary in the following sense: \\
(a) The localization of the branch points uses a simple bisection 
based on the change of ind$(u_H)$, cf.~\reff{inddef}, as a multiplier  
crosses the unit circle. Such bisections work well and robustly 
for bifurcations from steady branches 
(and can always be improved to high accuracy using the above extended 
systems for BPs and HPs), but the bisection for critical multipliers is 
often more difficult (i.e., less accurate) due to many multipliers $\ga_j$ 
close to the unit circle, and also often due to 
sensitive dependence of the $\ga_j$ on 
the continuation parameter $\lam$, and, moreover, on the numerical time 
discretization fineness. Typically, some trial and error is needed here. \\
(b) The computation of predictors for branch switching is currently 
based on the classical monodromy matrix $\CM$, see \reff{fl2} and \reff{fl3}. 
As explained 
in, e.g., \cite{lust01}, see also \cite{hotheo} and \S\ref{ocsec}, this may be 
unstable numerically, in particular for non--dissipative problems. 
The multiplier computations should then be based on 
the periodic QZ-Schur algorithm ({\bf FA2} algorithm in \pdep, see also 
\S\ref{ocsec}), but for 
the branch switching predictor computations this has not been 
implemented yet. 
}\eex\erem  

Despite Remark \ref{bsprem}, for 'nice' problems the branch--switching seems to be robust enough. 
Figure \ref{brfig1} shows some results from {\tt bru1dcmds\_b}, and 
Listing \ref{brul3} shows the basic commands. 
In line 2 we reload a point from the first Hopf branch and set {\tt p.hopf.bisec} 
which determines how many bisections 
are done to localize a BP after the detection of an index change. 
Here (and in other problems) {\tt p.hopf.bisec}=5 seems a reasonable 
value. We then continue further and find the two BPs on {\tt 1dh1} 
indicated in Fig.~\ref{brfig1} where the red and magenta branches 
bifurcate. The branch switching is done in lines 4-8. Typically this 
requires a rather large {\tt ds} (all this of course depends on scaling), 
and often one step with a large residual tolerance is needed to 
get on the bifurcating branch (which indicates that the predictor 
is not very accurate). However, once on the bifurcating branch, 
{\tt p.nc.tol} can and should be decreased again. The same strategy is used for the second (magenta) 
bifurcating branch. 

\begin{figure}[ht]
\bce 
\begin{tabular}{p{0.32\textwidth}p{0.6\textwidth}} 
{\small (a)  BD of secondary bifurcations from the first Hopf branch}&
{\small (b) Solution and Floquet plots at first two BPs and on the red branch} \\
\ig[width=0.3\textwidth,height=60mm]{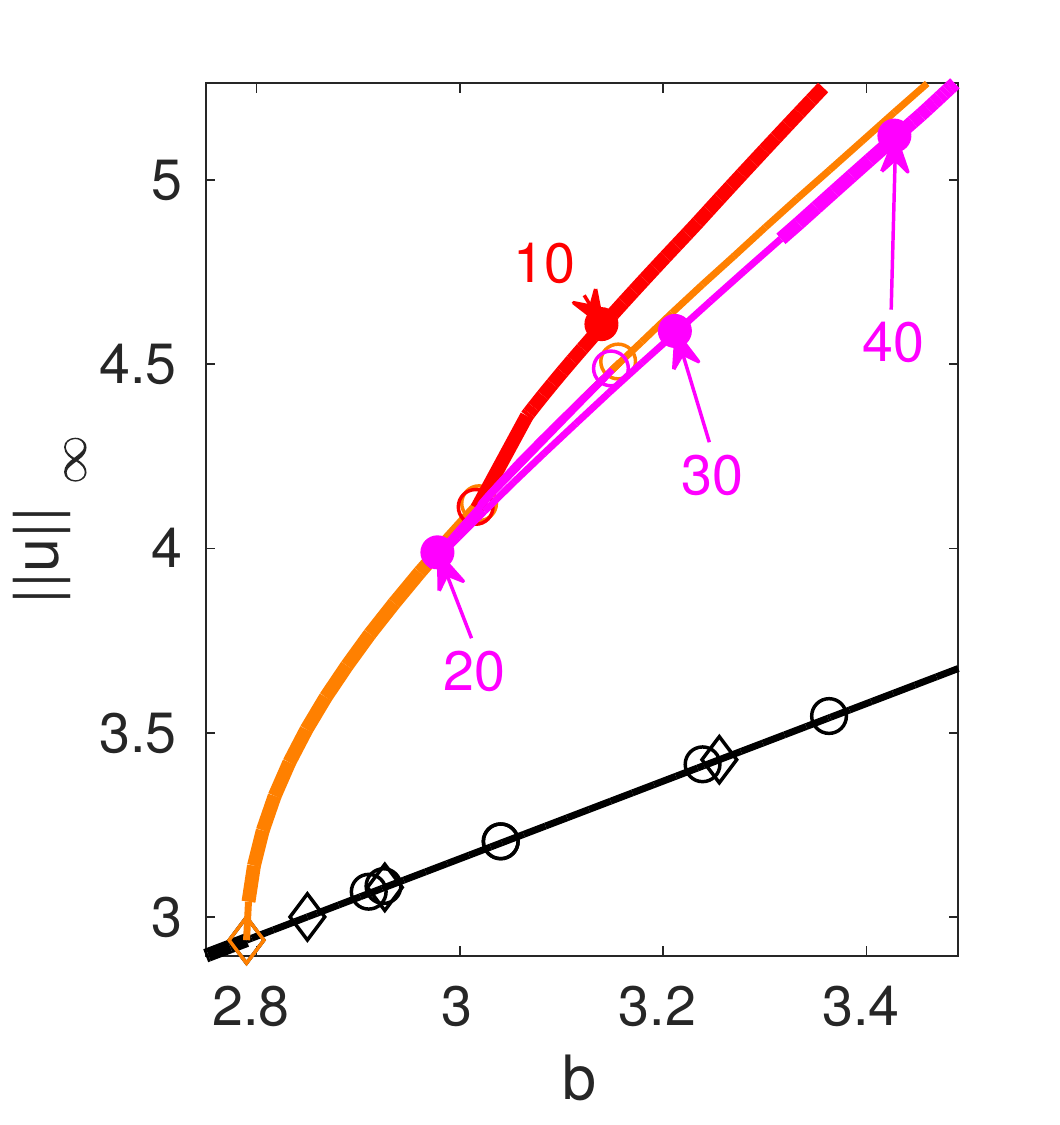}&
\raisebox{28mm}{\begin{tabular}{l}
\ig[width=0.19\textwidth]{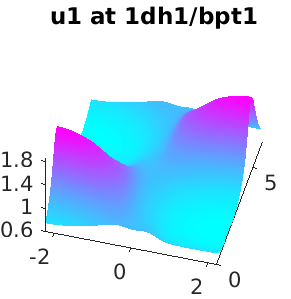}
\ig[width=0.19\textwidth]{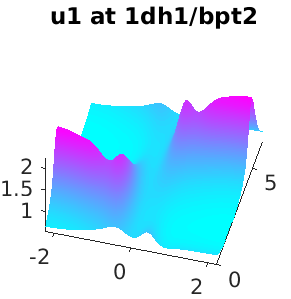}
\ig[width=0.19\textwidth]{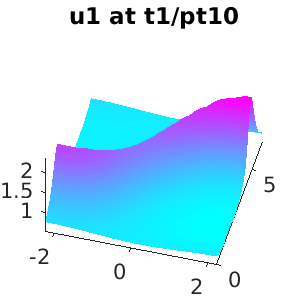}
\\
\ig[width=0.17\textwidth]{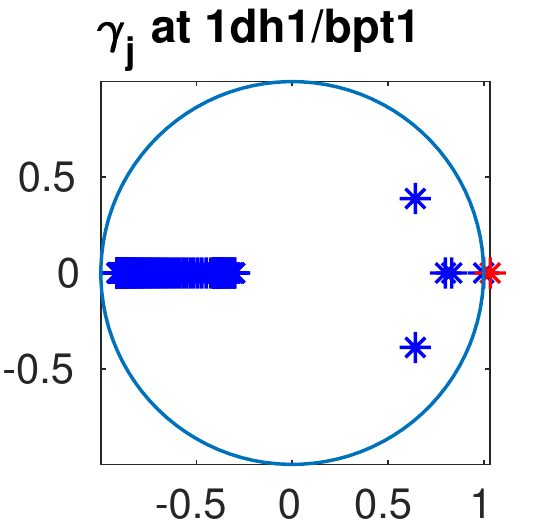}
\ig[width=0.19\textwidth]{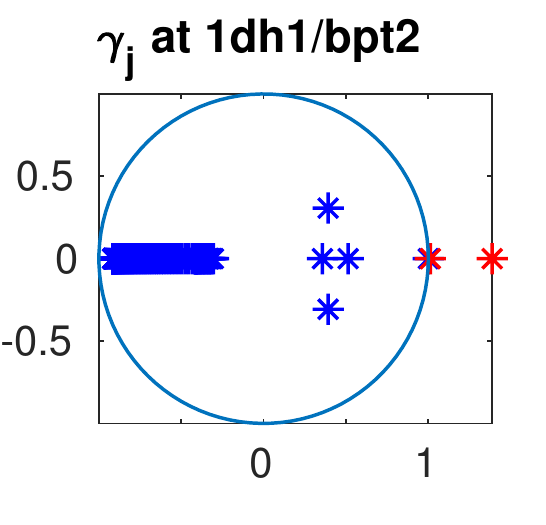}
\ig[width=0.17\textwidth]{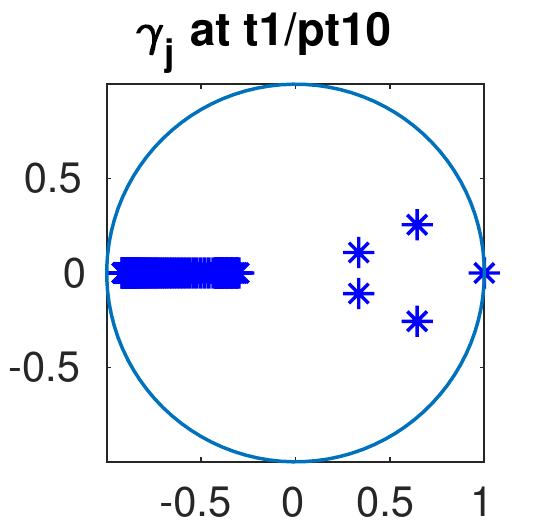}
\end{tabular}}
\end{tabular}\\
{\small (c) solutions and multipliers on the magenta branch}\\[2mm]
\ig[width=0.16\textwidth]{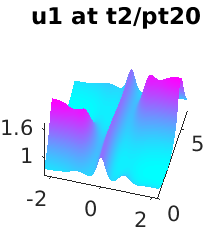}\ig[width=0.16\textwidth]{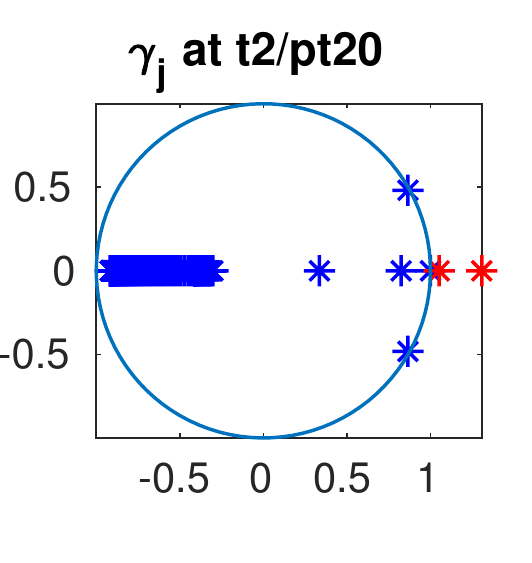}
\ig[width=0.16\textwidth]{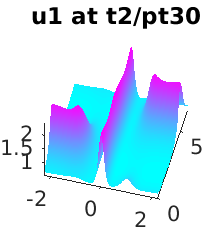}\ig[width=0.16\textwidth]{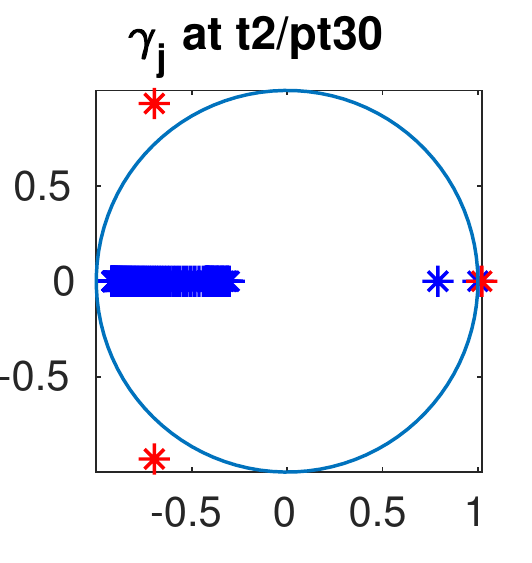}
\ig[width=0.16\textwidth]{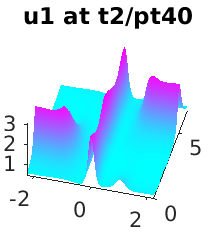}\ig[width=0.16\textwidth]{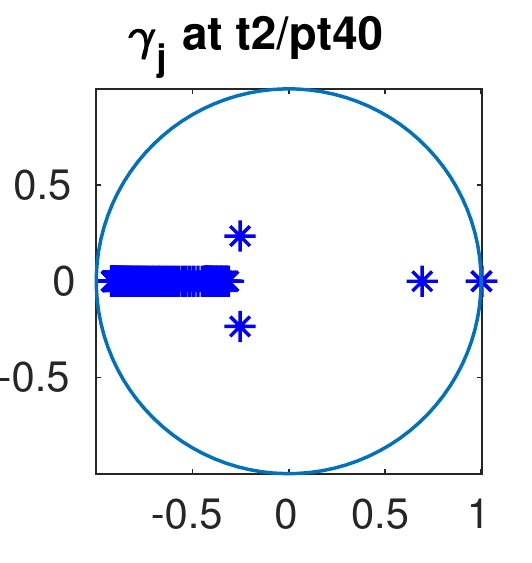}
\ece
\vs{-5mm}
   \caption{{\small Results for \reff{ebru} from 
{\tt bru1dcmds\_b.m}.   
(a) Bifurcation diagram, extenting Fig.~\ref{brfig1}(c) by 
the secondary pitchfork bifurcations from the primary Hopf branch (orange); 
first PD branch {\tt t1} in red, second PD branch {\tt t2} in magenta. 
(b,c) solution and Floquet plots. 
The magenta branch intermediately has index 3, before 
for increasing $b$ first two unstable multipliers come back 
inside the unit circle via a Neimark-Sacker scenario, and 
then the last unstable $\ga$ goes through $1$ and the orbit gains stability 
near $b=3.3$.  
  \label{brfig2}}}
\end{figure}

\hulst{caption={{\small {\tt \dname/bru1dcmds\_b.m} (first 7 lines). 
The remainder computes the 2nd bifurcating branch and deals with plotting.   }},label=brul3,language=matlab,stepnumber=5, linerange=3-9}{\dhome/bru1dcmds_b.m}

\subsection{2D}\label{br2dsec}
We close the discussion of \reff{ebru} with some comments on the continuation  of solution branches in 2D. See \cite[Fig.10]{hotheo} for example results, 
where we consider \reff{ebru} on $\Om{=}(-\pi/2,\pi/2){\times}(-\pi/8,\pi/8)$.  
Already on this rather small domain the linearization 
of \reff{ebru} around $U_s$ has many small real eigenvalues. 
Therefore, the Hopf eigenvalues (with imaginary parts near $\om_1=1$) are 
impossible to detect by computing just a few eigenvalues close to $0$, 
as illustrated in Fig.~\ref{bf1}(a,b). 
Thus, the preparatory step {\tt initeig} already used in 1D in line 5 
of Listing \ref{brul1} becomes vital. This uses a Schur complement algorithm to compute a guess for the spectral shift $\om_1$ near which we expect 
Hopf eigenvalues during the continuation of a steady branch, see \cite[\S2.1 and Fig.8]{hotheo} for illustration. 

\begin{figure}[t]
\bce{\small 
\begin{tabular}{llll}
(a) $\neig=200$&(b) $\neig=300$&(c) $|g|$ from \cite[(2.11)]{hotheo}&(d) \begin{tabular}{l}
$\neig=(3,3)$ with\\ $\om_1=0.9375$\end{tabular}\\
\ig[width=0.2\tew, height=35mm]{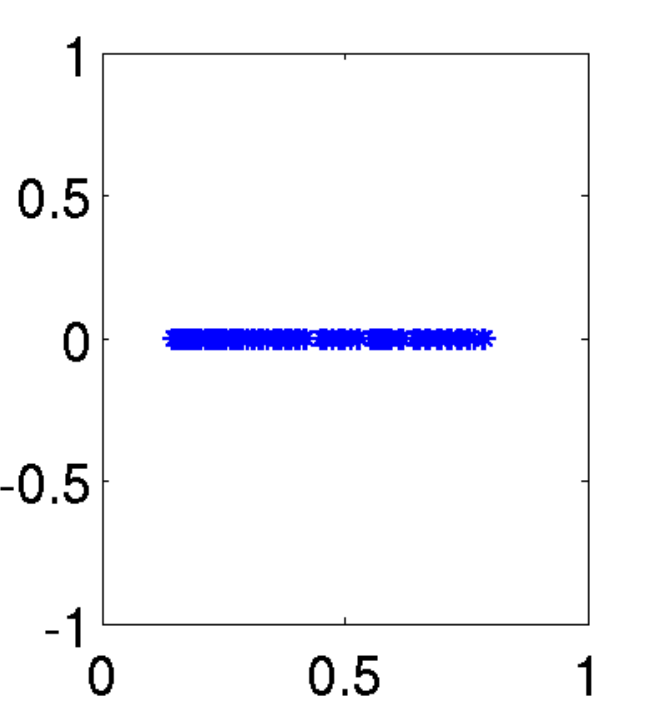}
&\ig[width=0.21\tew, height=35mm]{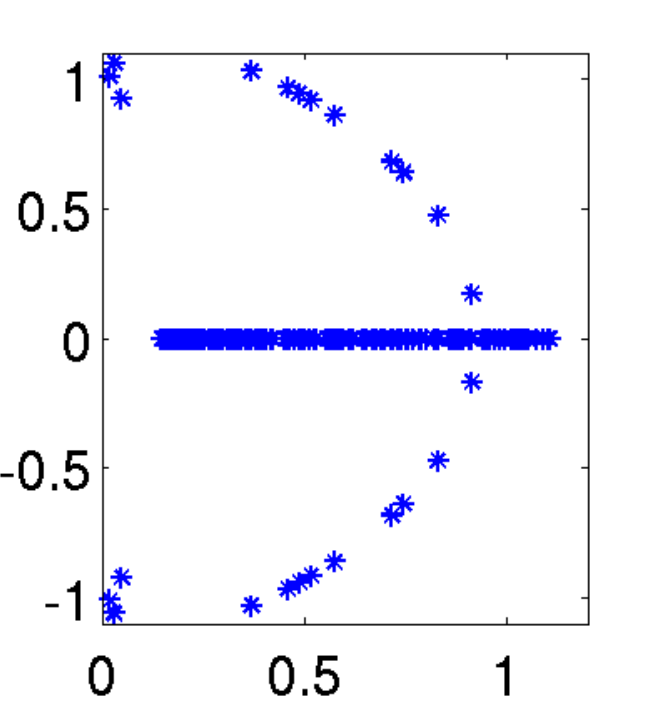}
&\ig[width=0.23\tew, height=38mm]{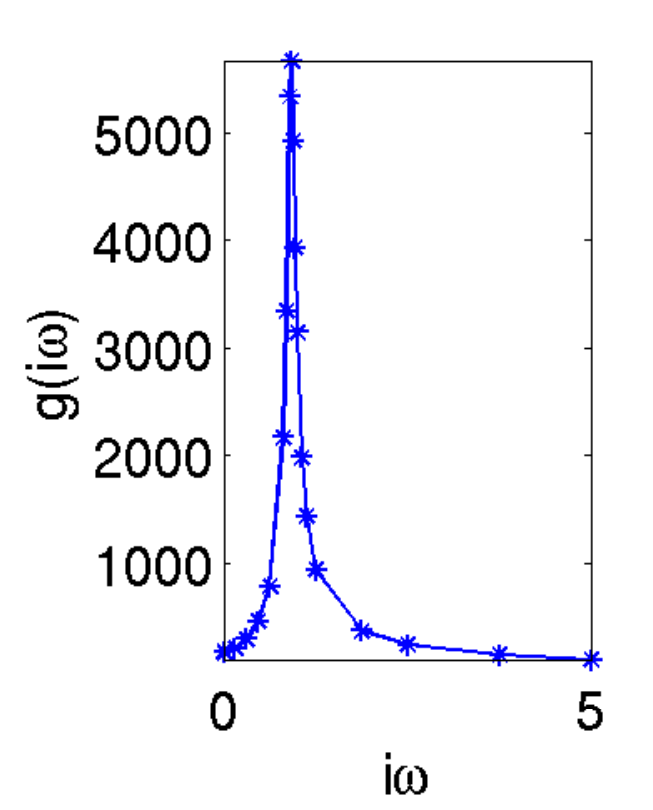}
&\raisebox{0mm}{\ig[width=0.19\tew, height=35mm]{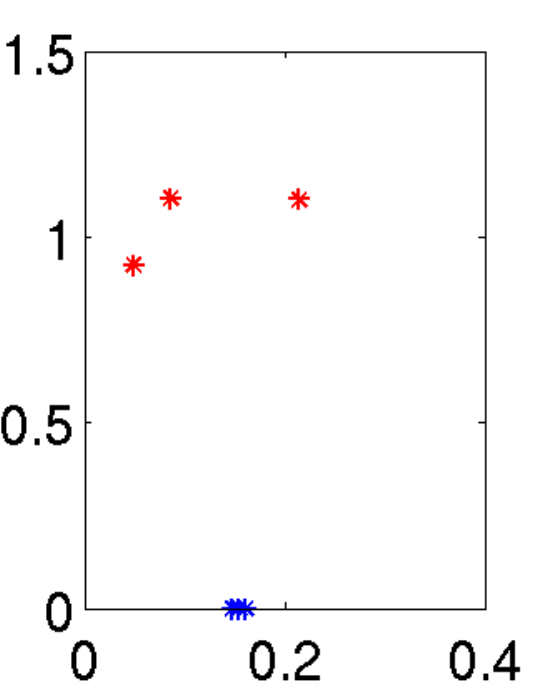}}
\end{tabular}}
\ece 

\vs{-5mm}
   \caption{{\small (a,b) $\neig$ smallest eigenvalues of the linearization 
of \reff{ebru} around 
$U_s$ at $b=2.75$, remaining parameters from \reff{brupar}; 
\heda\  with $\neig=200$  will not detect any Hopf points. 
(c) \reff{reso} yields a guess $\om_1=0.9375$ for the $\om$ value at 
Hopf bifurcation, and then \hedb\ with $\neig=(3, 3)$ is reliable and 
fast: (d) shows the three eigenvalues closest to $0$ in blue, 
and the three eigenvalues closest to $\ri \om_1$ in red. 
  \label{bf1}}}
\end{figure}

\hulst{caption={{\small {\tt \dname/bru2dcmds.m} (first Cell, 
i.e., initialization). 
The main issue is the preparatory step in line 6. This produces 
a (here quite accurate) guess 0.9375 for the candidate $\om$ for imaginary 
parts at Hopf bifurcations, which, together with $\om_0=0$, is put 
into {\tt p.nc.eigref}. The remainder of {\tt brucmds2.m} then 
continues the homogeneous branch and some bifurcating Hopf and Turing 
branches, including secondary bifurcations from the Turing branch to 
'spotted' Hopf branches. Here, an adaptive spatial mesh refinement 
is helpful to increase accuracy.}},label=brl4,language=matlab,stepnumber=5, firstnumber=1, linerange=2-7}{\hdhome/brussel/bru2dcmds.m}

\hulst{caption={{\small {\tt \dname/e2rsbru.m}. For a general 
discussion of error estimators in the \oop\ setting we refer to \cite{actut}. 
The only difference is that here we have a 3 component system, and thus 
we sum up the element wise errors over the components.}},label=brl7,language=matlab,stepnumber=5, firstnumber=1}{\hdhome/brussel/e2rsbru.m}

As indicated in the caption of Listing \ref{brl4},
 at the start of the (1D and 2D) Turing branches 
we do some adaptive mesh--refinement.  The used 
error estimator is given in Listing \ref{brl7}. The further BPs and HBPs then 
obtained are very close to the BPs and HBPs on the coarser mesh, 
but the resolution of the bifurcating Hopf branches becomes considerably 
better, with a moderate increase of computation time, which in any case 
is faster than starting with a uniform spatial mesh yielding a comparable 
accuracy. 


\section{A canonical system from optimal control: Demo {\tt pollution}}\label{ocsec}
\def\dname{pollution}
In \cite{U16,GU17}, \pdep\ has been used to study infinite time horizon distributed optimal control (OC) problems, 
see also \cite{octut} for a tutorial on OC computations 
with \pdep. 
As an example for such problems with Hopf bifurcations%
\footnote{which so far could not be found 
in the systems studied in \cite{U16, GU17}}  we consider, following 
\cite{wirl00}, a model in which the states $v_1=v_1(t,x)$ and 
$v_2=v_2(t,x)$ are the emissions of some firms and the pollution stock, and the control $k=k(t,x)$ models the abatement policy of the firms. The objective 
is to maximize 
\begin{subequations} \label{oc1}
\hual{
    &J(v_0(\cdot),k(\cdot,\cdot)):=\int_0^\infty\er^{-\rho t}
J_{ca}(v(t),k(t)) \dd t,
}
\end{subequations}
where $\ds J_{ca}(v(\cdot,t),k(\cdot,t))=\frac 1{|\Om|}
\int_\Om J_c(v(x,t),k(x,t))\dd x$ 
is the spatially averaged current value function, with 
local current value
$J_c(v,k)=pv_1-\beta v_2-C(k)$, $C(k)=k+\frac 1 {2\ga} k^2$, 
where $\rho>0$ is the discount rate. Using Pontryagin's Maximum 
Principle, the so called canonical system for the states $v$ and 
co-states (or Lagrange multipliers or shadow prices) $\lam$  can be formally derived as a first order necessary optimality condition, 
using the intertemporal transversality condition 
\huga{\label{tcond}
\lim_{t\to\infty}\er^{-\rho t}\int_\Om \spr{v,\lam}\dd x=0.
}

The canonical system reads 
\begin{subequations} \label{cs}
\hual{
\pa_t v&=D\Delta v+f_1(v,k), \quad v|_{t=0}=v_0, \\
\pa_t \lam&=-D\Delta\lam+f_2(v,k), 
}
where $f_1(v,k)=(-k,v_1-\al(v_2))^T$, $f_2(v,k)=(\rho\lam_1-p-\lam_2, (\rho+\al'(v_2))\lam_2+\beta)^T$, 
$\pa_{\bf n} \lam=0$ on $\pa\Om$, and where the control $k$ is given by 
\huga{\label{kform}
k=k(\lam_1)=-(1+\lam_1)/\ga.
}
\end{subequations}
For convenience we set 
$u(t,\cdot):=(v(t,\cdot), \lam(t,\cdot)): \Om\ra\R^{4}$,  
and write \reff{cs} as 
 \hual{\label{cs2}
&\pa_t u=-G(u):=\CD\Delta u+f(u),
}
where $\ds 
\CD=$diag$(d_1,d_2,-d_1,-d_2)$, 
$\ds f(u)=\biggl(-k, v_1-\al(v_2), 
\rho\lam_1-p-\lam_2, (\rho+\al'(v_2))\lam_2+\beta\biggr)^T$. 
Besides the boundary condition $\pa_{\bf n} u=0$ on $\pa\Om$ 
we have the initial condition $v|_{t=0}=v_0$ (only) for the states. 
A solution $u$ of the canonical system \reff{cs2} is called a 
\emph{canonical path}, 
and a steady state of \reff{cs2} (which automatically fulfills \reff{tcond}) 
is called a \emph{canonical steady state (CSS)}. 
Due to the backward diffusion in $\lam$, and since we only have 
initial data for half the variables, \reff{cs2} is {\em not} well posed 
as an initial value problem. Thus, one method for OC problems of type \reff{oc1} is to first study CSS,  and then canonical paths connecting 
some initial states to some CSS $u^*$. This requires the so-called saddle-point property for 
$u^*$, and if this holds, then canonical paths 
to $u^*$ can often be obtained from a continuation process in the initial 
states, see \cite{octut}. 

A natural next step is to search for 
time--periodic solutions $u_H$ of canonical systems, which obviously 
also fulfill 
\reff{tcond}. The natural generalization of the saddle point property 
is that 
\huga{\label{ddef2}
d(u_H):=\ind(u_H)-\frac{n_u}2=0,  
}
i.e., that exactly half of the Floquet multipliers are in 
the unit circle.
In the (low--dimensional) ODE case, there then exist methods to compute 
connecting orbits to (saddle type) 
periodic orbits $u_H$ with $d(u_H)=0$,   see \cite{BPS01, grassetal2008}, 
which require comprehensive information on the 
Floquet multipliers and the associated eigenspace of $u_H$. 
A future aim is to extend these methods to periodic orbits of PDE OC systems. 

However, in \cite[\S3.4]{hotheo} we only illustrate that 
Hopf orbits can appear as candidates for optimal solutions 
in OC problems of the form \reff{oc1}, and that 
the computation of Floquet multipliers via the periodic 
Schur decomposition {\tt floqps} can yield reasonable results, even when 
computation via {\tt floq} completely fails. 

For all parameter values, \reff{cs2} has an explicit spatially homogeneous 
CSS, see \cite{hotheo}, and by a suitable choice of parameters we obtain Hopf bifurcations to spatially homogeneous and 
spatially patterned time periodic orbits.  Concerning the implementation, Table \ref{ptab} gives an overview of the 
involved scripts and functions. Since we again use the \oop\ setting, 
and since we restrict to 1D, although \reff{cs2} is a four component system, much of this is very similar to the {\tt cgl} demo in 1D, 
with the exceptions that: (a) we also need to implement the objective value and 
other OC related features; (b) similar to {\tt brussel} it is useful to prepare the detection of HBPs via {\tt initeig}; (c) we need to use {\tt flcheck=2} throughout. 
Thus, in Listings \ref{l41}-\ref{l46} we comment on these points, and for 
plots illustrating the results of running {\tt pollcmds.m} refer to \cite[\S3.4]{hotheo}. 

\taskip
\begin{table}[ht]\caption{Main scripts and functions in {\tt hopfdemos/pollution}. 
\label{ptab}}
{\small 
\begin{tabular}{l|p{0.69\tew}} 
script/function&purpose,remarks\\
\hline
pollcmds&main script\\%
p=pollinit(p,lx,nx,par)&init function\\
p=oosetfemops(p)&set FEM matrices (stiffness K and mass M) \\
r=pollsG(p,u)&encodes $G$ from \reff{cs2}; we avoid implementing the Jacobian here and instead use {\tt p.sw.jac=1} \\
f=nodalf(p,u)&nonlinearity, called in {\tt sG}.\\
jc=polljcf(p,u)&the (current value) objective function \\
\end{tabular}
}
\end{table}\teskip

\hulst{caption={{\small {\tt \dname/pollinit.m} (first 4 lines). 
Additional to the rhs, in line 3 we set a function handle to the objective 
value, as usual for OC problems (see \cite{octut}). Similarly, in line 4 we set 
{\tt p.fuha.outfu} to a customized branch output, which combines features 
from the standard Hopf output {\tt hobra} and the standard OC output {\tt ocbra}. We 
do not set {\tt p.fuha.sGjac} since for convenience here we use numerical 
Jacobians ({\tt p.sw.jac=0} in line 1). The remainder of {\tt pollinit.m} is as usual. }},
label=l41,language=matlab,stepnumber=5, firstnumber=1,lastline=4}{\hdhome/pollution/pollinit.m}

\hulst{caption={{\small {\tt \dname/polljcf.m}, function to compute the current objective value. Called in {\tt pollbra} to put the value 
on the branch (for plotting and other post-processing). }},
label=l45,language=matlab,stepnumber=5, firstnumber=1}{\hdhome/pollution/polljcf.m}

\def\dname{pollution}
\hulst{caption={{\small {\tt \dname/pollcmds.m} (first 19 lines). In cell 1 we use {\tt initeig} to generate a guess for $\ri \om_1$ (the Hopf wave number), and set {\tt neig} to compute 5 eigenvalues near 0 and near $\om_1$. 
For the computation of multipliers here we need to use {\tt floqps}, see line 17. The remainder of {\tt pollcmds} deals with plotting. 
}},
label=l46,language=matlab,stepnumber=5, firstnumber=1,lastline=19}{\hdhome/pollution/pollcmds.m}

\section{Hopf bifurcation with symmetries}\label{hosymsec}
If the PDE \reff{tform} has (continuous) symmetries, then already for 
the reliable continuation of steady states it is often 
necessary to augment \reff{tform} by $n_Q$ suitable phase conditions, in the form 
\huga{\label{pcgen}
Q(u,\lam,w)=0\in\R^{n_Q}
}
where $w\in\R^{n_Q}$ stands for the required $n_Q$ additional active 
parameters,  
see \cite{symtut} for a review. For instance, if \reff{tform} is spatially 
homogeneous and we consider periodic BC, then we have a translational invariance, and (in 1D) typically augment \reff{tform} by the phase condition 
\huga{\label{pctrans}
\spr{\pa_x u^*,u}=0\in\R, 
}
where (for scalar $u,v$) $\spr{u,v}=\int_\Om uv\dd x$, and 
where $u^*$ is either a fixed reference profile or the solution from the 
previous continuation step. We thus have $n_Q=1$ additional equations, and consequently must free 1 additional parameter. 

Similarly, we must add phase conditions to the computation of Hopf orbits 
(additional to the phase condition \reff{pca} fixing the translational 
invariance in $t$). This is in general not straightforward, since \reff{pcgen}, 
with \reff{pctrans} as an example, is not of the form $\pa_t u=Q(u,\lam)$ and 
thus cannot simply be appended to \reff{tform}. Instead, the steady phase 
conditions \reff{pcgen} must be suitably modified and explicitly appended to the 
Hopf system, see \reff{fsac}. Examples for the case \reff{pctrans} have been discussed in \cite[\S4]{symtut}, namely 
the cases of modulated fronts, and of breathers. 

Here we give two more examples, and extend the breather example to compute 
period doubling bifurcations. The first example deals with Hopf orbits in a reaction diffusion system 
with mass conservation, and the second with Hopf orbits in the Kuramoto-Sivashinsky (KS) equation, where we need two phase conditions, 
one for mass conservation and one to fix the translational invariance. 
For both problems we restrict to 1D; like, e.g., the cGL equation, they both can immediately be transferred to 2D (where for the KS equation we need a third 
phase condition $q_3(u)=\spr{\pa_y u^*,u}=0$, cf.~(\ref{ksb}c)), but 
the solution spaces and bifurcations then quickly become ``too rich'', 
such that -- as often -- 2D setups only make sense if there are specific questions 
to be asked.

\subsection{Mass conservation: Demo {\tt mass-cons}}
\def\dname{mass-cons}
As a toy problem for mass conservation in a reaction diffusion system we consider 
\huga{\label{rdm}
\pa_t u_1=\Delta u_1+d_2\Delta u_2+f(u_1,u_2), \quad 
\pa_t u_2=\Delta u_2-f(u_1,u_2), \text{ in } \Om, 
}
$f(u_1,u_2)=\al u_1-u_1^3+\beta u_1u_2$, 
with parameters $d_2,\al,\beta\in\R$, and homogeneous Neumann BC. 
Then $m:=\frac 1 {|\Om|}\int u+v \dd x$ is conserved since $\ddt \int_\Om (u+v)\dd x
=\int_{\pa\Om}\pa_n u_1+(1+d_2)\pa_n u_2\dd S=0$. Given a steady state $(u,v)$ for some fixed $\al,\beta$, this always comes in a continuous family parameterized by the ``hidden'' parameter $m$. Thus, to study 
steady states and their bifurcations we use the mass constraint 
\huga{\label{qs} 
Q(u,\lam):=\frac 1 {|\Om|}\int u+v \dd x-m=0, 
}
where as usual $\lam$ stands for the vector of all parameters. 
Given this additional equation, we have the differential-algebraic system 
\huga{\label{mdae}
 M\dot u=-G(u,\lam), \qquad Q(u,\lam)=0, 
}
and to compute solution branches we need 2 parameters, which we choose as 
$\al,\beta$. If we restrict to $m=0$, then we have two explicit branches 
of homogeneous solutions, namely $u_2=-u_1$ and $u_1{=}-\frac \beta 2\pm \sqrt{\frac{\beta^2} 4+\al}$. We choose the initial point $(\al,\beta){=}(1,1)$, $u_1=-1/2-\sqrt{5/4}, u_2=-u_1$ and continue in $\al$. 

As a Hopf version of \reff{qs} we use 
\huga{\label{mcqH}
Q_H(u(\cdot,\cdot)):=\sum_{i=1}^m \left(\int_\Om (u_1(t_i,x)+u_2(t_i,x))\dd x-m\right)\stackrel!=0, 
} 
see Listings \ref{mcl3} and \ref{mcl4}. In \reff{mcqH} 
i.e., we require the 
average (in $t$) mass to be conserved. Theoretically it would be sufficient 
to require $\int_\Om (u_1(t_0,x)+u_2(t_0,x))\dd x-m=0$, but it turns out that 
\reff{mcqH} is more robust numerically, and that also with \reff{mcqH} 
we have $\ds \left|\int_\Om (u_1(t_i,x)+u_2(t_i,x))\dd x-m\right|\le {\tt tol}$ for all $i$, i.e., pointwise in $t$. 

The implementation of \reff{mdae} is rather straightforward, 
see Table \ref{mctab} for an overview, and Listings \ref{mcl2}--\ref{mcl4}. 
We fix $d_2=10$ and restrict to 1D, namely $\Om=(-\pi,\pi)$. 
Figure \ref{mcf1}(b) shows a basic bifurcation diagram, with various quantities as functions of $\al$.  The continuation of 
\reff{mdae} in $\al$ with fixed $m=0$ yields that the homogeneous solution $u$ stays fixed, i.e., $u_1=-1/2-\sqrt{5/4}, u_2=-u_1$ for all $\al$, and that only $\beta$ is adjusted, see the black lines in (b). (c) shows a number of 
Hopf orbits, where on each orbit we have $|Q(u(t,\cdot)|<10^{-8}$ 
(see also the last plot in (a) for the average $Q_H$), 
where the tolerance for the Hopf orbits is $10^{-6}$. These Hopf orbits 
are all unstable according to the associated Floquet multipliers,  
see also Fig.~\ref{mcf1b}(a), and thus it is interesting to see the evolution of solutions 
starting on a Hopf orbit (with the numerical error acting as a perturbation 
of the true point on a Hopf orbit). In Fig.~\ref{mcf1b}(b) we exemplarily show this for 
the case of $u(0)$ from {\tt h3/pt15}; here, as in all other cases we 
considered, the time evolution converges to another stable spatially homogeneous 
steady state. From this we may again start continuation in, e.g., 
$\al$ and $\beta$, and find that this branch again typically shows some Hopf bifurcations. 

\begin{figure}[ht]
\bce{\small
\begin{tabular}{l}
(a) BDs, parameter $\beta, \min(u_1), \max(u_2$ and mass as functions of $\al$.   \\
\ig[width=0.22\tew,height=45mm]{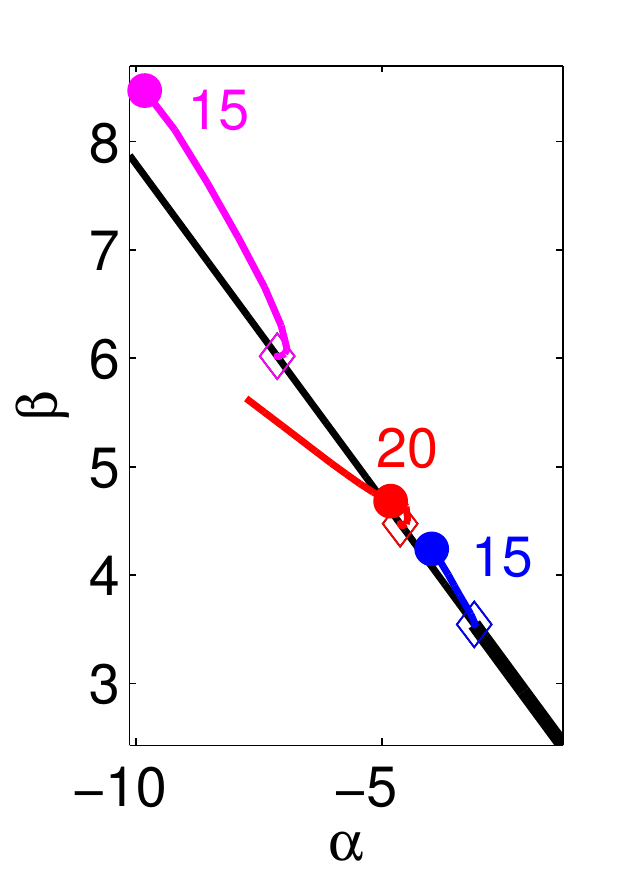}
\ig[width=0.22\tew,height=45mm]{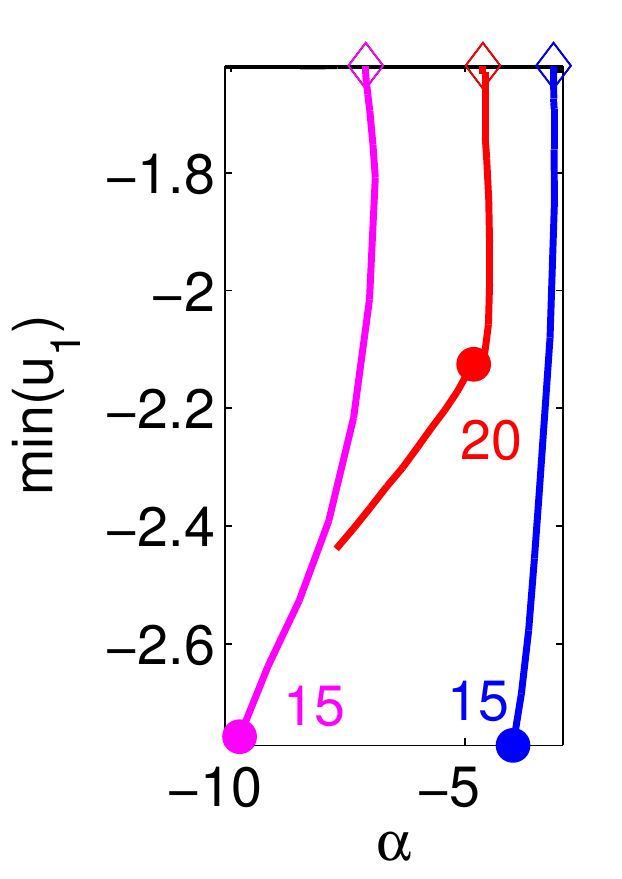}
\ig[width=0.22\tew,height=45mm]{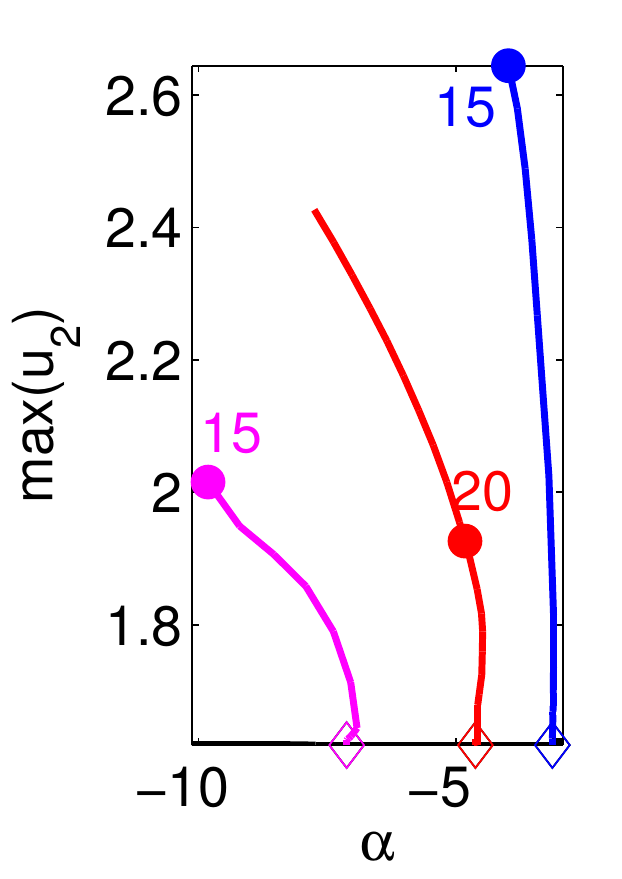}
\ig[width=0.22\tew,height=45mm]{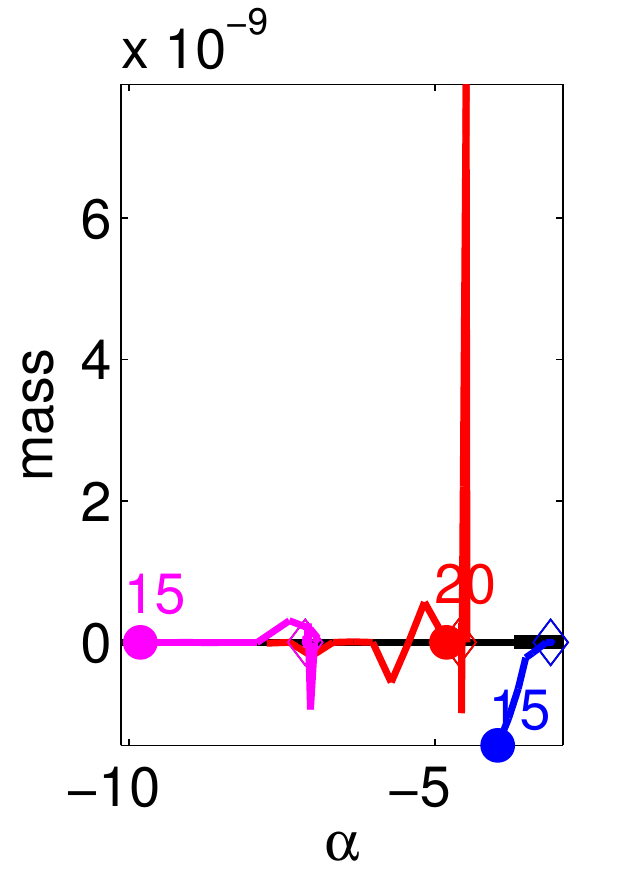}\\
(c) Selected solution plots (both components for {\tt h1/pt15}) \\
\ig[width=0.23\tew]{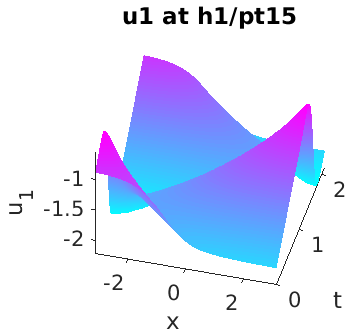}\hs{-2mm}
\ig[width=0.23\tew]{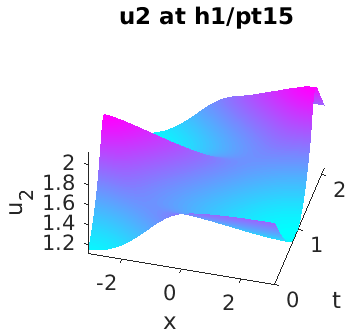}\hs{-2mm}
\ig[width=0.23\tew]{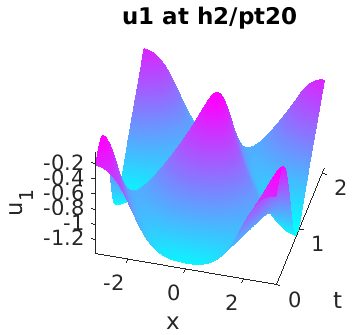}\hs{-2mm}
\ig[width=0.23\tew]{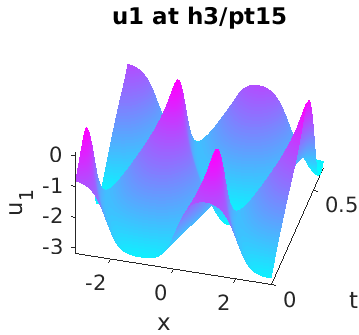}
\end{tabular}}
\ece 

\vs{-5mm}
   \caption{{\small Continuation in $\al$ for \reff{mdae} with $m=0$. 
(a) Branch data on the homogeneous branch (black) and on three Hopf branches h1 (blue), h2 (red), and h3 (magenta). (b) Example solution plots.
  \label{mcf1}}}
\end{figure}

\begin{figure}[ht]
\bce{\small
\begin{tabular}{l}
(a) Floquet multipliers for {\tt h1/pt15}\\ 
\hs{-0mm}\ig[width=0.8\tew]{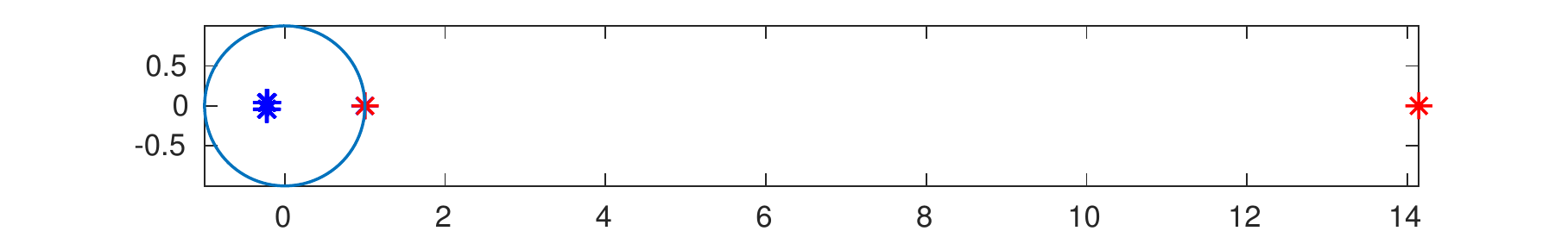}\\
(b) Time evolution of small perturbation of $u$ from {\tt h1/pt15} at $t=0$, 
left $u_1$, right $u_2$\\
\ig[width=0.25\tew]{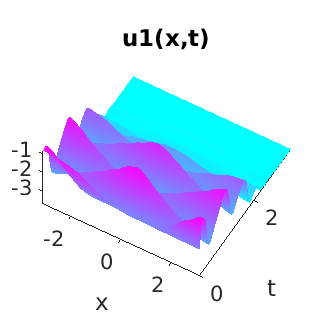}\ig[width=0.25\tew]{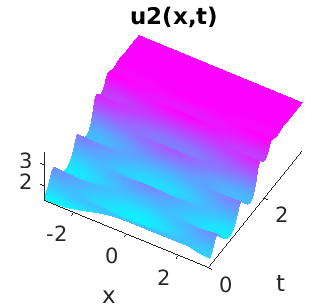}
\end{tabular}}
\ece 

\vs{-5mm}
   \caption{{\small (a) Instability of {\tt h1/pt15} as seen in its 
Floquet multipliers. (b) time integration, with convergence to another 
spatially homogeneous steady state. 
  \label{mcf1b}}}
\end{figure}

\taskip
\begin{table}[ht]\caption{Scripts and functions in {\tt hopfdemos/mass-cons}. 
\label{mctab}}
{\small 
\begin{tabular}{l|p{0.55\tew}}
script/function&purpose,remarks\\
\hline
cmds1d&main script \\%
mcinit, oosetfemops, sG, sGjac, nodalf&initialization, FEMops, rhs, Jac., and nonlinearity, as usual.\\
qf, qfjac&the phase condition \reff{qs}, and its Jacobian. \\
qfh, qfhjac&the Hopf version \reff{mcqH} of \reff{qs}, and its Jacobian.\\
\end{tabular}
}
\end{table}\teskip

\hulst{caption={{\small {\tt \dname/qf.m}; mass constraint for 
steady state computations. }},
label=mcl2,language=matlab,stepnumber=5, firstnumber=1}{\hdhome/mass-cons/qf.m}

\hulst{caption={{\small {\tt \dname/qfh.m}; Hopf setting of 
mass constraint. The summing up (in $t$) 
of the masses turns out to be more robust, with the mass-constraint 
actually fulfilled pointwise (in $t$).  }},
label=mcl3,language=matlab,stepnumber=5, firstnumber=1}{\hdhome/mass-cons/qfh.m}

\hulst{caption={{\small {\tt \dname/qfhder.m}, $u$--derivatives 
of $Q_H$, cf.~last line of 
\reff{news}, where the parameter derivatives are done automatically 
via finite differences. }},
label=mcl4,language=matlab,stepnumber=5, firstnumber=1}{\hdhome/mass-cons/qfhjac.m}

\hulst{label=mcl1a,language=matlab,stepnumber=5, 
linerange=2-20}{\hdhome/mass-cons/cmds1d.m}
\hulst{caption={{\small {\tt \dname/cmds1d.m} (with some omissions) 
C1 continues the 
homogeneous branch, giving a number of Hopf bifurcations; here $u_1=-(1+\sqrt{5})/2$ and $u_2=-u_1$ stay fixed, 
and only $\beta$ varies with $\al$. In C2 we follow the first three 
Hopf branches, where we replace the stationary Hopf constraint in {\tt qf} 
by the Hopf version {\tt qfh}, see Fig.~\ref{mcf1} for bifurcation 
diagrams and example Hopf solutions. All Hopf branches turn out to be unstable 
(from the Floquet multipliers), and thus in C3-5 
we exemplarily look into the time evolution from the first point ($t=0$) on 
the Hopf orbit {\tt h1/pt15}. This converges to a (stable) homogeneous solution again, but at larger amplitude. Finally in C6 we use this as 
a starting point for continuation in $\al$, and again find a number 
of Hopf bifurcations for decreasing $\al$.  }},
label=mcl1,language=matlab,stepnumber=5, linerange=31-36}{\hdhome/mass-cons/cmds1d.m}

\subsection{Mass and phase constraints: Demos {\tt kspbc4} and {\tt kspbc2}}\label{kssec}
\def\dname{kspbc4}
The Kuramoto-Sivashinsky (KS) equation \cite{kura-tsu76, siva77} is a canonical and much studies model for 
long--wave instabilities in dissipative systems, for instance in laminar flame 
propagation, or for surface instabilities of thin liquid films. Here we 
consider the KS equation in the form 
\huga{\label{KS}
\pa_t u=-\al\pa_x^4 u-\pa_x^2 u-\frac 1 2 \pa_x(u^2), 
}
with parameter $\al>0$, on the 1D domain $x\in(-2,2)$ with periodic BC. 
\reff{KS} is thus translationally invariant, and has the boost invariance 
$u(x,t)\mapsto u(x-ct)+c$, and we need two phase conditions, 
\begin{subequations}\label{kspc}
\hual{
\frac 1{|\Om|}\int_\Om u\dd x=m,&\text{ fixing the mass $m$,}\\
\spr{\pa_x u^*,u-u^*}=0,&\text{ fixing the translational invariance.}
}
\end{subequations}
Here fixing $m=0$, \reff{KS} shows bifurcations from the trivial solution $u\equiv 0$ to stationary spatially periodic solutions at $\ds \al_k=\left(\frac 2 {k\pi}\right)^2$, $k\in\N$. Next, for 
decreasing $\al$ we obtain 
secondary Hopf bifurcations from some branches of steady patterns, and for $\al\ra 0$ 
the dynamics become more and more complicated, making \reff{KS} a model 
for turbulence. 
In \cite{vV17}, a fairly 
complete bifurcation diagram (with $\al$ in the range $0.025$ to $0.4$) has been obtained for \reff{KS} 
on $\Om=(0,2)$ with {\em Dirichlet} BC, i.e., $u(0,t)=u(2,t)=\pa_x^2u(0,t)=\pa_x^2u(2,t)=0$, where in particular many  bifurcations have been explained analytically as hidden symmetries by 
extending solutions antisymmetrically to the domain $(-2,2)$ with periodic BC. 

Here we directly study \reff{KS} in this setting, giving us the opportunity to also explain how to setup 4th order equations and periodic BC in \pdep. 
For the latter we only need to call {\tt  p=box2per(p,1)}, which generates 
matrices {\tt fill} and {\tt drop} which are used to transform the FEM 
matrices such as $M$ and $K$ to the periodic setting, see \cite{pbctut}. 
In order to implement 4th order equations there are basically two options: 
\bci
\item[(i)] Since $-\pa_x^2 u=M^{-1}Ku$ in the FEM sense, 
\reff{KS} can be written in the \pdep\ FEM setting 
as $M\pa_t u=-\al KM^{-1}K u+Ku-\frac 1 2 K_x(u^2)$. For pBC, $K,M$ commute, 
and thus we can multiply by $M$ to obtain 
$M^2\pa_t u=-\al K^2 u+MKu-\frac 1 2 MK_x(u^2)$. 
Then letting $M_0{=}M$ and 
redefining $M{=}M^2$ we obtain $M\pa_t u{=}-\al K^2 u+M_0Ku-M_0K_x(u^2)$. 
To incorporate the phase conditions \reff{kspc} we introduce 
the parameters $s$ for phase-conservation and $\eps$  for mass conservation, 
and thus ultimately consider the system 
\begin{subequations}\label{ksb}
\hual{
M\dot u&=-\al K^2 u+M_0Ku-\frac 1 2M_0K_x(u^2)+sK_xu+\eps, \\
0&=q_1(u):=\frac 1 {|\Om|}\sum_{i=1}^{n_u}(M_0 u)_i-m, \\
0&=q_2(u):=\spr{\pa_x u^*,u-u^*}, 
}
\end{subequations}
where $\frac 1 {|\Om|}\sum_{i=1}^{n_u}(M_0 u)_i$ is the (Riemann sum) approximation of $\frac 1 {|\Om|}\int_{\Om} u\dd x$, and $u^*$ is a 
suitable reference profile. This set up is implemented in {\tt kspbc4}, 
see below. 

\item[(ii)] Alternatively we can rewrite the 4th order equation as a 2 component 2nd order system, for 
instance for \reff{KS} in the form 
\huga{\label{ks2}
\pa_t u=-\al\pa_x^2 v-\pa_x^2 u-\frac 1 2 \pa_x(u^2), \qquad 
0=-\pa_x^2 u+v. 
} 
By exploiting the mass matrix on the lhs of \reff{tformd}, \reff{ks2} can be straightforwardly implemented in \pdep\ in the form 
\huga{\label{ks2p2pa}
\CM \dot U{=}-G(U),\quad U{=}\bpm u_1\\ u_2\epm,\  
\CM{=}\bpm M&0\\0&0\epm, \  G(U){=}-\bpm K&\al K\\K&M\epm U{+}\bpm\frac 1 2 K_x(u_1^2)\\0\epm, 
}
where $M,K$ and $K_x$ are the scalar mass, stiffness and advection matrices. 
Importantly, the spectral picture and time evolution for \reff{ks2p2pa} are still fully equivalent to \reff{KS}. Adding phase conditions like 
(\ref{ksb}b,c), this is implemented in {\tt kspbc2}, and yields 
the same results as the {\tt kspbc4} set up, except for small differences 
wrt to Floquet multipliers, which in any case are somewhat delicate for constrained 
Hopf orbits, see Remark \ref{flremc}. 
\eci

The implementation of \reff{ksb} is rather straightforward. See Table \ref{kstab} 
for an overview, 
Listings \ref{ksl1}--\ref{ksl5} for pertinent sections from {\tt oosetfemops, 
cmds1d, 
sG, qf} and {\tt qfh}, while for {\tt cmds2.m} and Jacobians/derivatives of 
{\tt sG, qf} and {\tt qfh} we refer to the m-files  {\tt sGjac, qjac} and {\tt qfhjac}, respectively. 

\taskip
\begin{table}[ht]\caption{Scripts and functions in {\tt hopfdemos/kspbc4}. 
\label{kstab}}
{\small 
\begin{tabular}{l|p{0.79\tew}}
script/function&purpose,remarks\\
\hline
cmds1&main script, steady state branches, and associated Hopf bifurcations of standing waves \\%
cmds2&script for one traveling wave branch, and associated Hopf bifurcations 
of modulated traveling waves\\
ksinit, oosetfemops&initialization and FEMops; this is somewhat different from 
the other examples. {\tt ksinit} also contains the call {\tt p=box2per(p,1)} 
to set up the periodic BC; {\tt oosetfemops} contains calls of {\tt filltrafo}, 
and specifically the redefinition of $M$ as $M^2$. \\
sG, sGjac&rhs, Jacobian; again somewhat different from before due to 4th 
order derivatives.\\ 
qf, qfjac&the phase conditions (\ref{ksb}b,c), and the derivatives\\
qfh, qfhjac&the Hopf version of {\tt qf} and its derivative\\
\end{tabular}
}
\end{table}\teskip

\begin{figure}[ht]
\bce{\small
\begin{tabular}{lll}
(a) BD of steady branches&(b) Zooms into BD, including Hopf branches&
(c) Profiles at HBPs \\
\ig[width=0.32\textwidth,height=50mm]{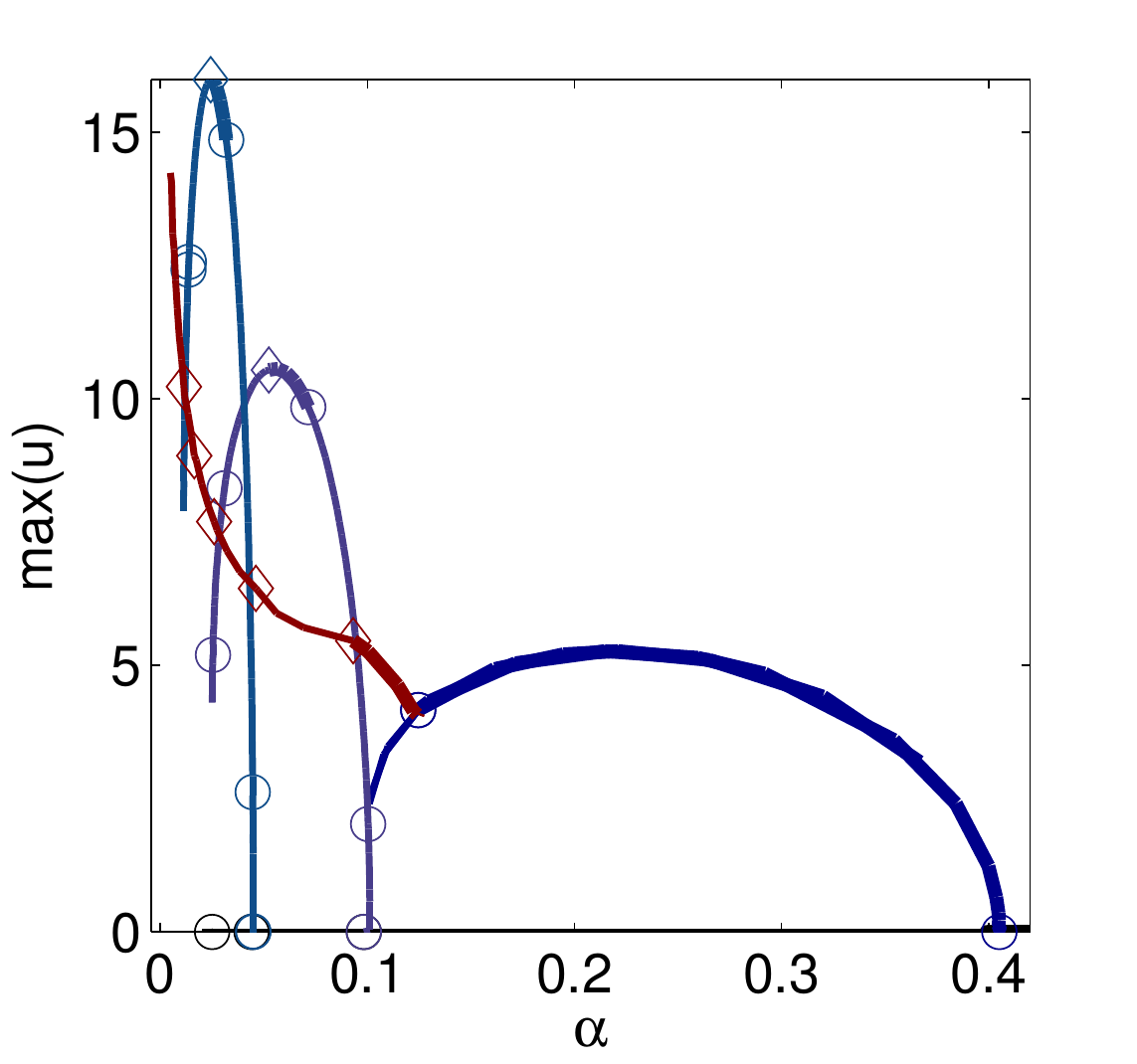}&
\hs{-4mm}\ig[width=0.24\textwidth,height=50mm]{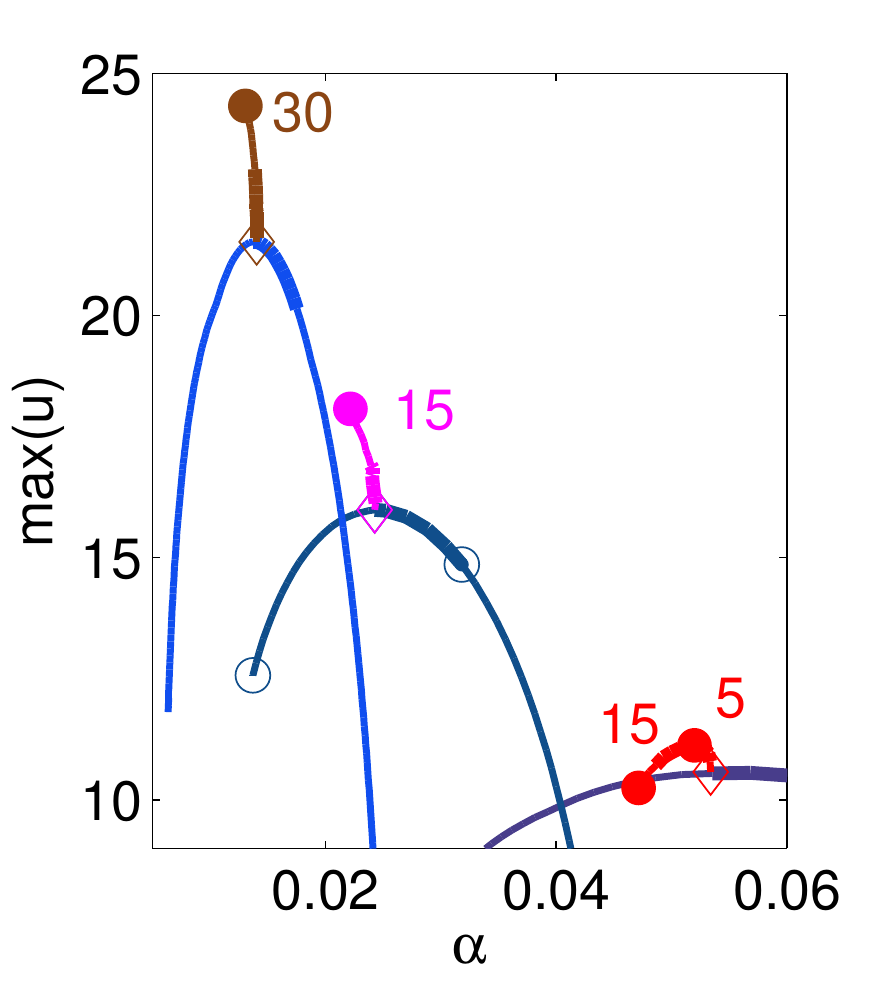}
\hs{-4mm}\ig[width=0.19\textwidth]{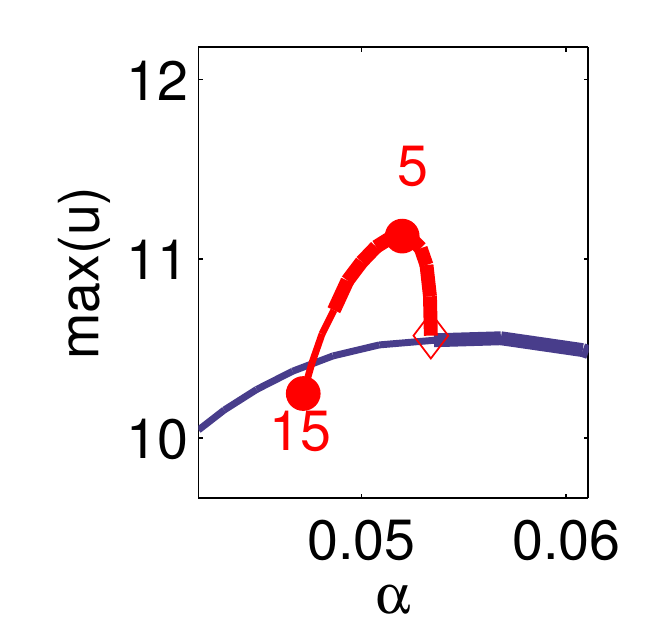}&
\hs{-12mm}\raisebox{25mm}{\begin{tabular}{l}
\ig[width=0.27\textwidth,height=14mm]{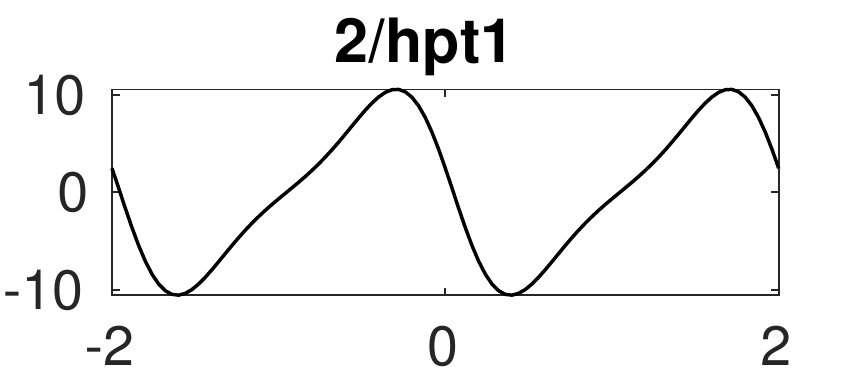}\\
\ig[width=0.27\textwidth,height=14mm]{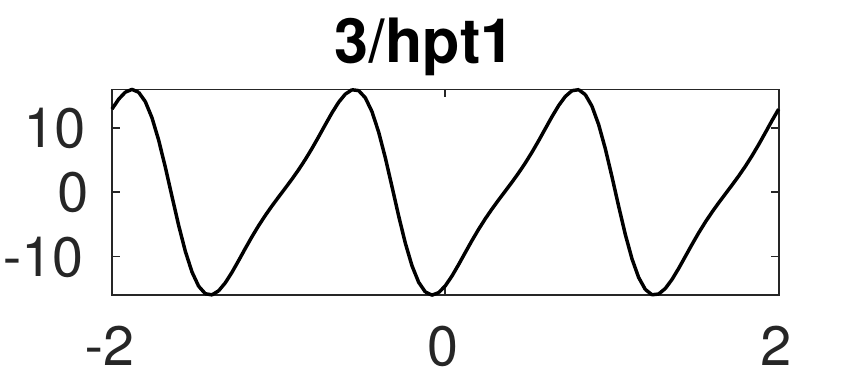}\\
\ig[width=0.27\textwidth,height=20mm]{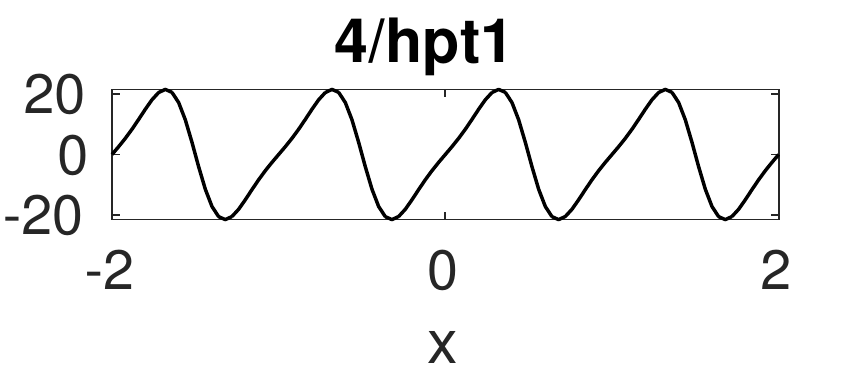}
\end{tabular}}
\end{tabular}
\begin{tabular}{l}
(d) Selected Hopf solutions\\
\hs{-5mm}\ig[width=0.25\textwidth,height=40mm]{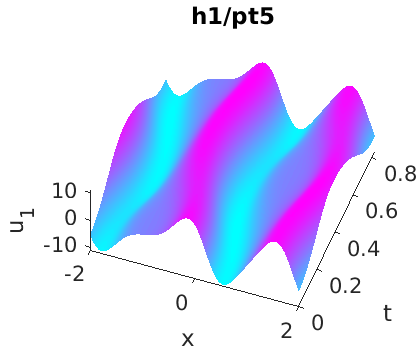}\hs{-2mm}
\ig[width=0.25\textwidth,height=40mm]{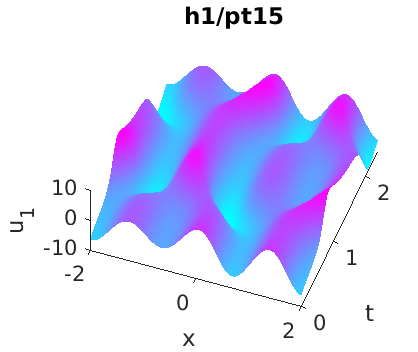}\hs{-2mm}
\ig[width=0.25\textwidth,height=40mm]{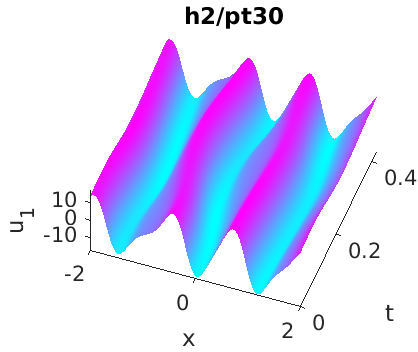}\hs{-2mm}
\ig[width=0.25\textwidth,height=40mm]{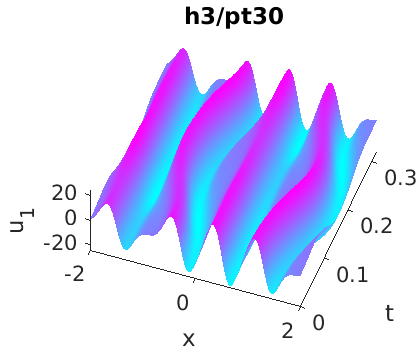}
\end{tabular}
}

\ece 

\vs{-5mm}
   \caption{{\small Results from {\tt kspbc4/cmds1.m}. Bifurcation diagrams 
of steady solutions (except for the brown branch of traveling waves, 
see Fig.~\ref{ksf2}) (a), with 
zoom in (b), including the 4th steady branch and 3 Hopf branches h1 (red), h2 (magenta) and h3 (brown). For all these branches $m=\eps=0$ 
(numerically $\CO(10^{-10})$, and except for the brown branch in (a) 
also $s=0$. 
Profiles at the Hopf bifurcation points the steady branches in (c). In (d) we plot selected Hopf orbits and multiplier 
spectra. {\tt h1} loses stability at $\al\approx 0.0486$ 
 in a pitchfork (a multiplier becoming unstable at $\mu=1$), 
and the largest multiplier of {\tt h1/pt15} is $\mu_2\approx 4000$. Also h2 
is initially stable, but 
looses stability in a pitchfork at $\al\approx 0.024$, i.e., rather close 
to bifurcation, and a similar behaviour occurs on {\tt h3}. See 
Fig.~\ref{ksf1b} for multiplier plots, and  
{\tt cmds2.m} and Fig.~\ref{ksf2} for further plots, for instance of 
solutions on the secondary brown branch in (a), and the Hopf bifurcations 
from this branch.  \label{ksf1}}}
\end{figure}

Figure \ref{ksf1}(a) shows a basic bifurcation 
diagram of steady states, including one branch of traveling waves, 
obtained from {\tt cmds2.m}. As predicted, at $\al_k$ we find supercritical pitchforks of steady branches. The first one starts out stable, and looses 
stability in another supercritical pitchfork around $\al=0.13$ to 
a traveling wave branch (brown), which then looses stability in a Hopf bifurcation, see Fig.~\ref{ksf2}. However,  
here we first focus on Hopf bifurcations from the  2nd and 3rd primary 
branches, which first gain stability at some rather large amplitude, then loose it again in Hopf bifurcations, with the 
solution profiles at the HBP in (c). (b) zooms into the BD at low $\al$, 
including the 4th steady branch, and three Hopf branches, 
while (d) shows selected Hopf orbits. 

These results all fully agree with those in \cite{vV17} (by extending 
the solutions from \cite{vV17} antisymmetrically), who however 
proceed further by also computing some (standing) Hopf 
branches bifurcating in pitchforks and period doublings from the above (standing) Hopf branches. Naturally, these bifurcations are also 
detected in \pdep, but already their localization requires some 
fine tuning, e.g., small stepsizes. Moreover, with the given 
discretizations the branch switching then still often fails. 
This will be further studied elsewhere, and instead we illustrate 
period doubling with a model with better scaling properties in \S\ref{pdsec}.  
On the other 
hand, for our periodic BC on the larger domain we also have 
traveling waves and Hopf bifurcations to modulated traveling waves. 
Some examples for these are considered in {\tt cmds2.m}, see 
Fig.~\ref{ksf2}.

\hulst{caption={{\small {\tt \dname/cmds1.m}. Cell 1 deals with initialization and continuation of the trivial branch. Since the eigenvalues 
$\mu_k=-\al(k\pi/2)^4+(k\pi/2)^2$ of the linearization around $u\equiv 0$ 
have a rather large spacing, in line 7 we set $\mu_{1,2}$ (see \reff{mus}) 
to rather large values. 
In Cell 2 we compute the first 4 branches of steady patterns. The phase condition 
$\spr{\pa_x u^*,u}=0$ (2nd component of {\tt qf}, see Listing 
\ref{ksl4}) is only switched on after a few initial steps and then 
setting the reference profile $\pa_xu^*={\tt p.u0x}$, because 
it only makes sense for $u^*$ not spatially homogeneous. C3 computes 
the first Hopf branch {\tt h1}, bifurcating from steady branch 2. 
We use a rather large Floquet tolerance {\tt p.hopf.fltol}, see \reff{inddef}, 
because the Floquet computations do not remove the neutral directions, 
cf.~Remark \ref{flremc}. Moreover, for this problem {\tt floqps} for the multiplier computations via periodic Schur decomposition sometimes fails 
(for unknown reasons), 
while {\tt floq} (for unknown but maybe related reasons) seems somewhat unreliable for the small multipliers; the large multipliers (and hence the stability information) however 
always seem correct. The remainder 
of {\tt \dname/cmds1.m} deals with the Hopf branches {\tt h2} and {\tt h3}, and with plotting.}},
label=ksl2,language=matlab,stepnumber=5, firstnumber=1, linerange=3-26}{\hdhome/kspbc4/cmds1.m}

\begin{figure}[ht]
\bce{\small
\begin{tabular}{lll}
(a) BD of steady branches&(b) Hopf orbits&
(c) Stability \\
\hs{-5mm}\raisebox{0mm}{\begin{tabular}{l}
\ig[width=0.25\tew,height=40mm]{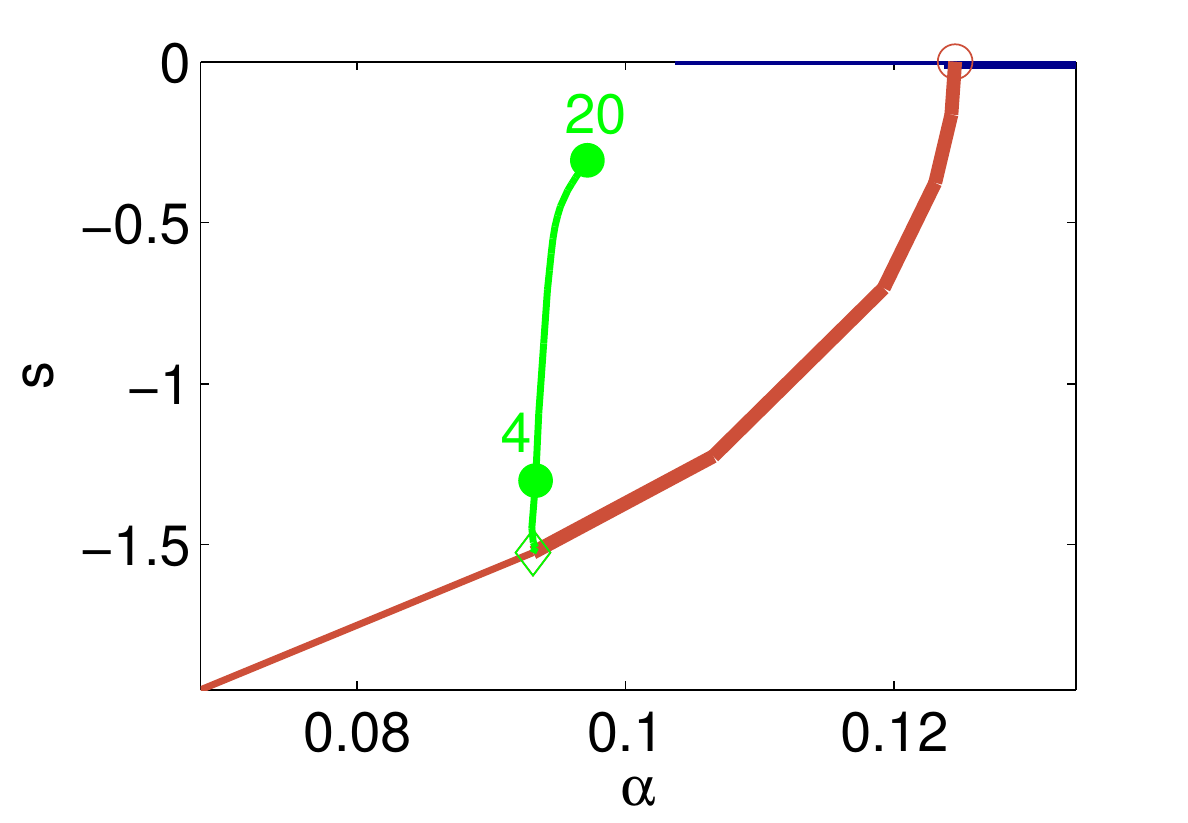}\\
\ig[width=0.25\tew,height=20mm]{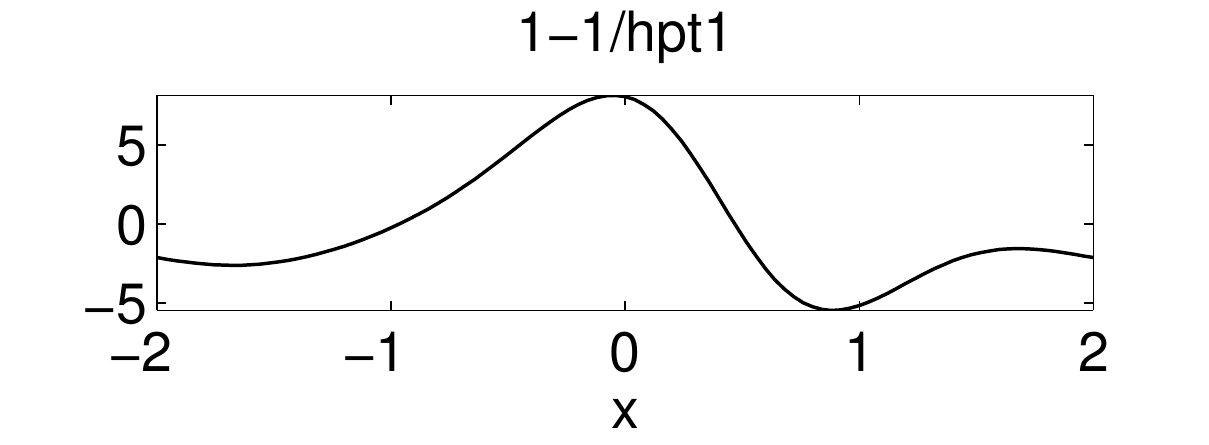}
\end{tabular}}
&
\begin{tabular}{l}
\hs{-12mm}\ig[width=0.28\tew,height=50mm]{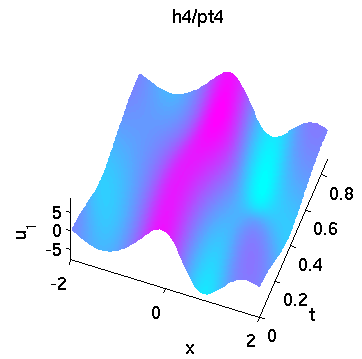}\hs{-2mm}
\ig[width=0.28\tew,height=50mm]{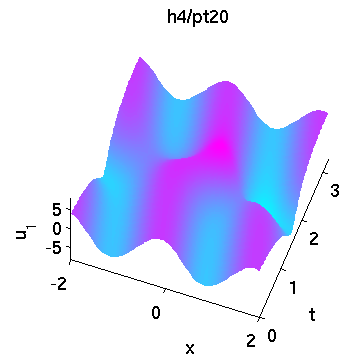}
\end{tabular}
&
\hs{-5mm}\raisebox{0mm}{\begin{tabular}{l}
\ig[width=0.18\tew]{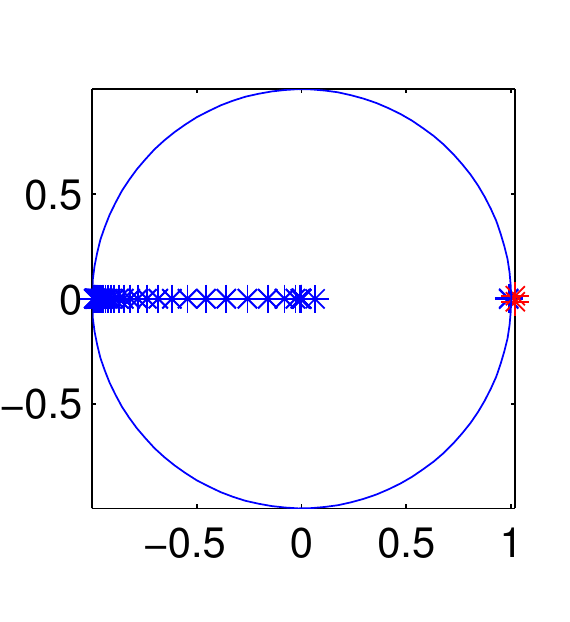}\\\ig[width=0.18\tew]{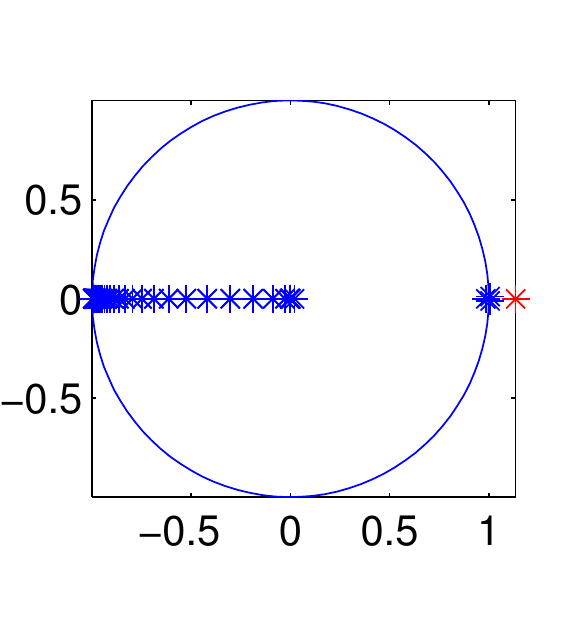}
\end{tabular}}
\end{tabular}
}\ece 

\vs{-5mm}
   \caption{{\small Results from {\tt kspbc4/cmds2.m}. (a) Bifurcation diagram 
($s$ over $\al$) 
of the traveling wave branch from Fig.~\ref{ksf1}(a), and of 
the first bifurcating modulated traveling wave branch (green), 
with the solution profile at bifurcation at the bottom. 
(b) shows two Hopf orbits on the green branch (in the frames moving 
with speeds $s$ from (a), respectively), and (c) the 
associated Floquet multipliers. The bifurcation is subcritical, 
and the modulated traveling waves are (mildly) unstable.  \label{ksf2}}}
\end{figure}

\hulst{
label=ksl1,language=matlab,stepnumber=50, firstnumber=1}{\hdhome/kspbc4/oosetfemops.m}

\hulst{caption={{\small {\tt \dname/oosetfemops.m} and {\tt \dname/sG.m}. The mass matrix {\tt p.mat.M} is redefined in {\tt oosetfemops} to $M^2$, and the proper mass matrix is stored in {\tt p.mat.M0}.}},
label=ksl3,language=matlab,stepnumber=5, firstnumber=1}{\hdhome/kspbc4/sG.m}

\hulst{
label=ksl4,language=matlab,stepnumber=5, firstnumber=1}{\hdhome/kspbc4/qf.m}

\hulst{caption={{\small  {\tt \dname/qf.m} and {\tt \dname/qfh.m}. 
The phase conditions for the steady and for the Hopf case. }},
label=ksl5,language=matlab,stepnumber=5, firstnumber=1}{\hdhome/kspbc4/qfh.m}

\subsection{Period doubling of a breather (demo {\tt symtut/breathe})}\label{pdsec}
\def\dname{symtut/breathe}\def\dhome{./breathe}
In \cite[\S4.2]{symtut} we studied the RD system 
\huga{\label{iim}
\pa_t u=\pa_x^2 u+f(u,v), \quad \pa_t v=D\pa_x^2 v+g(u,v), 
}
with homogeneous Neumann BC, 
$f(u,v)=u(u-\al)(\beta-u)-v$, $g(u,v)=\del(u-\ga v)$, with 
$\al,\beta,\ga>0$, and $0{<}\del{\ll}1$. For suitable parameters, this model has standing and traveling 
pulses (and traveling fronts), and for $\del{\ra}0$ we find a Hopf 
bifurcation to breathers, see the bottom row of Fig.~\ref{breathefig1} 
for examples. In \cite{symtut} 
this served as  an example for the usefulness of constraints, here 
regarding the approximate translational invariance for the case 
of narrow breathers. It turns out that the 
'primary breather branch' (red in (a)) looses 
stability in a period doubling bifurcation, yielding the magenta branch 
in (a), which starts out stable, and then looses stability in a torus 
\mbox{bifurcation, see (b).} 

\begin{figure}[ht]
\bce 
\begin{tabular}{p{0.48\tew}l} 
{\small (a) steady states (black), breather ({\tt h1}, red) and PD ({\tt pd1}, magenta) branch}&{\small (b) Multipliers} \\
\ig[width=0.22\tew]{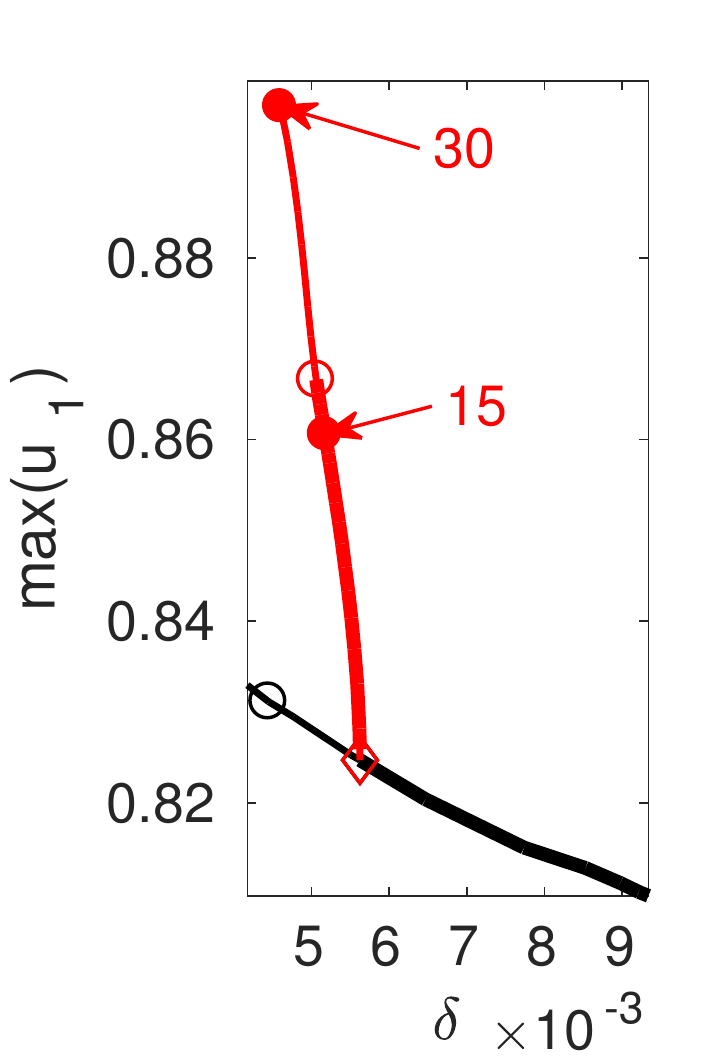}
\hs{-2mm}\ig[width=0.25\tew]{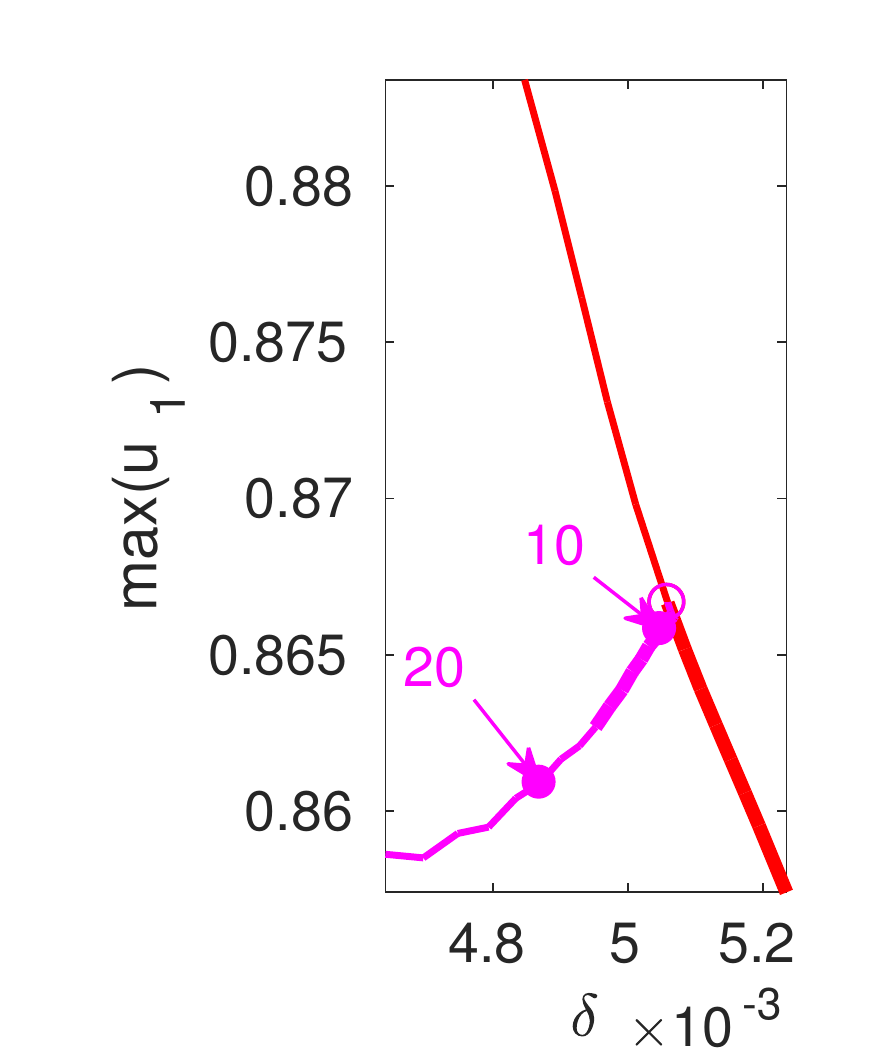}&
\hs{-0mm}\raisebox{28mm}{\begin{tabular}{l}
\ig[width=0.2\tew]{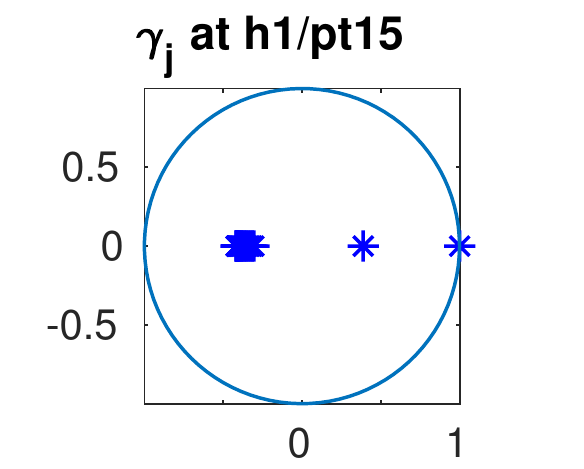}\\
\ig[width=0.18\tew]{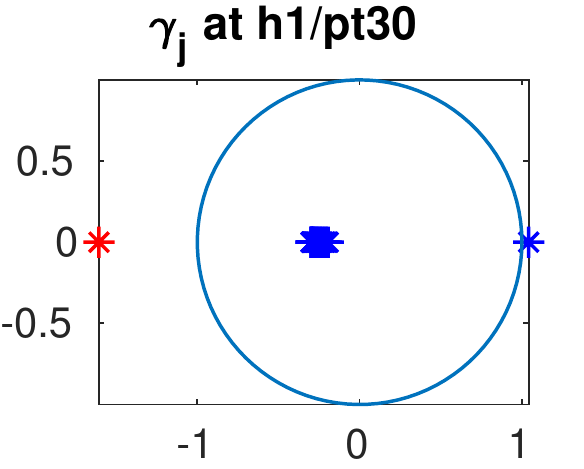}
\end{tabular}
\hs{-4mm}
\begin{tabular}{l}\ig[width=0.2\tew]{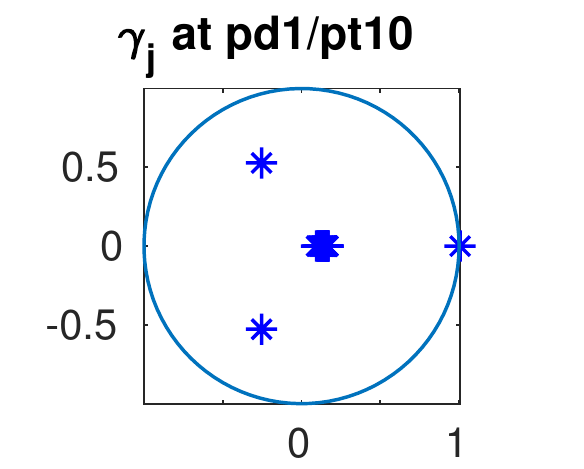}\\
\ig[width=0.18\tew]{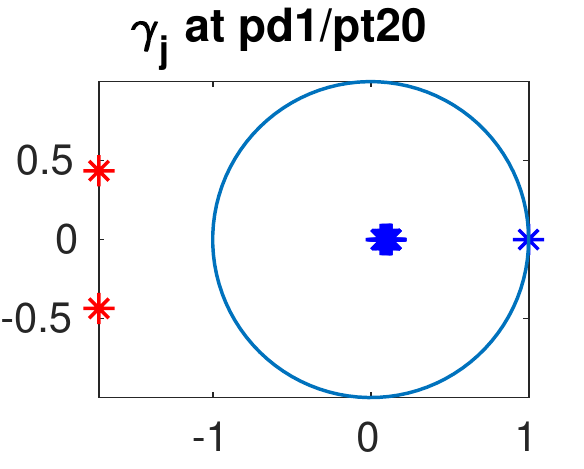}
\end{tabular}}\\
(c) Solution plots\\
\ig[width=0.22\tew]{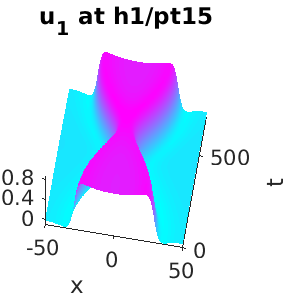}
\ig[width=0.22\tew]{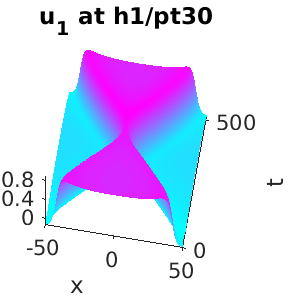}&
\ig[width=0.22\tew]{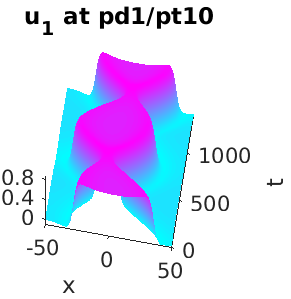}
\ig[width=0.22\tew]{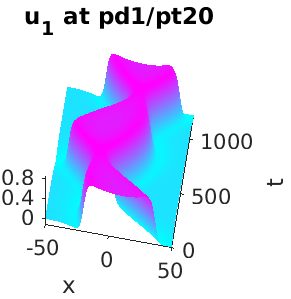}
\end{tabular}
\ece
\vs{-5mm}
   \caption{{\small Period doubling bifurcation in \reff{iim}, $(\al,\beta,\ga)=
(0.11,1,6)$, $D=2$.}
  \label{breathefig1}}
\end{figure}

\hulst{caption={{\small {\tt symtut/breathe/cmds1.m}, 
commands for the period doubling branch. }},
label=brel1,language=matlab,stepnumber=5, linerange=32-36}{\dhome/cmds1.m}

Figure \ref{breathefig1} is computed in the script {\tt symtut/breathe/cmds1.m}, see Listing \ref{brel1} for the relevant code snippet. 
For more background on the demo {\tt symtut/breathe} we refer to \cite{symtut}, 
and here only remark that: 
\bci 
\item A good localization of the PD bifurcation point on the breather 
is crucial; here $\ga_{{\rm crit}}\approx -1.1$ obtained using {\tt p.hopf.bisec=5} 
is good enough if we allow large residuals at startup of the magenta branch. 
After 5 steps we set {\tt p.nc.tol=1e-6} again. 
\item For the magenta branch we switch off the translational constraints, 
i.e., set {\tt p.hopf.nqh=0}.  While the constraint is useful 
for narrow breathers (at the start of the red branch), the wider 
breathers interact strongly enough with the boundary and the 
constraint can be dropped.   
The magenta branch can also be computed with {\tt p.hopf.nqh=1} but this becomes more 
expensive. Most importantly, due to a poorly localized BP 
(critical multiplier $-1.12$) we start with a very large 
tolerance {\tt tol=0.5} 
to get onto the period-doubled branch, but we can subsequently 
decrease the tolerance to 1e-8 as usual. 
\eci

\section{O(2) equivariance: traveling vs standing waves, and relative periodic orbits}\label{o2sec}
In \S\ref{cglsec} we considered the cGL equation over domains which lead to simple HBPs, i.e., boxes with NBC or DBC, where 
moreover in 2D and 3D we chose suitable side-lengths $l_x,l_y, l_z$, 
in particular $l_x\ne l_y$. If  for 
instance in 2D we instead chose a square domain, then naturally the 2nd HBP would be double, with oscillating 'horizontal' and 'vertical' stripes as 
two Hopf eigenfunctions. 

For steady bifurcations, the higher multiplicities 
of BPs due to discrete symmetries and the associated multiple 
bifurcating  branches can be dealt with systematically in \pdep, as 
described in \cite{pftut}. 
For HBPs of higher multiplicity we do not 
yet provide similar routines, but rather treat them in an ad hoc way. 
Moreover, for Hopf problems multiple branches due to continuous 
symmetries are probably even more 
important than multiple branches due to discrete symmetries. 
In particular, $\rO(2)$ equivariant Hopf bifurcations arise in a variety 
of settings, for instance for translational 
invariant problems (due to pBC) with reflection symmetry, and similarly 
for problems on circular domains, 
where the role of translational invariance is played by spatial 
rotations. 
Thus, here we first consider the cGL equation in boxes with pBC, 
and in a disk domain (demos {\tt cglpbc} and {\tt cgldisk}), 
and then review the demo {\tt gksspirals} dealing with a RD model from \cite[\S3.2]{hotheo} in the unit disk.  

For background on $\rO(2)$ equivariant Hopf bifurcation see, e.g., \cite{GoS2002}, 
and the references therein. Loosely said, the main result is that for 
double HBPs we generically obtain three bifurcating branches of Hopf orbits: 
left/right traveling waves (TWs), and standing waves (SWs), which correspond 
to equal amplitude superpositions of TWs. Importantly, the TWs are 
steady solutions in an appropriate co-moving frame, and are thus much cheaper 
to compute than general Hopf orbits. 

\subsection{The cGL equation in boxes with pBC: demo {\tt cglpbc}}
\label{pcglsec}
\def\dname{cglpbc}\def\dhome{./hopfdemos/cglpbc}
\subsubsection{1D}\label{pcgl1dsec}. 
We consider a variant of \reff{cAC}, namely 
\hual{\label{cAC2} 
\pa_t& \bpm u_1\\ u_2\epm =-G(u,\lam)\\
&:=\bpm \Delta+r&-\nu\\\del^2\nu&\Delta+r\epm
\bpm u_1\\ u_2\epm-(u_1^2+u_2^2)\bpm c_3 u_1-\mu u_2\\ 
\mu u_1+c_3 u_2\epm-c_5(u_1^2+u_2^2)^2\bpm u_1\\ u_2\epm+s\pa_x\bpm u_1\\u_2\epm, \notag
}
first on the interval $\Om=(-\pi,\pi)$ with pBC. The additional parameter 
$\del$ can be used to break the phase invariance $\ds u\mapsto \bpm \cos\phi&-\sin\phi\\\sin\phi&\cos\phi\epm u$ (i.e.~$u\mapsto \er^{\ri \phi}u$ 
in complex notation) of \reff{cAC0}, see below for further comments, 
and the parameter $s\in\R$ describes a frame moving with speed $s$, useful later for the continuation of TWs. 
As in \S\ref{cglsec} we fix $c_3=-1, c_5=1$, $\nu=1$, 
(and initially  $\mu=0.5$, $\del=1$ and $s=0$), and use $r$ as the primary bifurcation 
parameter. 

The HBPs from the trivial solution are still 
\huga{\label{hp17}
\text{$r_k=k^2, k\in\N$,\ \ eigenvalues $\pm \ri\om$ with $\om=\nu$, }
}
but now are double for $k>0$. Two 'natural' two eigenfunctions are 
(in complex notation) 
\huga{\label{aefu} 
\phi_1(t,x)=\er^{\ri(\om t-kx)}\text{ and } \phi_2(t,x)=\er^{\ri(\om t+kx)}, 
}
and the ansatz for bifurcating periodic orbits is 
\huga{\label{bans}
u=z_1\phi_1+z_2\phi_2+\hot, \quad (z_1,z_2)\in\C^2. 
} 
Thus, $\phi_1$ corresponds to a right TW with speed $\om/k$, and 
$\phi_2$ to a left TW with speed $-\om/k$. For $s=0$, \reff{cAC2} 
is $\rO(2)$ equivariant, i.e., $G(\ga u)=\ga G(u)$ for all $\ga\in\Ga=\rO(2)$. 
Here $\ga=(m,\xi)$, $m\in\Z_2=\{\pm 1\}$, $\xi\in \sO(2)=[0,2\pi)$, and 
the action of $\ga$ on $x$ and $u(x)$ is given 
by $(\ga u)(x)=u(m(x+\xi))$ (reflection and translation). 
For each $k\in\N$, the subspace $X_k:={\rm span}\{\phi_1,\phi_2\}$ 
is $\Ga$ invariant, and the action of $\ga$ on $(z_1,z_2)$ is 
\huga{
m:(z_1,z_2)\mapsto (z_2,z_1), \quad \xi:(z_1,z_2)\mapsto (\er^{\ri \xi}z_1,\er^{-\ri\xi}z_2).
}
The equivariant Hopf theorem (\cite[Thm 4.9]{GoS2002} or \cite[Thm 4.6]{hoyle})  yields that generically for each $k$ we have exactly three bifurcating branches, 
namely 
\huga{\label{eHo}
\left\{\barr{ll}
\text{rTW}&u(x,t)=z_1\phi_1+\hot\\
\text{lTW}&u(x,t)=z_2\phi_2+\hot\\
\text{SWs}&u(x,t)=z(\phi_1+\phi_2)+\hot, 
\earr\right. 
}
where as usual $\hot$ stands for higher order terms. Moreover, 
the TWs are solutions of the form 
\huga{\label{eHo2} \text{$u(x,t)=v(x-st)$ with some speed $s\in\R$, 
}
}
i.e.~{\em relative equilibria}, which means that we can find them 
as steady solutions of \reff{cAC2} with a suitable $s$. 
In fact, from the phase invariance for $\del=1$ 
we already have an explicit formula for TWs, 
namely 
\hual{\label{cACs}\text{$u(x,t)=R\er^{\ri (\om t-k\cdot x)}$, }\ 
|R|^2{=}-\frac{c_3}{2c_5}\pm\sqrt{\frac{c_3^2}{4c_5^2}+r-k^2}, 
\ \ \om{=}\om(k,r){=}\nu-\mu|R|^2, 
}
but \reff{eHo2} is the more general result. 
Moreover,  on the spaces $X_k$ the additional SO(2) phase symmetry 
$\vt: u\mapsto \er^{\ri \vt} u$ acts like time-shifts $u\mapsto u(t+\vt/\tau)$, and hence is not an additional symmetry for the Hopf bifurcation 
and does not need to be considered further.

\brem\label{hmrem}{\rm In summary, at each $r_k=k^2$ we have the bifurcation of TWs and SWs, 
and our aim is to compute these numerically (even if for $\del=1$ we know the TWs analytically from \reff{cACs}). Thus we face a similar problem 
like for steady bifurcations of higher multiplicity, discussed in detail 
in \cite{pftut}: When computing the eigenfunctions associated to the eigenvalue 
$\ri \om$ of $G_u$ we in general do not obtain the 'natural' ones 
$\phi_{1,2}|_{t=0}=\er^{\pm\ri k x}$ from \reff{aefu}, but {\em some} 
linear independent 
$\tilde\phi_{1,2}\in\spani\{\er^{\ri kx},\er^{-\ri kx}\}$. Thus, even though 
we know the analytic form \reff{eHo} of the bifurcating branches, 
this only applies to the natural basis $\phi_{1,2}$ of the center eigenspace. 
In principle we could compute 
the (3rd order) amplitude system on the center manifold associated to the 
basis 
$\tilde\phi_{1}\er^{\ri \om t}, \tilde\phi_{2}\er^{\ri \om t}$, and from this the coefficients $\tilde z_{1,2}$ 
for TWs and SWs. However, in contrast to the steady case, for which 
we provide routines to do so, we refrain from implementing this 
for the Hopf case in \pdep, 
because 
the general case of multiple Hopf bifurcations becomes significantly 
more complicated, see, e.g., \cite{kiel79}, or \cite{mei2000} 
for the case where additional mode interactions with steady modes come into  play. Instead, we proceed ad hoc, and require user input of coefficients 
$z_1, z_2$ (and $z_3,\ldots,z_m$ in case of still higher multiplicity $m$, 
see \S\ref{cgl2sec}). In practice this works quite well. 
}\eex\erem 

\hulst{caption={{\small {\tt cglpbc/cmds1d.m} (first 2 cells). 
In C1, the only new command is {\tt box2per}, which switches on the pBC. In C2 
we compute a TW branch and a SW branch bifurcating from Hopf point 2 by 
'guessing' and passing on to {\tt hoswibra} coefficients {\tt aux.z}. 
For the SW branch we additionally set the average 
translational PC {\tt qfh} as in \S\ref{kssec}, with reference 
profile $u_0(x)=u_{{\rm pred}}(0,x)$, where $u_{{\rm pred}}$ is the 
predictor for the Hopf orbit, and 
speed parameter $s$ (given in {\tt par(6)}). }},
label=cll1,language=matlab,stepnumber=11, linerange=2-19}{\dhome/cmds1d.m}

In Listing \ref{cll1} we give the start of {\tt cglpbc/cmds1d.m} for \reff{cAC2} on 
$\Om=(-\pi,\pi)$ with pBCs, which are switched on via {\tt box2per}, 
see \cite{pbctut}. We ignore the first (spatially homogeneous $k=0$) 
Hopf branch, and in C2 compute one TW branch and one SW branch 
bifurcating at the 2nd HBP, 
corresponding to wave number $k=1$. Here we 'guess' by some trial and 
error the coefficients $z_1,z_2$ for each of these branches.  Additionally, 
for the SW branch we set a translational PC (with $s=0$). It turns out that this PC is 
usually enough to force SWs, even if the guess for the coefficients {\tt aux.z} rather corresponds to a TW. See Fig.~\ref{cfig1} for some results. 
The SW branch stays unstable up to $r=2$ and beyond. 
The TW branch starts unstably (as expected, as the trivial branch is already unstable) with $\ind(u_H)=5$, which turns into $\ind(u_H)=4$ at the fold, into $\ind(u_H)=2$ 
shortly after, and $u_H$ becomes stable near $r\approx 1.25$. 
Both crossings of unstable multipliers into the unit circle 
are of torus type, and hence the bifurcating branches 
in this form currently can not be computed with \pdep. But as already said, the TWs can also be computed as relative equilibria, 
i.e., as steady states in a frame comoving with speed 
\huga{s=\om/k  \text{ at bifurcation,}} 
where $k$ is the spatial wave number. The pertinent branch switching is implemented in 
$$\text{
{\tt p=twswibra(dir,fname,spar,kwnr,newdir,aux)},}
$$ where {\tt spar} is 
the index of $s$ in the parameter vector, ${\tt kwnr}=k$, and 
{\tt aux.z} again can be used to pass the coefficients $z_{1,2}$ 
for the predictor guess. Complementing this with the 
PC {\tt qf}, i.e., $\spr{\pa_x u_0,u}=0$, 
we obtain the same TW as in C2 with a small error in the period $T$ 
between the two methods, which vanishes if we increase the temporal 
resolution for the {\tt hoswibra} solution. To obtain a 
space--time plot of TWs, use {\tt twplot}.

\begin{figure}[ht]
\bce 
\begin{tabular}{lll} 
{\small (a) }&{\small (b) }&{\small (c) } \\
\ig[width=0.17\textwidth]{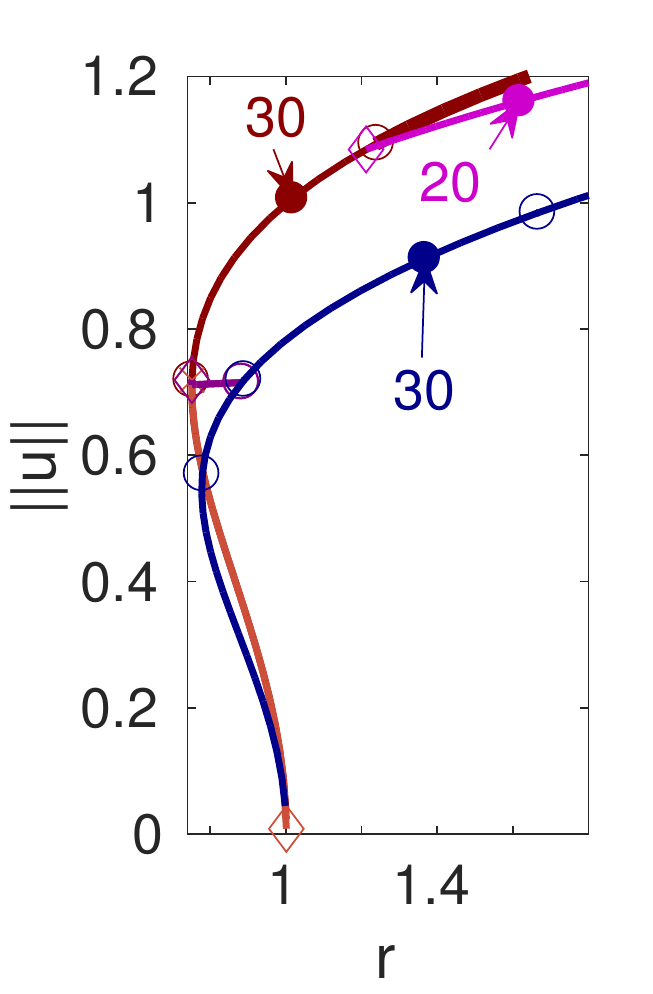}
\hs{-2mm}\ig[width=0.16\textwidth]{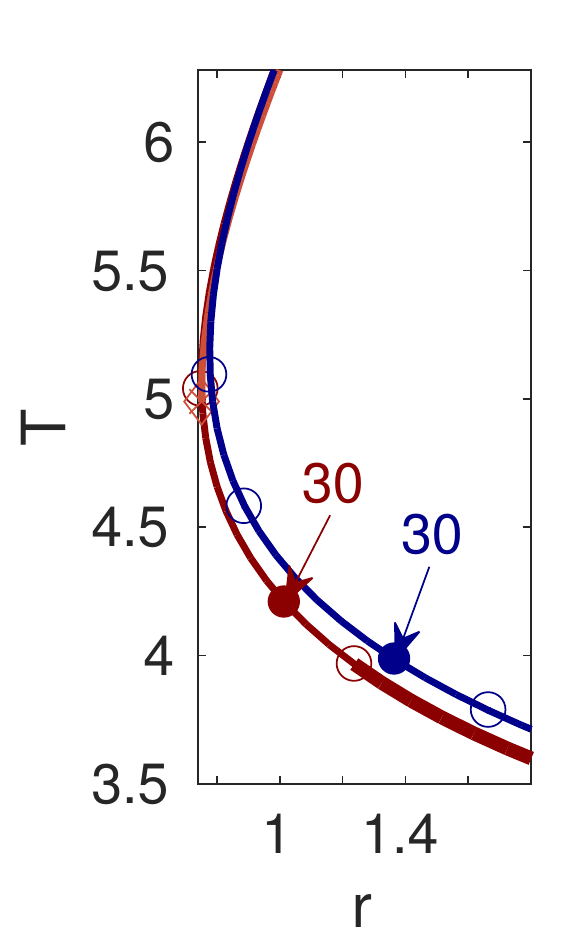}&
\hs{-3mm}\raisebox{23mm}{\begin{tabular}{ll}
\ig[width=0.16\textwidth]{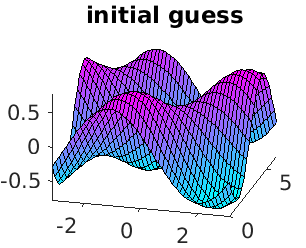}&
\hs{-3mm}\ig[width=0.17\textwidth]{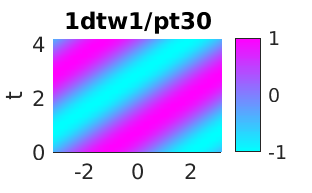}\\
\ig[width=0.16\textwidth]{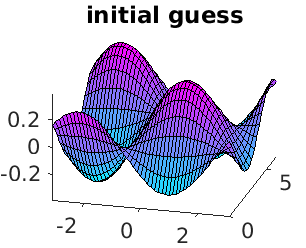}&
\hs{-3mm}\ig[width=0.17\textwidth]{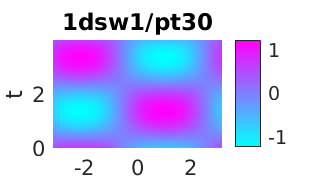}
\end{tabular}}&
\hs{-7mm}\raisebox{0mm}{\ig[width=0.16\textwidth]{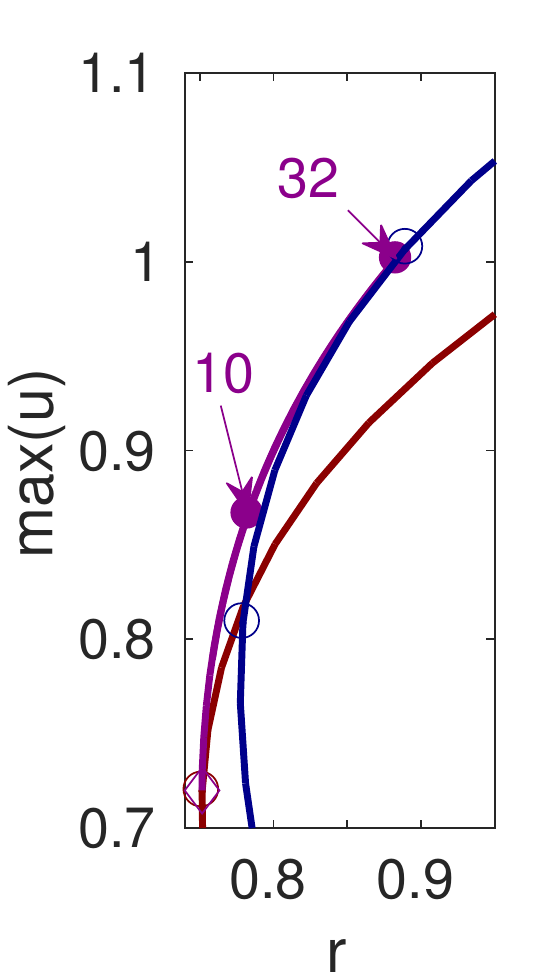}
\hs{-3mm}\ig[width=0.16\textwidth]{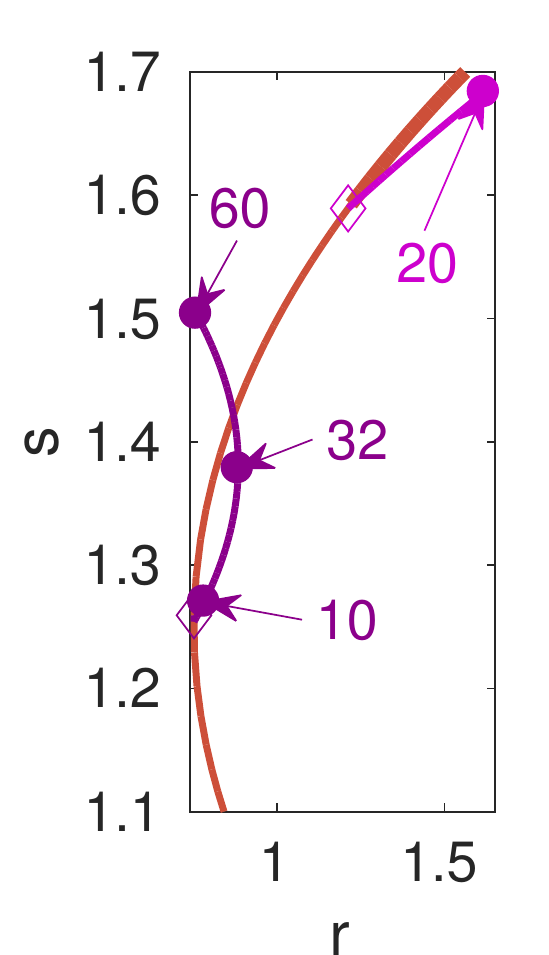}}
\end{tabular}\\
\begin{tabular}{ll} 
{\small (d) }&{\small (e) }\\
\hs{-2mm}\raisebox{23mm}{\begin{tabular}{l}
\ig[width=0.18\textwidth]{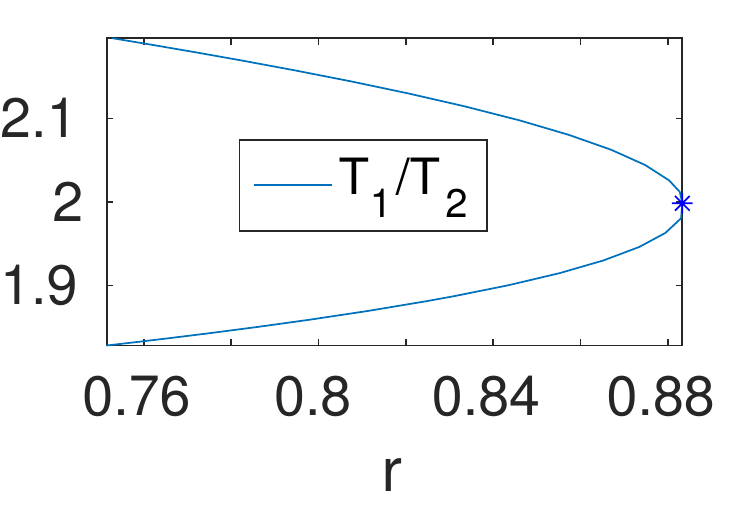}\ig[width=0.2\textwidth]{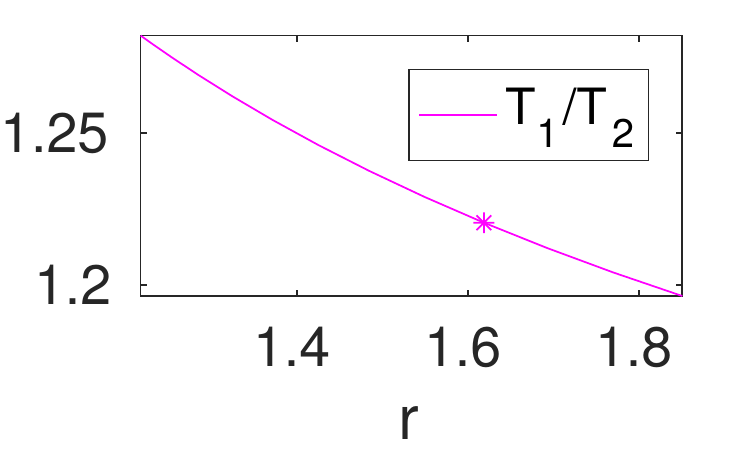}\\
\hs{-5mm}\ig[width=0.42\textwidth]{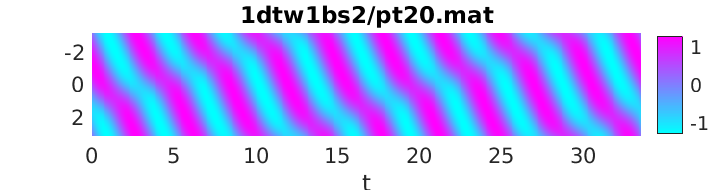}
\end{tabular}}&
\hs{-6mm}
\raisebox{28mm}{\begin{tabular}{lll}
\ig[width=0.18\textwidth]{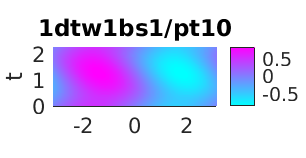}&
\hs{-2mm}\ig[width=0.18\textwidth]{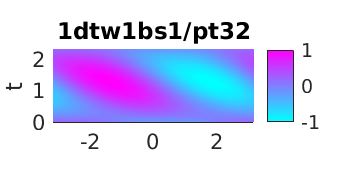}&
\hs{-2mm}\ig[width=0.18\textwidth]{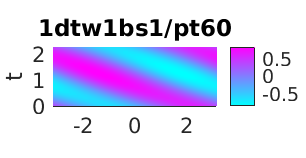}
\\
\ig[width=0.18\textwidth]{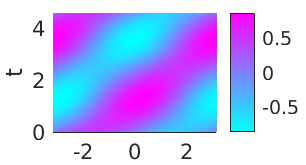}&
\hs{-2mm}\ig[width=0.18\textwidth]{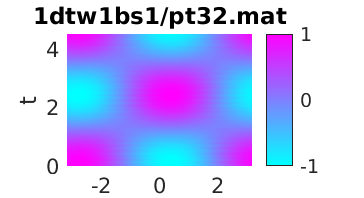}&
\hs{-2mm}\ig[width=0.18\textwidth]{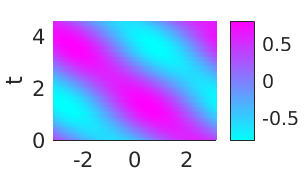}
\end{tabular}}
\end{tabular}
\ece
\vs{-8mm}
   \caption{{\small \reff{cAC2} on $\Om=(-\pi,\pi)$ with pBC, $(\nu,\mu,c_3,c_5,\del)=(1,0.5,-1,1,1)$. (a) BD of TW (brown) and SW (blue) 
branches, and secondary bifurcation from TW branch (dark and light magenta). 
(b) intial guesses for TW and SW branches, and two solution plots. 
(c) Zoom into BDs near TW fold, including secondary branches. 
(d) Quotients (top) 
of periods $T_1=2\pi/s$ and $T_2$ where $s$ is the frame speed 
and $T_2$ the period in the comoving frame. Left 
on 'connecting branch', right on magenta 
modulated TW. Bottom: solution plot in lab frame (bottom) at marked magenta 
point. 
(e)  solution plots at marked points on dark magenta branch in comoving frame (top) and lab frame (bottom). At the fold (pt32) there is the resonance $T_1=2T_2$, and the combination of left traveling in the moving frame 
and the motion of the frame yields the SW with period $T=2T_2$. 
\label{cfig1}}}
\end{figure}

\hulst{caption={{\small {\tt cglpbc/cmds1d.m} (cells 3 and 4). 
In C3 we compute 
the TW branch as a relative equilibrium via {\tt twswibra}. On this 
branch we find HBPs, and in C4 we compute secondary Hopf branches 
bifurcating from this relative equilibrium; see Fig.~\ref{cfig1} for 
plots (as obtained from {\tt plotcmds.m}, and text for further comments.  
}},
label=cll1,language=matlab,stepnumber=50, linerange=20-30}{\dhome/cmds1d.m}

The continuation of the TW as a relative equilibrium yields HBPs 
on this branch, 
at the locations where the continuation as periodic orbits yields the 
(torus) BPs, see also Remark \ref{tw1srem}(a). 
 Now we can use {\tt hoswibra} to compute 
the bifurcating modulated TWs as {\em relative periodic orbits} (relPO). 
Figure \ref{cfig1}(c) shows a 
zoom near the fold; the bifurcating branch seems to connect to the 
SW branch near $r\approx 0.89$, see also the solution plots in (e), 
where the bottom row shows the solutions in the lab frame, i.e., 
by shifting back $x\mapsto x+st$. Here we plot over $\Om\times [0,2T_2)$, 
where $T_2\approx 2.28$ is the period in the moving frame. In general, Hopf orbits 
bifurcating from relative equilibria (i.e., in the moving frame) 
correspond to quasiperiodic solutions in the lab frame, and the 
quotient $T_1/T_2$ (with $T_1=L/(ks)$, $s$ the comoving speed, $L$ 
the domain size, and $k$ the wave number) 
varies continuously. At pt32 
the solution is (approximately) $2T_1$ periodic in the lab frame 
and corresponds to the SW. 
Similarly, solution 20 on the magenta branch 
is approximately $11T_2=33.36$--periodic in the lab frame, see the bottom 
panel of (d).  

\brem\label{tw1srem}{\rm  
(a) The multipliers of the periodic orbit $u(x,t)$ in the lab frame 
are given by $\ga_j=\er^{-\mu_j T}$,  
where the $\mu_j$ are the eigenvalues of the linearization around the 
TW in the comoving frame. In detail, the ansatz $u(x,t)=v(x-st,t)$ yields 
$\pa_t v=-G(v)+s\pa_\xi v$, with linearization $\pa_t v=-G_u(v_0(\xi))v+s\pa\xi v=:-Lv$. As $L$ is independent of $t$, the linear flow yields 
$v(T)=\sum_j c_j\er^{-\mu_j T}\phi_j$, where $v(0)=\sum_j c_j\phi_j$, and 
where for simplicity we assumed 
semisimple eigenvalues $\mu_j$ of $L$ with associated eigenvectors (eigenfunctions) 
$\phi_j$, $j=1,\ldots,n_u$. \\
(b) To plot the modulated TWs in the lab frame we use 
the function 
{\tt lframeplot(dir,pt, wnr,cmp,aux)}. This first aims to determine the minimal 
$m\in\N$ such that $msT=qL$ for some (minimal) $q\in\N$, 
where $T$ is the period in the frame moving with speed $s$ 
and $L$ is the domain size. 
Equivalently, $m=qL/sT$ for some (minimal) integer $q$, and the minimal is 
period $T^*=mT=\frac{q}{s}L$. Of course, $m=qL/sT\in\N$ numerically means $|m-\lfloor m\rfloor|<$tol, where tol can be passed as {\tt aux.pertol}. 
This should not be taken too small, i.e., on the order of the 
expected error in the speed $s$ and the period $T$. Naturally, this also 
ignores the fact that generically $T_1/T_2\not\in\Q$, and that 
the associated orbits are quasiperiodic rather than periodic (with a 
possibly large period). Alternatively, 
an integer $m$ can be passed in {\tt aux.m} to force the plot over 
$[0,mT]$. {\tt lframeplot} at the end also reports 
the final (integer) grid-point shift of the transformation to the moving 
frame is given, which should be $0$. 
}\eex\erem 

\subsubsection{2D box with pBC in $x$}\label{cgl2sec} 
In Fig.~\ref{cfig2b} we give some 
very introductory results for \reff{cAC2} over the 2D square box $\Om=(-\pi,\pi)^2$, with pBC in $x$ and homogeneous Neumann BC in $y$. The 
2nd HBP at (analytically) $r=1$ is then triple, with Hopf eigenfunctions 
(in complex notation)  
\huga{\label{3ker}
\er^{\ri(\om t-x)}, \er^{\ri(\om t-x)}, \cos(y)\er^{\ri\om t} 
}
and modulo spatial translation we may expect at least five 
primary bifurcating branches: SW and TW (twice) in $x$, 
SW in $y$, and a mixed SW mode of the form $b(t)\sin(x)\cos(y)$. Four such 
branches are computed in {\tt cmds2d.m} via 'educated' guesses of the 
three coefficients for the three numerical eigenfunctions replacing 
\reff{3ker}. Naturally, the TW branches 
can again also be computed as relative equilibria, and we find several 
secondary bifurcations. See {\tt cmds2d.m}, but here we refrain from 
giving the further details. 

\begin{figure}[ht]
\bce 
\begin{tabular}{ll} 
{\small (a) }&{\small (b) } \\
\ig[width=0.18\tew]{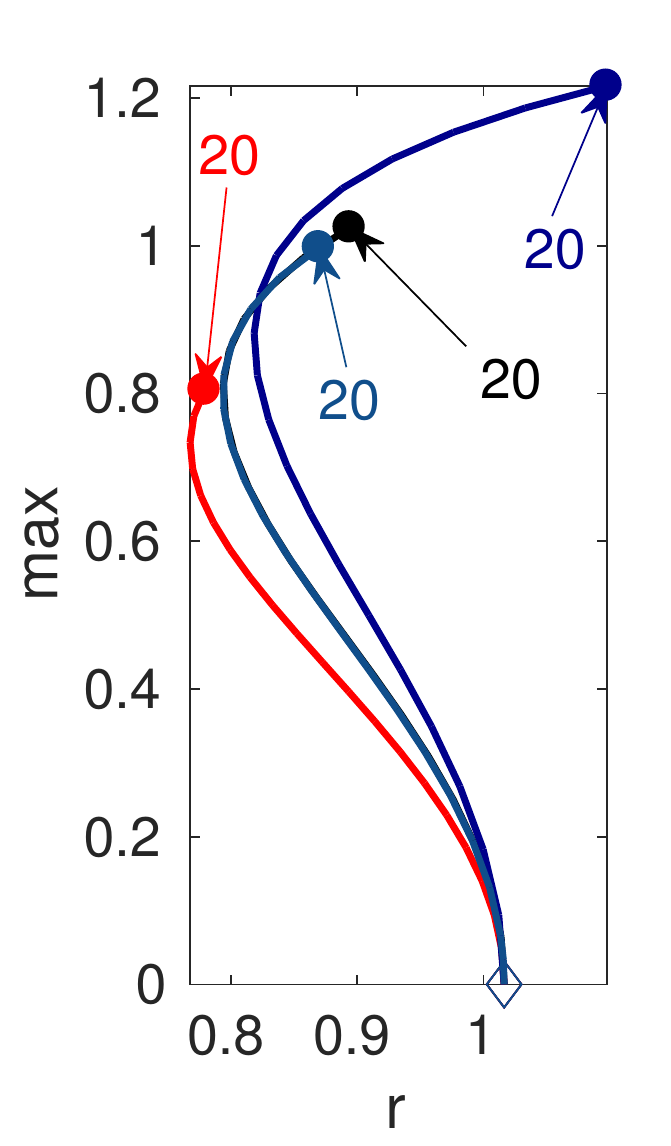}
\hs{-1mm}\ig[width=0.175\tew]{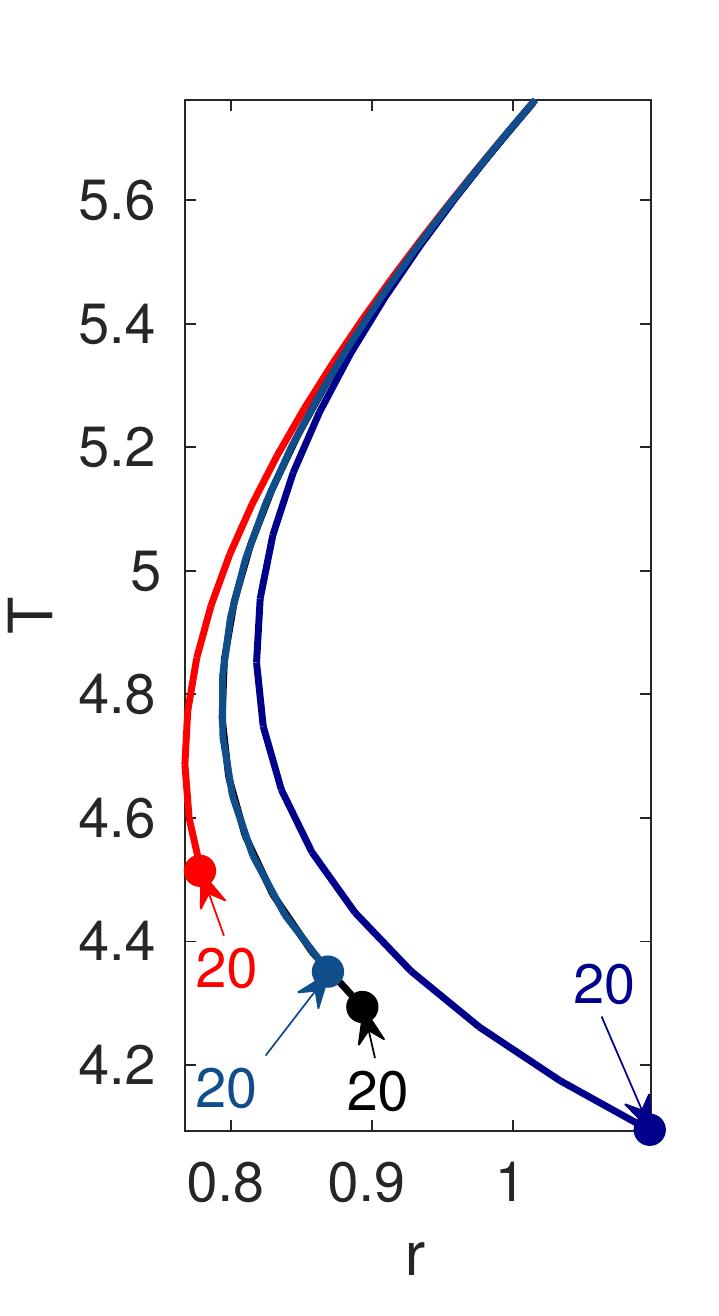}&
\hs{-10mm}\raisebox{28mm}{\begin{tabular}{l}
\ig[width=0.55\tew]{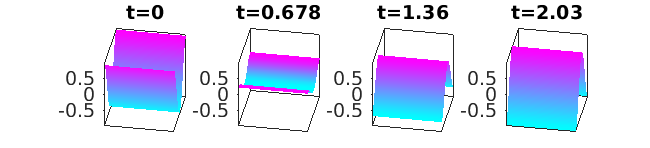}\\
\ig[width=0.55\tew]{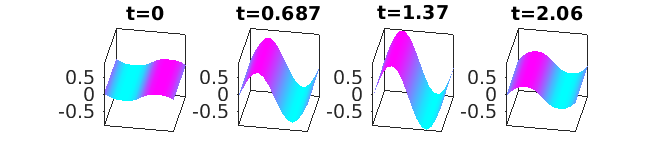}
\end{tabular}}\\
(c)\\
\hs{-8mm}\ig[width=0.55\tew]{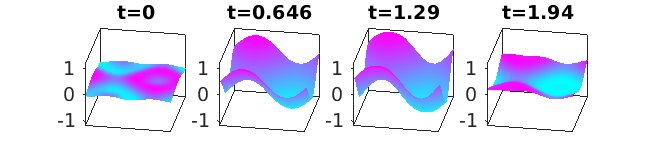}&
\hs{-10mm}\ig[width=0.55\tew]{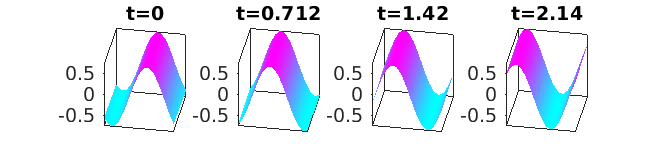}
\end{tabular}
\ece
\vs{-5mm}
   \caption{{\small \reff{cAC2} on $\Om=(-\pi,\pi)^2$ with pBC in $x$ and 
Neumann BC in $y$, $(\nu,\mu,c_3,c_5,\del)=(1,0.5,-1,1,1)$. (a) BD of 
4 branches bifurcating from the 2nd HBP $r=1$: swy (black), tw (red), swx-y (dark blue), swx (light blue). (b,c) example solution plots (roughly half a 
period). swy (top) and (swx) in (b), swx-y (left) and tw (right) in (c). 
 }  \label{cfig2b}}
\end{figure}

\subsection{The cGL equation in a disk: demo {\tt cgldisk}}\label{cgldsec}
\def\dname{cgldisk}\def\dhome{./hopfdemos/cgldisk}
\def\Krot{\tt K_{{\rm rot}}}
A situation very similar to the 1D-pBC case arises for \reff{cAC2} in a disk 
with Neumann BC (or other rotationally invariant BCs).  
The symmetry group is again $\rO(2)$ where the role 
of spatial translations is now played by spatial rotations 
\huga{\bpm x\\ y\epm\mapsto 
R_\vt\bpm x\\y\epm:=\bpm \cos\vt&-\sin\vt\\\sin\vt&\cos\vt\epm\bpm x\\y\epm, 
}
$\vt\in [0,2\pi)$. The generator $\pa_\vt R_\vt|_{\vt=0}$ of the associated Lie algebra acts on $u(x,y)$ as 
$$
\pa_\vt u(R_\vt(x,y))|_{\vt=0}=\bigl[\pa_x u\pa_\vt x(\vt) 
+\pa_y u\pa_\vt y(\vt)\bigr]|_{\vt=0}=-y\pa_xu+x\pa_y u=:K_{{\rm rot}} u.
$$ 
Hence the rotational phase condition reads 
$\spr{K_{{\rm rot}} u_0, u}=0$, where $u_0$ is a suitable profile, 
typically set at bifurcation. Similarly, the rotating wave (RW) ansatz 
$u((x,y),t)= \uti(R_{-st}(x,y),t)$ yields 
$$
\pa_t u=-sK_{{\rm rot}}\uti+\pa_t \uti=-G(u)=-G(\uti), 
\text{ hence }\pa_t \uti=-G(\uti)+sK_{{\rm rot}}\uti, 
$$
after which we drop the\, $\widetilde{}$\, again. 

For the implementation, in {\tt oosetfemops} we generate $K_{{\rm rot}}$ 
via 
\def\Krot{K_{{\rm rot}}}
\huga{\label{krots}
\text{{\small {\tt po=getpte(p); x=po(1,:); y=po(2,:); p.mat.Krot=convection(fem,grid,[-y;x])}}.}
} 
The PC $q=\spr{\Krot u_0,u}=0$ is implemented as 
\huga{\label{pcr}
\text{{\small {\tt
function q=qf(p,u); q=(p.mat.Krot*p.u0)'*u(1:p.nu); end }}}
}
and the pertinent modification of {\tt sG} reads 
{\tt r=K*u-p.mat.M*f+s*Krot*u}. 

The eigenfunctions $v$ of $\Delta$ with Neumann BC 
have the form 
$u_{n,j}(x,y)=B_{n,j}(lr)g_n(\phi)$, where $g_n(\phi)=\er^{\ri n\phi}$, 
$l$ is a scaling factor, and $B_{n,j}$ is a Bessel function. 
These can be used to explicitly compute the HBPs from $u\equiv 0$, and 
to see that the HBPs are simple for $n=0$ and double for $n\ne 0$ 
(with $\sin(n\phi)$ and $\cos(n\phi)$ the two basis functions in the 
angular direction). 
For $n=0$, there is no angular dependence, and hence switching on the 
PC on such branches (with $u_0=u_0(r)$) 
leads to a singular Jacobian because  $K_{{\rm rot}}u_0=0$. 
For coarse meshes, the rotational invariance is sufficiently broken 
for the continuation to work also without PC, but for finer meshes 
the PC becomes vital for robust continuation. 

\taskip
\begin{table}[ht]\caption{Short overview of scripts and functions in {\tt hopfdemos/cgldisk}; see sources for details. 
\label{cgldtab}}
{\small 
\begin{tabular}{l|p{0.64\tew}}
script/function&purpose,remarks\\
\hline
cmds2d, plotcmds, cGLinit&main script, plotting commands, init function as usual\\
sG, sGjac, nodalf,njac&rhs, Jacobian, and nonlinearity, as usual\\
qf,qjac,qfh,qfhjac&phase conditions, based on {\tt Krot} generated in 
{\tt oosetfemops}, and versions for Hopf orbits\\
\hline
hoplott, plottip, lfplottip, gettip&some additional customized plot commands, and 
helper functions\\
\end{tabular}
}\end{table}
\teskip

Table \ref{cgldtab} gives a short overview of the files for the implementation. 
In {\tt cmds2d} we consider \reff{cAC2} in a disk with radius $\pi$, 
with base parameters $(\nu,\mu,c_3,c_5,\del)=(1,-5,-1,1,1)$, and 
bifurcation parameter $r$ as before. We changed $\mu$ in order to have 
pronounced 'spirals' \cite{KH81, bark90, Bark95, scheel98, SSW99} as RWs, 
see Fig.~\ref{cfig2}, and \S\ref{rotsec} for 
further more general comments. Relative periodic orbits then typically 
correspond to 'meandering spirals' which are 
most easily characterized by the motion of the tip. In order to compute the spiral tips with reasonable accuracy,  
at the start of {\tt cmds2d} we locally refine the mesh near $(x,y)=0$, 
see Fig.~\ref{cfig2}(a), 
leading to a mesh with about 1400 grid points. The temporal 
resolution for POs will be 30 gridpoints, and thus we will have 
about 84000 DoF, which we find a reasonable compromise between 
accuracy and speed; see also Remark \ref{cgldrem}.  
\def\wt{\tilde}

\begin{figure}[H]
\bce 
\begin{tabular}{lll} 
{\small (a) }&{\small (b) }&{\small (c)} \\
\hs{-2mm}\raisebox{2mm}{\ig[width=0.21\textwidth]{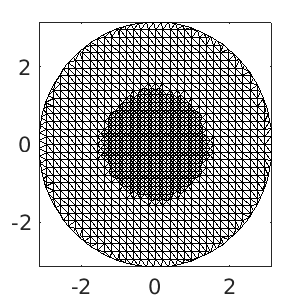}}
&\hs{-1mm}\ig[width=0.23\textwidth]{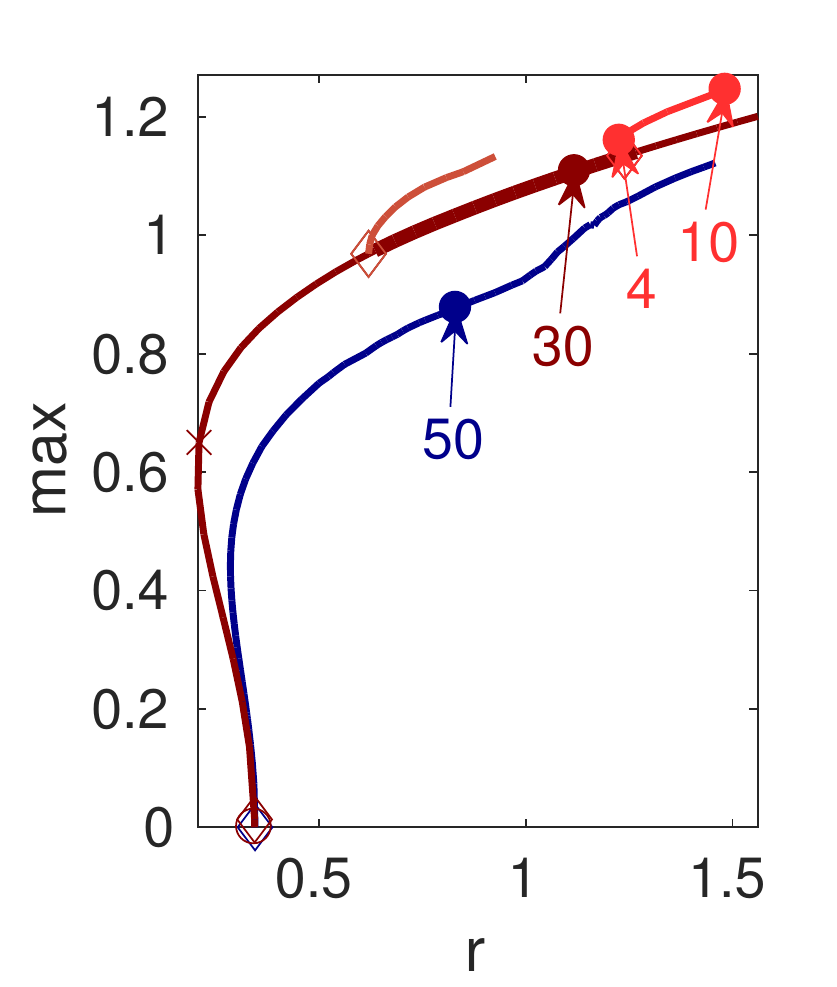}
&\hs{-6mm}\raisebox{20mm}{\begin{tabular}{l}
\ig[width=0.55\textwidth]{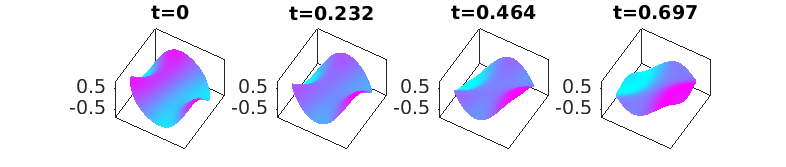}\\
\ig[width=0.5\textwidth]{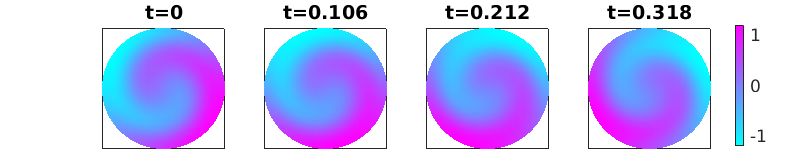}
\end{tabular}}
\end{tabular}
\begin{tabular}{llll} 
{\small (d) }&{\small (e) }&{\small (f)}&{\small (g)} \\
\hs{-5mm}
\raisebox{0mm}{\begin{tabular}{l}
\ig[width=0.22\textwidth]{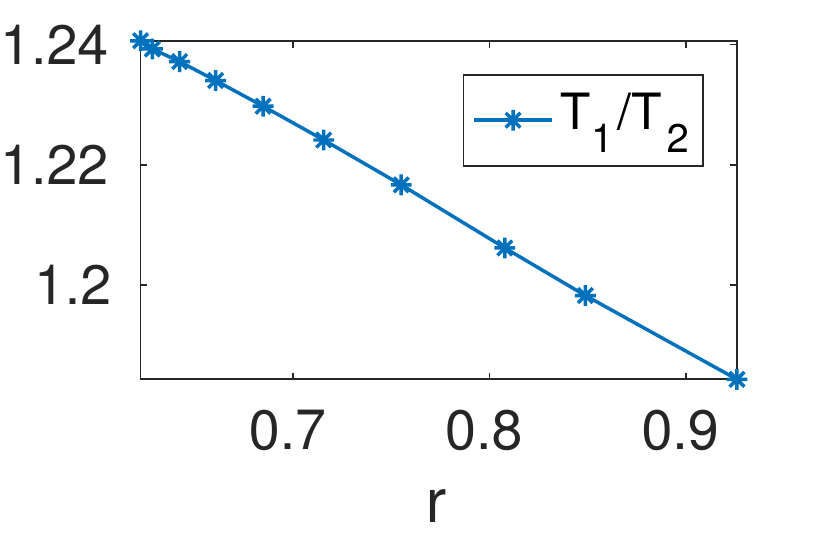}\\
\ig[width=0.22\textwidth]{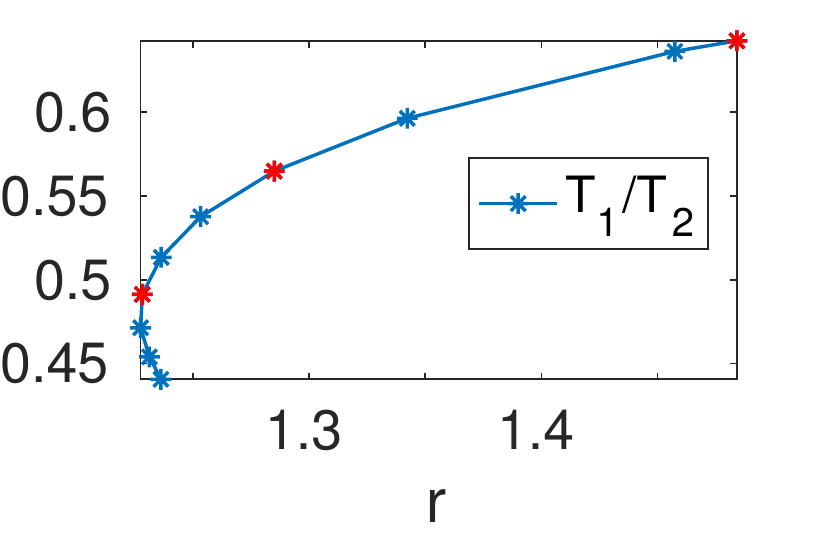}
\end{tabular}
}
&
\hs{-8mm}\raisebox{2mm}{\begin{tabular}{l}
\ig[width=0.16\textwidth]{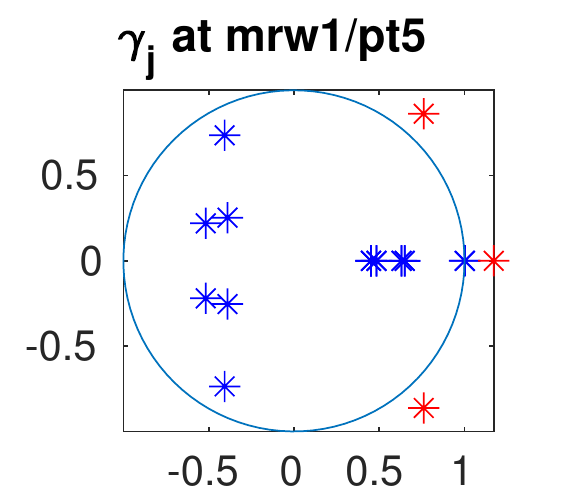}\\
\ig[width=0.16\textwidth]{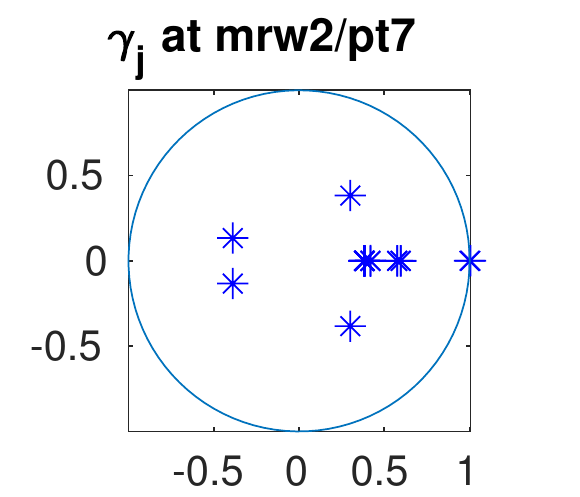}
\end{tabular}}
&
\hs{-12mm}\raisebox{2mm}{\begin{tabular}{l}
\ig[width=0.35\textwidth]{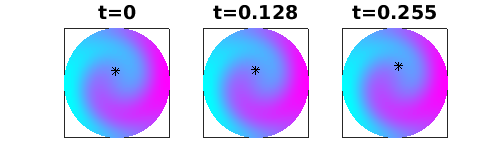}\\
\ig[width=0.35\textwidth]{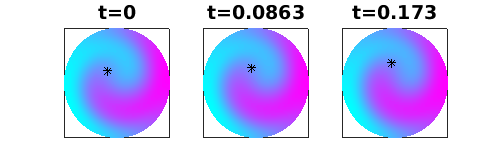}
\end{tabular}}
&
\hs{-10mm}\raisebox{2mm}{\begin{tabular}{l}
\ig[width=0.13\textwidth]{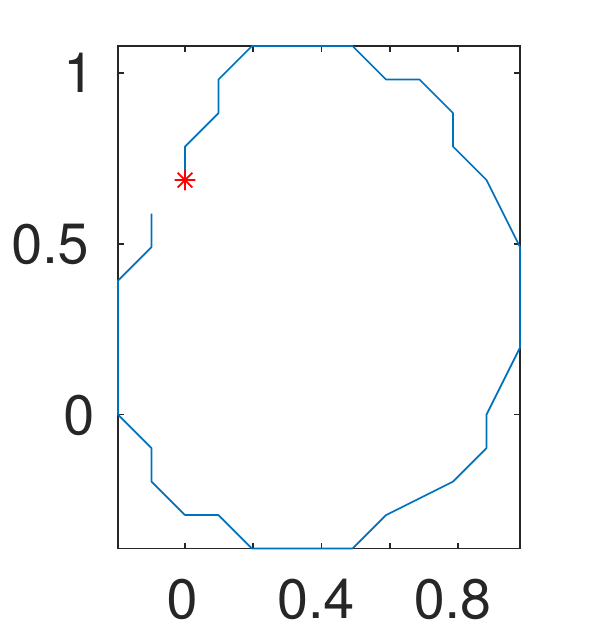}\ig[width=0.15\textwidth]{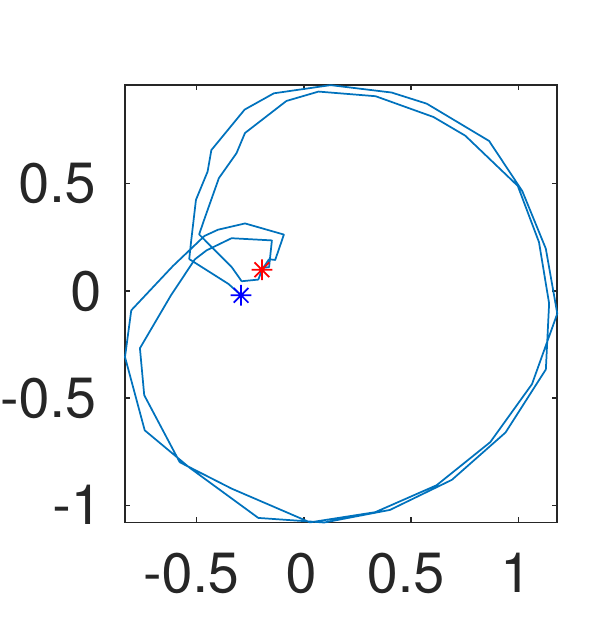}\\
\ig[width=0.14\textwidth]{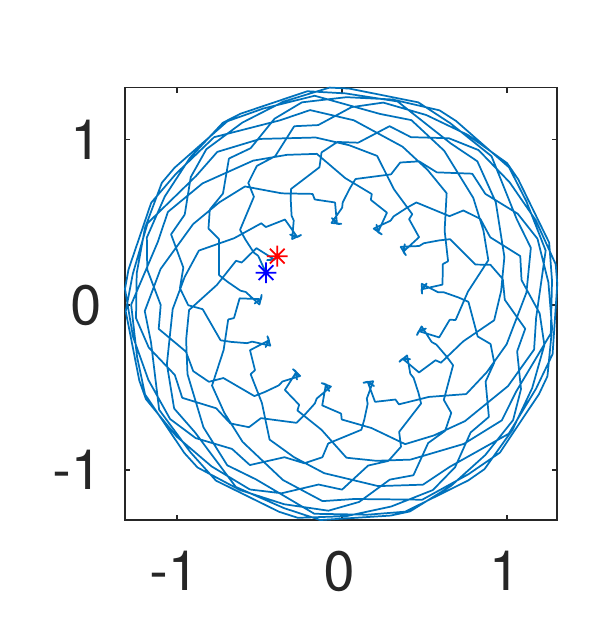}\ig[width=0.14\textwidth]{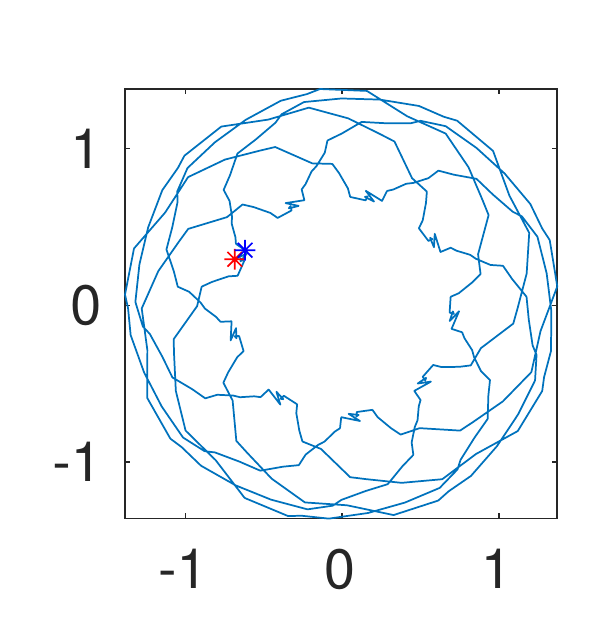}
\end{tabular}}
\end{tabular}
\ece
\vs{-5mm}
   \caption{{\small \reff{cAC2} on a disk with radius $\pi$ and homogeneous 
NBC, $(\nu,\mu,c_3,c_5,\del)=(1,-5,-1,1,1)$. (a) locally (near $\rho=0$) 
refined mesh. (b) basic BD of SWs (blue), RWs (brown), and two branches of modulated RWs, mRW1 (light brown, unstable), and mRW2 (red, stable). 
(c) example plots (snapshots 
from roughly the first quarter period) of SW and RW. (d) Period quotients  
on mRW1 (top) and mRW2 (bottom) branches, $T_1=2\pi/s$ (lab-frame time period of RWs), 
$T_2=$period on mRW2 in the rotating frame. The mRW1 branch is unstable 
and hence not discussed further, but the mRW2 branch is stable after the 
fold, and the red dots in the bottom panel correspond to the solutions used 
in (f,g). (e) Floquet spectra of selected mRW1 and mRW2. 
(g) Rotating frame plots of 
example solutions pt4 (top) and pt10 (bottom) from mRW2 branch (first $\frac 1 5$th of period), including computed tip positions. 
(g) Top: tip-path for mRW2/pt4 in comoving (left) and lab frame (right, $m=2$). 
Bottom: lab frame paths of tips for  mRW2/pt7 (left, $m=13$) 
and mRW2/pt10 (right, $m=9$). 
}  \label{cfig2}}
\end{figure}

Figure \ref{cfig2} 
(b) shows a basic BD of SWs, of RWs bifurcating at the second HBP, 
and two branches mRW1 and mRW2 of modulated RWs. 
As in Fig.~\ref{cfig1} we omit the (stable)  primary spatially 
uniform SW branch bifurcating at $r=0$. 
The RW branch (dark brown) is stable between the first and second HBP, 
where mRW1 and mRW2  bifurcate. The SW branch 
is unstable. We now focus on mRW2, which bifurcates (slightly) 
subcritically, and is stable after the fold, 
which however for efficiency we check a posteriori. (d) 
shows the ratio of periods $T_1=2\pi/s$ (lab-frame period of RWs) and $T_2$ (period on mRW2 in the rotating frame). 

To visualize the solutions in the lab frame one often uses the motion of 
the spiral tip $(\wt x(t),\wt y(t))$, where, 
if $(x(t),y(t))$ are the coordinates in the rotating frame, then 
\huga{\label{tipt} 
\bpm \wt x(t)\\\wt y(t)\epm=R_{st}\bpm x(t)\\y(t)\epm. 
}
There are different options how to define the tip coordinates 
$(x,y)(t)$.  In the far field, in radial direction 
the spirals behave like plane waves, and thus the basic 
idea is to use the intersection of isolines of $u_1$ and $u_2$ 
to define $(x,y)(t)$. Alternatively, for 'good' spirals one can 
use the maximum of $|\nabla u_1\times \nabla u_2|$, following \cite{JSW89}. 
In \cite{bark90}, a definition based on isolines of the nonlinearity $(f_1(u),f_2(u))$ has been used as numerically robust definition of the tip 
for the case of a FHN like model. 
All these definitions typically give approximately the same tip 
positions. Here we choose the intersection of the $f_1(u)=c_1$ and $f_2(u)=c_2$ 
level curves (with a small tolerances), with $c_1=c_2=0.25$, 
which for all our spirals 
gives a unique tip which makes good sense visually, see 
Fig.~\ref{cfig2}(f) for example solutions, and {\tt gettip} for 
the implementation. Moreover, with the locally 
refined spatial mesh and a sufficient time resolution, the tip paths 
$t\mapsto (x(t),y(t))$ in the comoving frame are 'reasonably smooth', 
see the first plot in (g) for an example. 
To plot the tip paths in the lab frame, we again choose a (minimal) 
$m\in\N$ such 
that $|mT_1-qT_2|<$tol for some $q\in\N$. Using tol=0.025 and \reff{tipt} 
for $0\le t\le mT_1$ generates the flower--like patterns plotted in (g). 
In particular, at pt4 (top right of (g)) we are near $1:2$ resonance. 

Similarly, the top panel of Fig.~\ref{cfig2}(d) shows the period 
quotient on the first mRW branch mRW1. On mRW1 we have a similar 
tip motion as on mRW2, but since these mRWs are unstable (see top panel 
of (e)), we skip 
further plots and discussion. 

\brem\label{cgldrem}{\rm 
(a) The RWs can be plotted using {\tt rwplot}. As in Remark \ref{tw1srem}(a), the multipliers of the RWs can be obtained 
from exponentiation of the linearization around the RW in the lab-frame. 
See below, in particular Remark \ref{spsprem}, for comments on such 
spiral spectra. \\ 
(b) With the given settings, leading to about 84000 DoF for the Hopf 
orbits, the computation of the mRW branches (10 steps on each) 
takes about 400s. The tip paths in Fig.~\ref{cfig2}(g), computed 
on the given mesh (with no interpolation involved), 
suggest that the numerics are somewhat pushed to the limit here, i.e., 
the meshes (in $x$ and $t$) may be somewhat under resolved. However, we did 
check that the results stay qualitatively the same when increasing the 
temporal resolution to 50. 
On the other hand, omitting the preparatory step of local mesh refinement near $(x,y)=0$ 
the flower patterns in (g) become (even) more ragged, but other 
quantities such as the BPs and the periods only change slightly, 
which indicates that we have a correct general picture. 
To obtain smoother tip paths, 
alternative methods allowing finer meshes should be more appropriate, e.g., time-simulation 
\cite{bark90}, or shooting methods for Hopf orbits \cite{SN10, SN16}. 
}\eex\erem

\subsection{Reaction diffusion in a disk: Demo {\tt gksspirals}}
\label{rotsec}
\def\dname{gksspirals}
Spirals like in \S\ref{cgldsec} occur in a variety of settings, i.e., in 
various excitable or oscillatory 2D RD systems. 
In particular, the cGL equation can be considered 
as a normal form near a Hopf bifurcation, and is an example of 
a so--called $\lam-\om$ system 
\huga{\label{lamom}
\pa_t u=\Delta u+\lam(u,v)u-\om(u,v)v, \quad \pa_t v=\Delta v+\lam(u,v)v+\om(u,v)u,
}
where $\lam$ and $\om$ are some functions of $u^2+v^2$. See \cite{PET94} 
for further discussion, where in particular 
$\lam(u,v)=1-u^2-v^2$ and $\om(u,v)=1+q(u^2+v^2)$, and where then 
the role of $q$ (corresponding to $\mu$ in \reff{cAC2}) as a perturbation parameter for the existence theory 
of spiral waves  and its role to determine their shape is discussed. 

In \cite[\S3.2]{hotheo}, with associated demo {\tt hopfdemos/gksspirals} we consider a two-component reaction diffusion 
system from \cite{GKS00} on the unit disk with somewhat non--standard Robin--BC, namely 
\hual{\begin{split}
\pa_t u&=d_1\Delta u+(0.5+r) u+v-(u^2+v^2)(u-\al v), \\
\pa_t v&=d_2\Delta v+rv-u-(u^2+v^2)(v+\al u), 
\end{split}\label{spir1} \\
&\pa_{{\bf n}} u+10 u=0,\quad \pa_{{\bf n}} v+0.01 v=0, \label{spirbc}
}
where ${\bf n}$ is the outer normal. In \cite[\S3.2]{hotheo}, our focus was on the computation 
of the primary SW and RW branches, and we did not compute the 
RWs as relative equilibria, and thus also skipped the computation of modulated RWs 
as relative POs. Here we include this and thus extend the presentation in \cite[\S3.2]{hotheo}. 

The eigenfunctions of the linearization 
around $(u,v)=(0,0)$ are again build from Fourier Bessel functions 
\huga{\label{n11} 
\phi_m(\rho,\vt,t)=\Re(\er^{\ri(\om t+m\vt)}J_m(q\rho)),\quad m\in\Z, 
} 
where $(\rho,\vt)$ are polar-coordinates, and with, 
due to the BC \reff{spirbc}, in general complex 
$q\in\C\setm\R$. Then the modes are growing in $\rho$, 
which is a key idea of \cite{GKS00} to find modes bifurcating 
from $(u,v)=(0,0)$ which resemble spiral waves near their core. 
The trivial solution $(u,v)=(0,0)$ is stable up to $r_1\approx -0.21$ 
where a Hopf bifurcation with angular wave number $m=0$ in \reff{n11} occurs, 
and then further Hopf-bifurcations occur with $m=1,2,0,3,\ldots$. 

Like \reff{cAC2}, \reff{spirbc} has spatial symmetry group $\rO(2)$, acting by  
rotations and reflections in $x$. The modes \reff{n11} with $m\ne 0$ 
have the symmetry of RWs, and the associated HBPs are double. 
Thus, as in \reff{bans} we use a modified branch switching 
\huga{\label{eqswibra} 
u(t)=u_0+2\eps\al\re(z_1\er^{-\ri\om_Ht}\psi_1+z_2\er^{\ri\om_H t}\psi_2)
}
with user provided $z_1,z_2\in\C$. If we apply \reff{eqswibra} at a double 
Hopf bifurcation with $(z_1,z_2)=(1,0)$, then it turns out that 
the initial guess is sufficiently close to a RW for the subsequent Newton loop 
to converge to this RW. On the other hand, $z_1=1,z_2=\ri$ at the HBPs with $m\ge 1$ yields bifurcation to SW. Alternatively, setting the PC \reff{pcr} 
(and initializing $s=0$) forces SWs, and moreover makes the continuation 
of SWs more robust.%
\footnote{The continuation of SW works also works without PC, 
because the FEM discretization 
destroys the (strict) rotational invariance, 
but it initially needs small stepsizes due to 
small eigenvalues, and in the initial continuation steps the angular phase 
of the SW pattern slightly shifts.}

First (script files {\tt cmds1} and {\tt cmds1sw}) 
we follow \cite{GKS00} and set 
 $\al=0$, $d_1=0.01$, $d_2=0.015$, and take $r$ as 
the main bifurcation parameter. Then (script file {\tt cmds2}) 
we set $\al=1$ (where $\al$ corresponds to $q$ from \cite{PET94}), 
see the comments following \reff{lamom}, let 
\huga{(d_1,d_2)=\del(0.01,0.015), 
}
and also vary $\del$ which corresponds to changing the domain size 
by $1/\sqrt\del$. Note that the operator/matrix $\pa_\phi=\Krot$ from \reff{krots}, 
and hence the rotating wave speed $s$ is independent of $\del$. 
In the end, this gives a model with distinguished 
bifurcations of spiral waves from the trivial solution $(u,v)\equiv 0$. 

\begin{figure}[ht]
{\small 
\bce
\begin{tabular}{ll}
(a) Bifurcation diagram&
(b) {\small snapshots from rotating waves rw2,rw3 (top) and rw5,rw6  (bottom)}\\
\hs{-4mm}\ig[width=48mm, height=65mm]{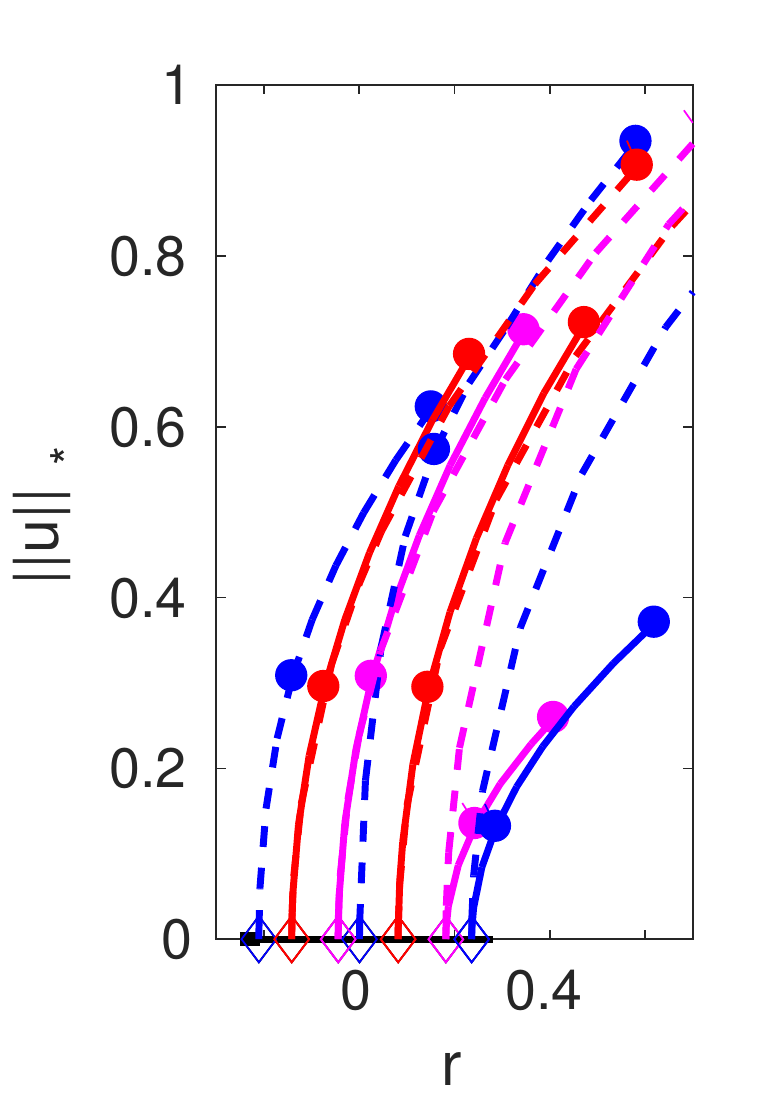}&
\raisebox{35mm}{
\begin{tabular}{l}
\hs{-10mm}\ig[width=68mm]{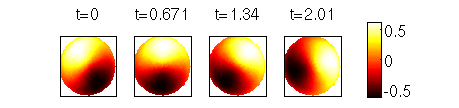}\hs{-6mm}\ig[width=68mm]{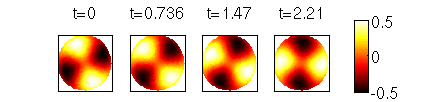}\\
\hs{-10mm}\ig[width=68mm]{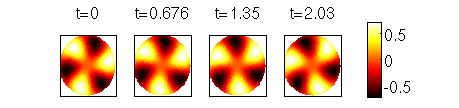}\hs{-6mm}\ig[width=68mm]{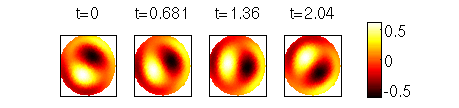}\\
\hs{-2mm}{\small (d) snapshots from standing waves sw2 and sw3} \\
\hs{-10mm}\ig[width=68mm]{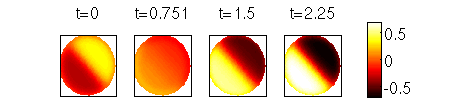}\hs{-6mm}\ig[width=68mm]{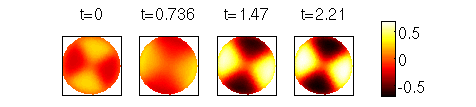}
\end{tabular}}
\end{tabular}
\ece
}

\vs{-5mm}
   \caption{{\small  
(a) Basic bifurcation diagram for rotating waves (full lines rw2, rw3, rw5, rw6, rw7), and standing waves (dashed lines sw1,\ldots, sw7) for \reff{spir1}, 
\reff{spirbc}, 10 continuation 
steps for each. On sw1 and the RW branches we mark the points 5 and 10. 
 (c)  Snapshots of $u$ from the RW branches at the last points, $t=0, T_j/9,2T_j/9,T_j/3$, with $T_j$ the actual period. (c) Snapshots of $u$ from the RW branches sw2 and sw3. Run {\tt cmds1} and {\tt cmds1sw.m} to obtain these 
(and other) plots as in \cite[\S3.2]{hotheo}. 
  \label{spf1}}}
\end{figure}

Figure \ref{spf1}(a) shows a basic bifurcation diagram for \reff{spir1}, 
\reff{spirbc}, and (b) and (c) illustrate the difference between 
rotating and standing waves. Otherwise we refer to \cite[\S3.2]{hotheo} 
for further discussion of this and other plots generated in {\tt cmds1, cmds1sw} and {\tt cmds2}. Instead, here we focus on the 
implementation, and comment on results from {\tt cmds3} about modulated RWs
(not discussed in \cite[\S3.2]{hotheo}). 

\taskip
\begin{table}[ht]\caption{Scripts and functions in {\tt hopfdemos/gksspirals}. 
\label{rottab}}
{\small 
\begin{tabular}{l|p{0.7\tew}}
script/function&purpose,remarks\\
\hline
cmds1, cmds1sw&main script for $\al=0$ in \reff{spir1}: 
continuation of RW (SW) branches\\ 
cmds2&similar to cmds1, but for $\al=1$ in \reff{spir1}\\
cmds3&continuation of RW (1--armed spiral) for $(\al,\del)=(1,0.25)$, 
branch--switching to mRW, and continuation in $\del$.\\
\hline
geo=circgeo(r,nx)&defining a circular domain\\
p=rotinit(p,nx,par)&initialization, as usual\\
sG, sGjac, nodalf&rhs, Jacobian, and nonlinearity, as usual\\
nbc&BC for \reff{spirbc}, using the convenience function {\tt gnbc}
\\\hline
auxcmds&making movies, using customized plotting from {\tt homovplot(..)}\\
hoplotrot, proplot, levplot&some additional customized plot commands\\
\end{tabular}
}
\end{table}\teskip

Table \ref{rottab} lists the pertinent files 
in {\tt gksspirals}. The general setup is very similar to {\tt cgldisk}, but a difference is that here we use the 'legacy' 
\ptool\ setting. 
As a consequence, there is no file {\tt oosetfemops.m}. 
Instead, in line 5 of {\tt rotinit.m} we define a diffusion tensor, and 
the system matrices are then generated via the \pdep\ function {\tt setfemops}. 
Additionally, we set up a function {\tt nbc.m} defining the boundary 
conditions. See Listings \ref{rl1}-\ref{rl3}. 
Afterwards, the files {\tt sG.m} and {\tt sGjac.m} encoding 
\reff{spir1},\reff{spirbc}, and 
the script files follow standard rules, where for the plotting 
we set up some customized functions derived from the default 
{\tt hoplotsol}. Additional to the discussion in 
\cite[\S3.2]{hotheo}, in C4 of {\tt cmds1} we compute the RWs 
via {\tt twswibra} as relative equilibria. For this we need to 
compute $\Krot$ in the \ptool--setting, see {\tt setfemops} for 
a pertinent local modification of the standard {\tt setfemops} function. 
To compute level curves of spirals and from these the tips of spirals 
we take advantage of \ptool\ plotting routines (see Fig.~\ref{spf2}(b) for 
an example). 
The script {\tt cmds1sw} differs 
from {\tt cmds1} only in setting the PC \reff{pcr}. 
{\tt cmds2} is similar to {\tt cmds1} but with $\al=1$ 
and with continuation in the domain size to have more pronounced spirals. 

\hulst{caption={{\small {\tt \dname/rotinit.m} (first 10 lines). 
The crucial difference to the other Hopf demos is that this is based on 
the old \ptool\ setting. This leads to changes in lines 5-8, i.e., the setup 
of the tensors and BC needed by {\tt setfemops} 
in line 10 (which does {\em not} call a function {\tt oosetfemops} if {\tt p.sw.sfem}$\ne -1$).  }},label=rl1,language=matlab,stepnumber=5, firstnumber=1,lastline=10}{\hdhome/gksspirals/rotinit.m}


\hulst{caption={{\small {\tt \dname/nbc.m}, using the 
'generalized Neumann BC' convenience function  
{\tt gnbc}. }},label=rl3,language=matlab,stepnumber=5, firstnumber=1}{\hdhome/gksspirals/nbc.m}

In {\tt cmds3} we set $(\al,\del)=(1,0.25)$ 
and start by continuation 
of the one-armed spirals as relative equilibria. Here we again strain the 
numerics with a mesh of about 3000 points, and 20 points in time 
and hence about 120000 DoF for the continuation of the mRW.%
\footnote{The finer mesh is again needed to somewhat accurately 
compute the tip paths, and, here also to compute the HBPs from 
the RWs with good accuracy.} 
Figure \ref{spf2} shows some basic results from {\tt cmds3}. The RW1 branch 
stabilizes shortly after bifurcation, where a branch of mRWs (mRW1) bifurcates.  
As expected, this is somewhat similar to the cGL case in Fig.~\ref{cfig2}, 
but there are also interesting differences: Again, the branch of RWs 
stabilizes at the HBP due to the subcritical bifurcation of mRW1,  
which itself are unstable. However, in contrast to the cGL case, here the RW1 branch then stays stable up 
to large $r>10$, i.e., upon further continuation in $r$ there does not 
seem to be a further bifurcation to mRWs. 
Moreover, there now seems to be a period locking between RW1 and mRW1, i.e., $T_1/T_2=1$ (within numerical accuracy%
\footnote{Increasing 
the temporal resolution mildly from $20$ to $40$, we obtain 
qualitatively the same behavior, and the quotient $T_1/T_2$ in 
the bottom panel of (a)
moves closer to $1$.  See also Remark \ref{cgldrem}}), 
where $T_1$ is the period on 
RW1 in the lab frame, and $T_2$ is the period on 
mRW1 in the frame rotating with speed $s$. Panels (c), (d) show tip--paths 
in the rotating and lab frames, and additionally, in (d) we also show 
the 20 largest multipliers at mRW1/pt25, of which 1 is unstable. 
(e) shows a continuation of {\tt rw1/pt15} in $\del$, and the associated 
spectrum at $\del=0.1$ (see Remark \ref{spsprem} for further comments). The period $T$ only depends weakly 
on the domain size $1/\sqrt{\del}$, as it should for reasonably 
well developed spirals for which the BC play a negligible role. 

\begin{figure}[ht]
{\small 
\bce
\begin{tabular}{lll}
(a)&(b)&(c)\\
\begin{tabular}{l}
\hs{-4mm}\ig[width=44mm]{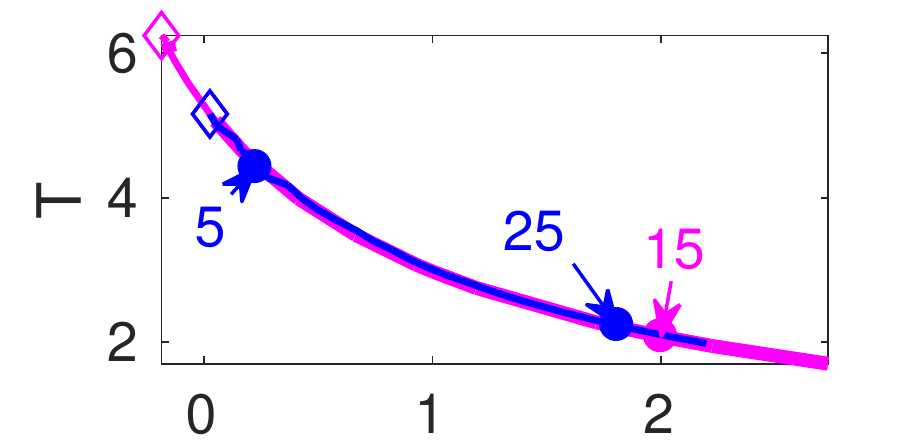}\\
\hs{-4mm}\ig[width=44mm]{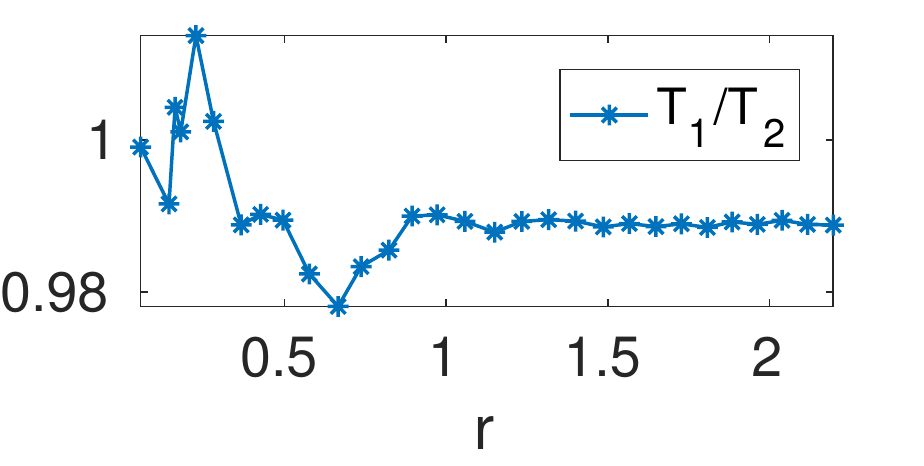}
\end{tabular}
&
\raisebox{0mm}{
\begin{tabular}{l}
\hs{-5mm}\ig[width=26mm]{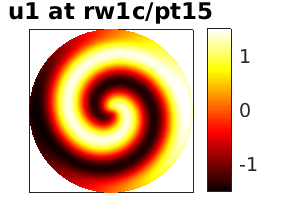}
\raisebox{-2mm}{\hs{0mm}\ig[width=26mm]{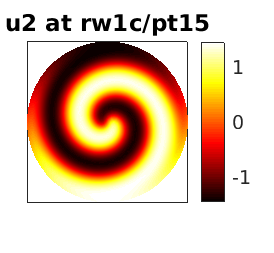}}\\
\hs{-5mm}\ig[width=34mm]{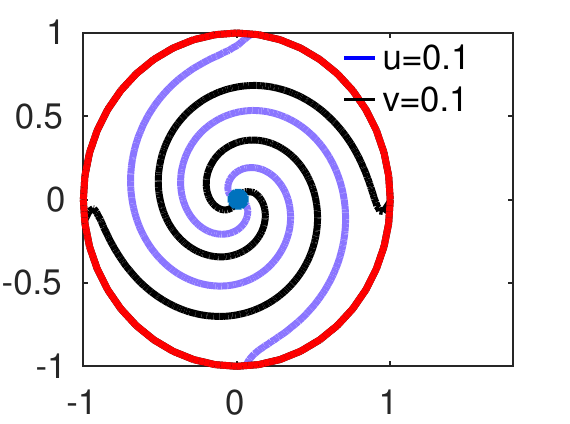}
\end{tabular}}&
\raisebox{0mm}{
\begin{tabular}{l}
\hs{-3mm}\ig[width=25mm]{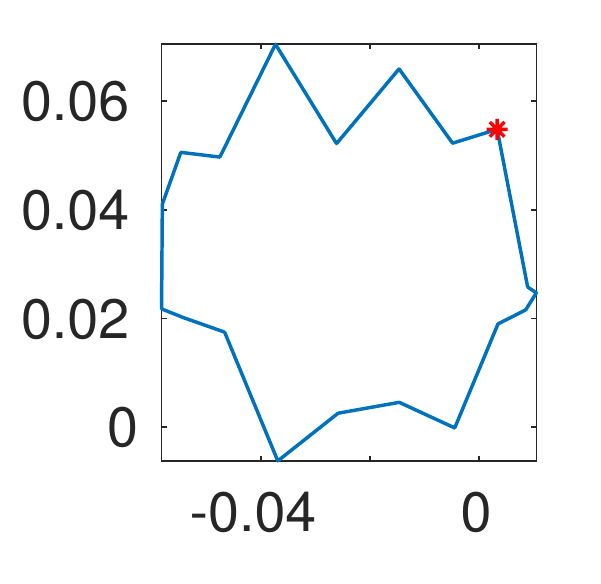}
\hs{0mm}\ig[width=25mm]{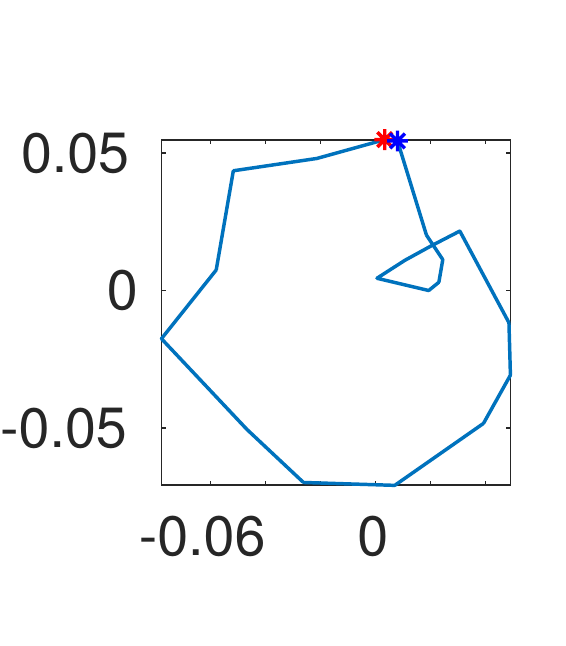}\\
\hs{-3mm}\ig[width=40mm]{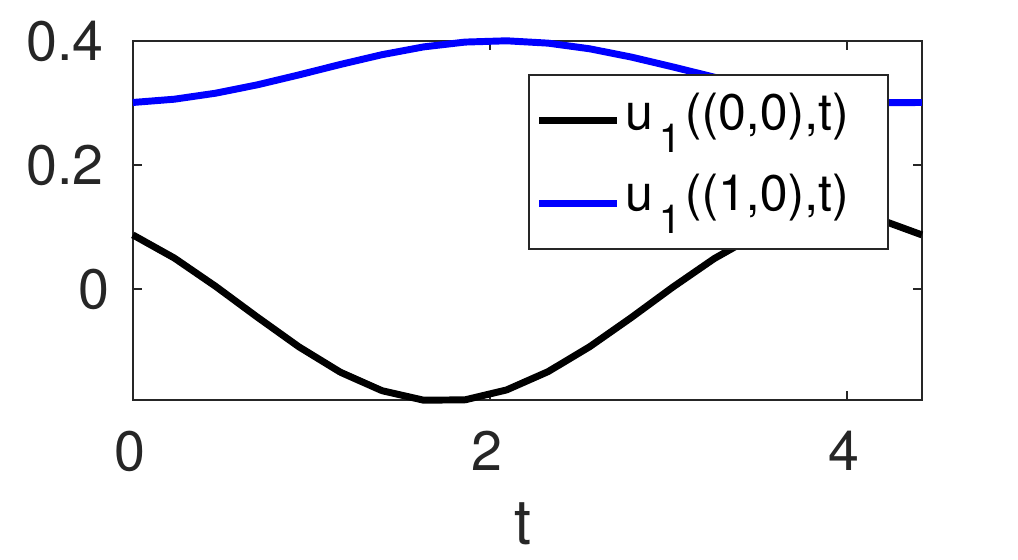}
\end{tabular}}
\end{tabular}\\
\begin{tabular}{ll}
(d)&(e)\\
\raisebox{20mm}{\begin{tabular}{l}
\hs{-3mm}\ig[width=25mm]{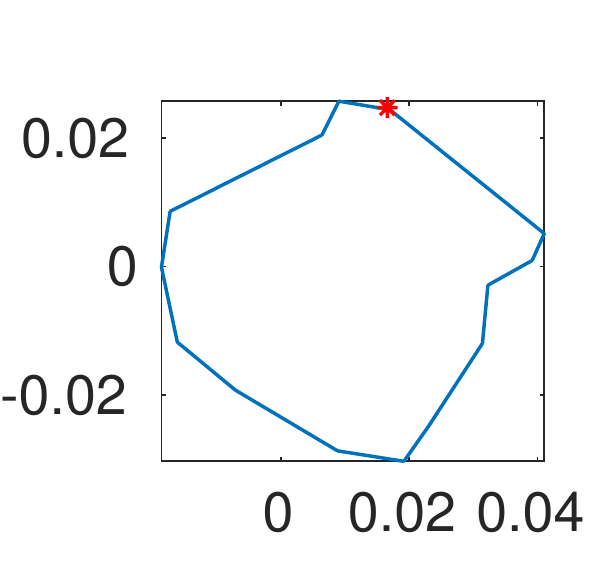}
\hs{0mm}\ig[width=27mm]{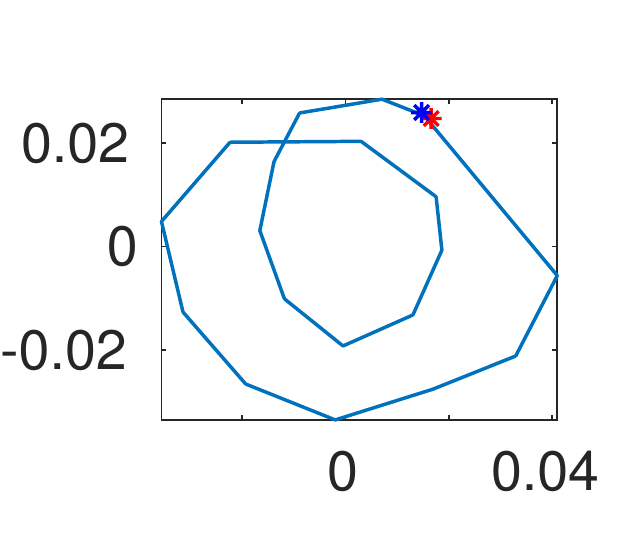}\\
\hs{0mm}\ig[width=30mm]{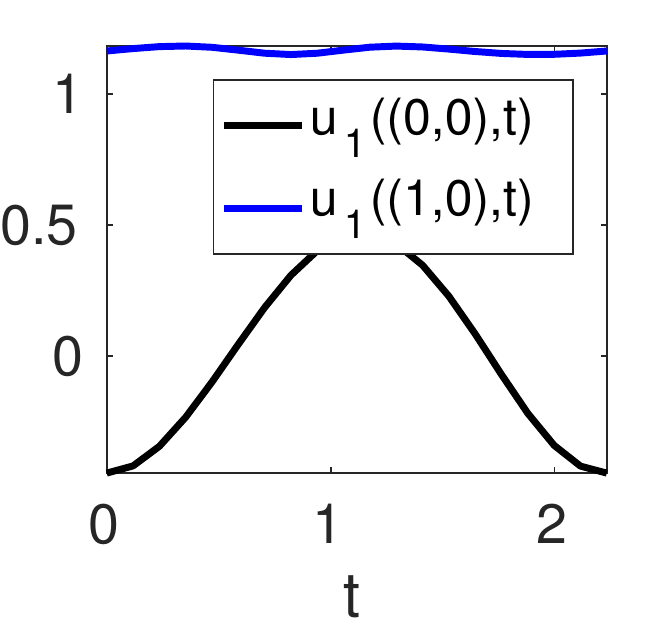}\ig[width=28mm]{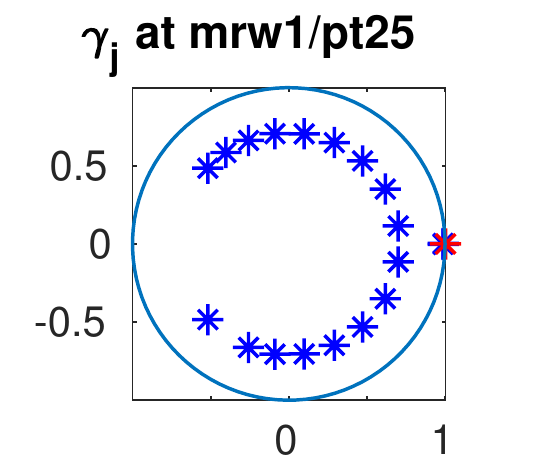}
\end{tabular}}
&\raisebox{20mm}{\begin{tabular}{l}\ig[width=28mm]{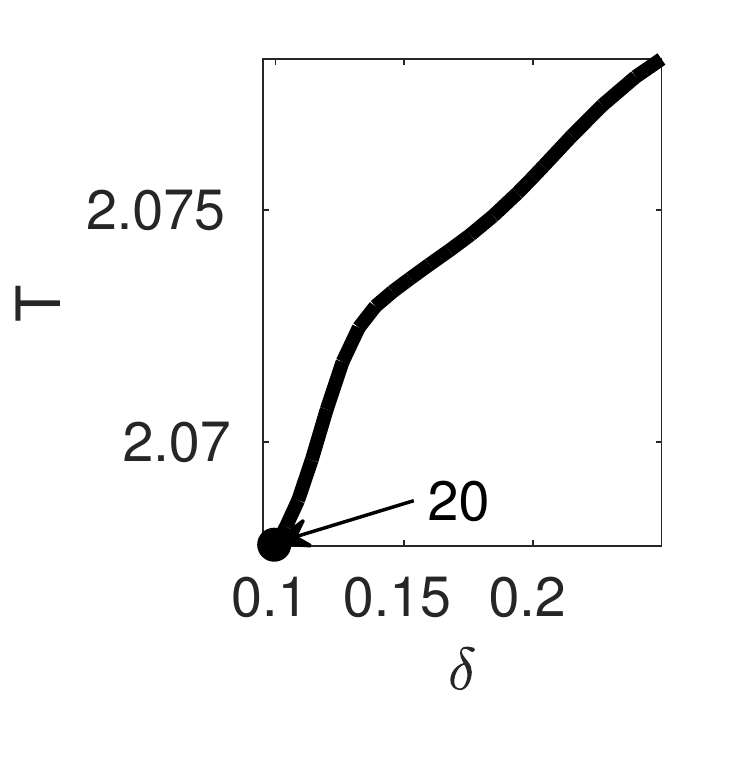}\\
\ig[width=26mm]{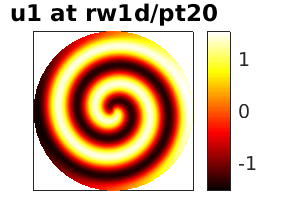}
\end{tabular}}
\ig[width=25mm]{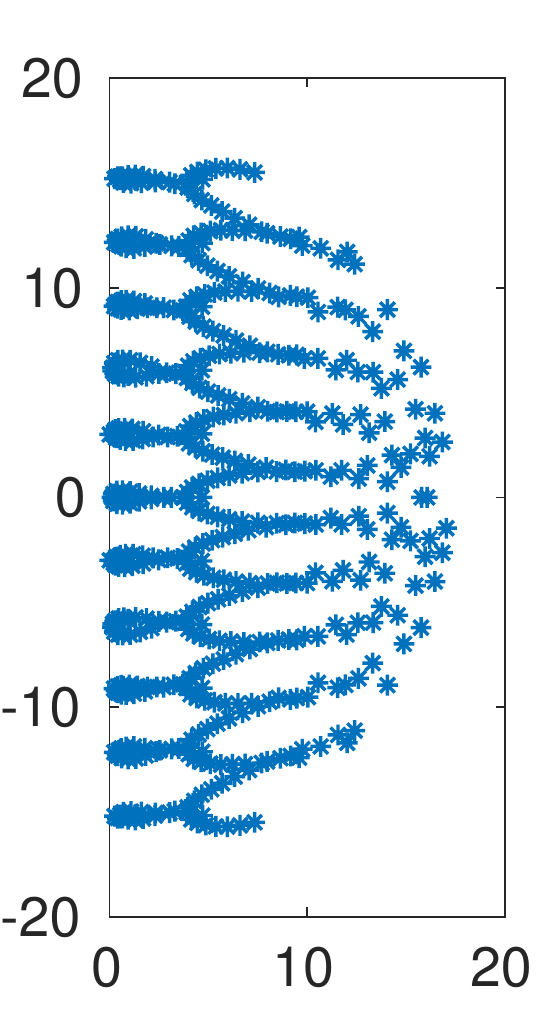}
\ig[width=25mm]{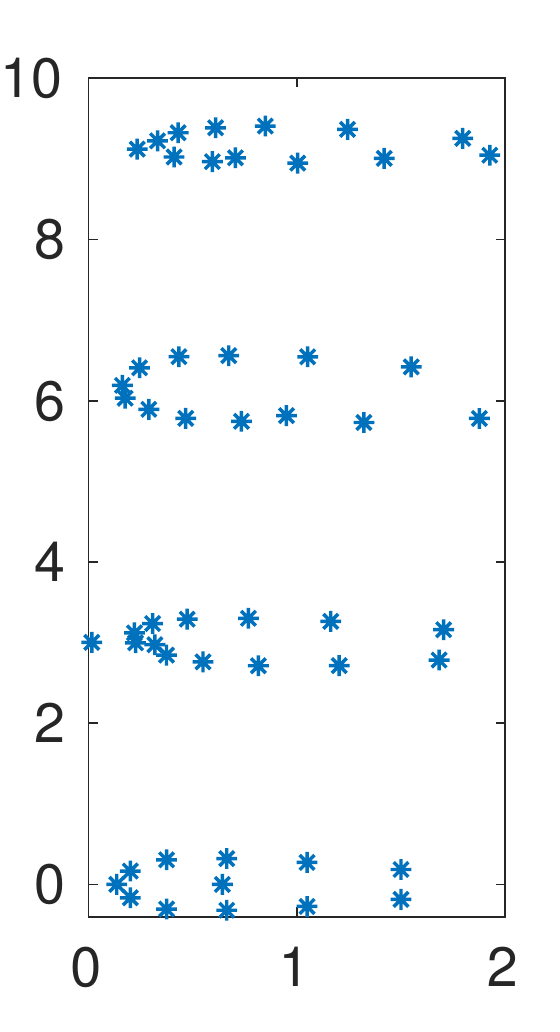}
\end{tabular}\\
\ece
}

\vs{-5mm}
   \caption{{\small  Results from {\tt cmds3.m}. 
(a) top: BD  ($\al=1, \del=0.25$) of RW1 (magenta, rotating (spiral) wave) and mRW1 (blue, modulated (meandering) spiral wave bifurcating 
at HBP1 on RW1). The periods $T_1$ 
(of RW1 in the lab frame) and $T_2$ (of mRW1 in the rotating frame) are 
equal (within numerical accuracy), and also the speeds $s_1$ (of RW1) and 
$s_2$ (of mRW1, in the average sense) agree. Bottom: quotients of periods. 
(b) snapshots from RW1/pt15, 
and contours $u_1=u_2=0.1$, used to define the spiral tip in (c,d). 
(c) top: tip paths for mRW1/pt5 in rotating frame (left) and lab frame (right); 
bottom: time series at $(x,y)=(0,0)$ (black) and $(x,y)=(1,0)$ (blue), 
illustrating that the modulation acts near the tip. (d) analogous 
data from mRW1/pt25, and additionally the 20 largest multipliers. 
(e) continuation of {\tt rw1/pt15} to $\del=0.1$, 
and illustration of the typical spiral wave spectrum: 500 smallest 
eigenvalues (middle) and zoom near imaginary axis (right); $s\approx 3.04$. 
  \label{spf2}}}
\end{figure}

\brem\label{spsprem}{\rm 
As in Remark \ref{tw1srem}(a), the multipliers $\ga_j$ of a RW can 
be obtained by exponentiation of the eigenvalues of the linearization 
in the co--rotating frame. Such spiral-wave spectra \cite{SS06, WB06, SS07} 
show some interesting 
structure, e.g., near the imaginary axis they are periodic with period 
$s$, see Fig.~\ref{spf2}(e) for an example. Heuristically, this can be explained 
as follows. In polar coordinates $(\rho,\phi)$, the linearization around 
\def\urw{u_{{\rm RW}}}
$\urw$ on the infinite plane reads 
\huga{\label{rwlin}
Lu=D(\pa_\rho^2u+\frac 1 \rho \pa_\rho u+\frac 1 {\rho^2}\pa_\phi^2 u)
+s\pa_\phi u+\pa_uf(\urw(\rho,\phi),\lam)u,
}
where following for instance \cite{SS07} we used the shorthand $D\Delta u$ for the diffusion terms 
with diffusion matrix $D$, and $f$ for the remaining terms (without spatial derivatives). 
If we then let $\rho\ra\infty$ in the eigenvalue problem $Lu=\mu u$ 
we formally obtain 
\huga{
D\pa_\rho^2u+s\pa_\phi u+\pa_uf(\urw(\rho,\phi),\lam)u=\mu u, 
}
and if $(u,\mu)$ is an eigenpair, so is $(u\er^{\ri\ell \phi},\mu+\ri s\ell)$, 
for each $s\in\Z$. This formal computation for the essential 
spectrum is explained in more 
detail in \cite{SS06, WB06, SS07} and the references therein, where 
moreover it is explained that the eigenvalues on large but finite 
domains accumulate near the so called absolute spectrum, 
and how the (formal) order 
$\CO(1/\rho)$ and $\CO(1/\rho^2)$ terms in \reff{rwlin} may generate 
isolated point spectrum. 
}\eex\erem

\subsection{Extensions: fixed period $T$, and non--autonomous cases}
\label{cglextsec}
\def\dname{cglext}
All our examples so far, and most of our research applications, 
deal with the autonomous case $M\pa_t u=-G(u)$ with no explicit 
time dependence of $G$, and with a free period $T$ which is computed as 
a part of the solution.  In the demo {\tt cglext} we explain how to 
modify the setup to 
\bci 
\item treat an explicit $t$--dependence of $G$, and/or 
\item free an additional parameter to deal with a fixed period $T$. 
\eci 
This also includes the option to compute and continue POs that are 
not generated in a Hopf bifurcation via the function {\tt poiniguess}.  
The toy problem is again \reff{cAC0}, in 1D, i.e., 
\huga{\label{pcgl} 
\pa_t u=\pa_x^2 u+(r+\ri\nu)u-(c_3+\ri \mu)|u|^2 u-c_5|u|^4u, \quad 
u=u(t,x)\in\C, 
}
but a {\em multiplicative forcing} as the 5$^{th}$ order 
coefficient $c_5$ is possibly $t$ and $x$ dependent, namely 
\huga{\label{fofu}
c_5=c_5^*+\al\tanh(10((t-\beta T){\rm mod} T))\sin x, 
}
where $\al,\beta$ are additional constants which we put into {\tt par(8)} 
and {\tt par(9)}, respectively. 
Table \ref{pcgltab} shows the new/modified files in {\tt cglext}, 
and Listings \ref{pcgll1}--\ref{pcgll4} show the pertinent files. 

For the cGL equation, a natural way to fix the period of POs is to 
free the dispersion parameter $\nu$ in \reff{pcgl} by setting  
{\tt p.hopf.ilam=2}. Then we set {\tt p.hopf.freeT=0} (default is 1), 
which tells \pdep\ to remove $T$ from the list of unknowns, and 
the column $(\pa_T \CG,0,(1-\xi_H)w_T\tau_T,\pa_T Q_H)$ from the Jacobian $\CA$, 
see \ref{news}. For the case of POs from {\tt p=hoswibra(\ldots,aux)}, 
{\tt p.hopf.freeT=0} 
should be switched on by setting {\tt aux.freeT=0}, see Listing \ref{pcgll3}. 
For the case of 
POs from {\tt poiniguess}, this can similarly be done by setting 
{\tt p.hopf.freeT} and {\tt p.hopf.ilam=2} by hand, see Listing \ref{pcgll4}. 
Independent of whether $T$ is free or not, to code a time--dependent $G$, 
we can access {\tt p.t} (and {\tt p.T}) 
filled in the Hopf interface function {\tt hosrhs}, see Listing \ref{pcgll1} 
for an example. 

\taskip
\begin{table}[ht]\caption{{\small Selected scripts and functions in 
{\tt hopfdemos/cglext} (cGLinit, oosetfemops, sG, \ldots as before)  }
\label{pcgltab}}
\bce 
\vs{-5mm}
{\small 
\begin{tabular}{l|l}
script/function&purpose,remarks\\
\hline 
cmds1\&2&scripts; cmds1 with NBCs, cmds2 with pBCs\\
fofu&the forcing function from \reff{fofu}, called in {\tt nodalf} and {\tt sGjac}. 
\end{tabular}
}
\ece 
\end{table}\teskip

\hulst{caption={{\small {\tt \dname/fofu.m}, using the current period $T$ and the current time $t$ as set into {\tt p.T} and {\tt p.t} in the interface 
function {\tt hosrhs} for calling {\tt sG} for PO continuation. }},
label=pcgll1,language=matlab} 
{\hdhome/cglext/fofu.m}
\hulst{caption={{\small {\tt \dname/nodalf.m}, 'nonlinearity' in \reff{pcgl}, 
calling {\tt fofu} in l5.}},
label=pcgll2,language=matlab}{\hdhome/cglext/nodalf.m}
\hulst{caption={{\small Selection from {\tt \dname/cmds1.m}: 
Branch switching to the 1st Hopf orbit, fixing $T=2\pi$ to the value 
computed at bifurcation, and instead freeing $\nu$.  
}},
label=pcgll3,language=matlab, linerange=11-15}
{\hdhome/cglext/cmds1simple.m}
\hulst{caption={{\small 2nd selection from {\tt \dname/cmds1.m}, here 
using {\tt poiniguess} to go to the third branch, again with fixed 
$T=2\pi$ and free $\nu$.  
}},
label=pcgll4,language=matlab, linerange=37-48}
{\hdhome/cglext/cmds1simple.m}

Figure \ref{pff1}(a,b) shows sample results from {\tt cmds1} for 
\reff{pcgl} with $(c_3,\mu,c_5,\al,\beta)=(-1,0.1,1,0.5,0.5)$ 
on $\Om=(-\pi,\pi)$ with NBCs, with as usual $r$ as primary continuation 
parameter, and $\nu$ as additional parameter as we through fix $T=2\pi$. 
We compute the first two bifurcating PO branches (black and red) 
via {\tt hoswibra} from the trivial branch, and the third (blue) 
by {\tt poiniguess} at finite amplitude. 
For non--autonomous systems we may in general expect to drop 
the (temporal) phase condition \reff{pca}, as the $t$--translational 
invariance is lost. However, for \reff{pcgl} with the forcing 
\reff{fofu} we still have continua of (approximately) ``time--shifted'' 
POs: a different phase at bifurcation yields a different action 
of the quintic terms, but this effect is small at small amplitude, 
and altogether also at larger amplitude a small phase shift yields a 
slightly different solution (including a different $\nu$). Thus we have 
continua of solutions parametrized by phase, and it is correct and 
necessary to keep the phase condition \reff{pca}. 

Figure \ref{pff1}(c,d) shows analogous sample results from {\tt cmds2} 
where instead of NBCs we use pBCs. The main difference is that now 
the second HBP is double,  and for $\al=0$ 
there bifurcate TWs and SWs, cf.~\S\ref{pcgl1dsec}. The same happpens 
for the multiplicative forcing \reff{fofu}, except that the TWs 
(blue branch in (c)) become modulated TWs. Moreover, although the forcing 
in principle breaks the translational invariance in $x$ for the 
SWs (red branch) the continuation is more robust by keeping 
the $x$--PC as in (\ref{kspc}b). At the end of {\tt cmds2.m} we again 
use {\tt poiniguess} to compute PO branches from suitable initial guesses. 

\begin{figure}[ht]
\bce{\small
\begin{tabular}{ll}
(a) NBCs, BDs, $\max$ and $\nu$ over $r$&(b) sample solutions\\
\hs{-2mm}\ig[width=0.2\tew]{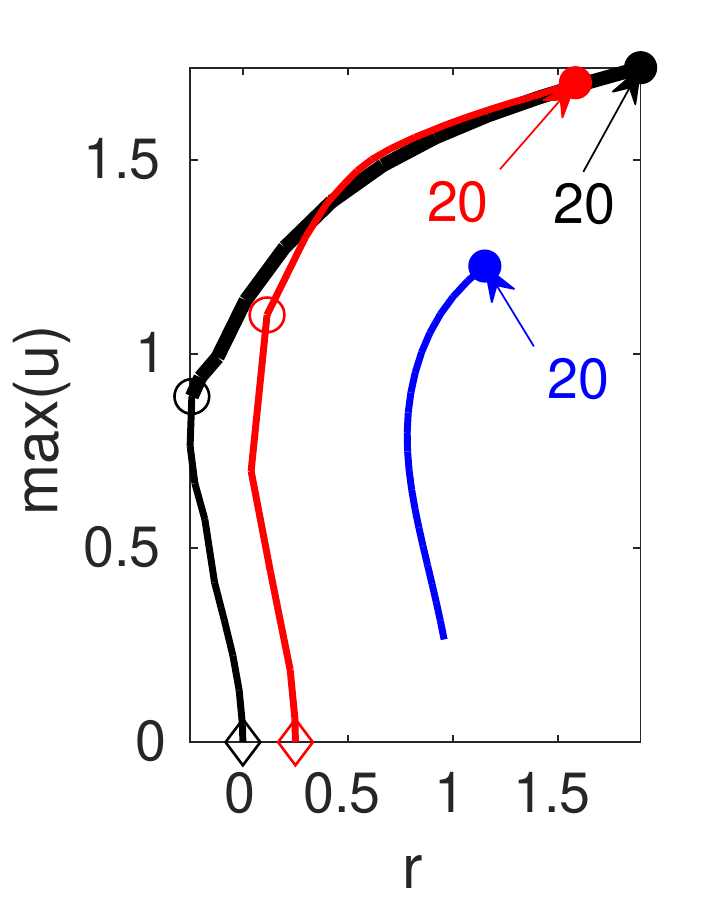}\ig[width=0.2\tew]{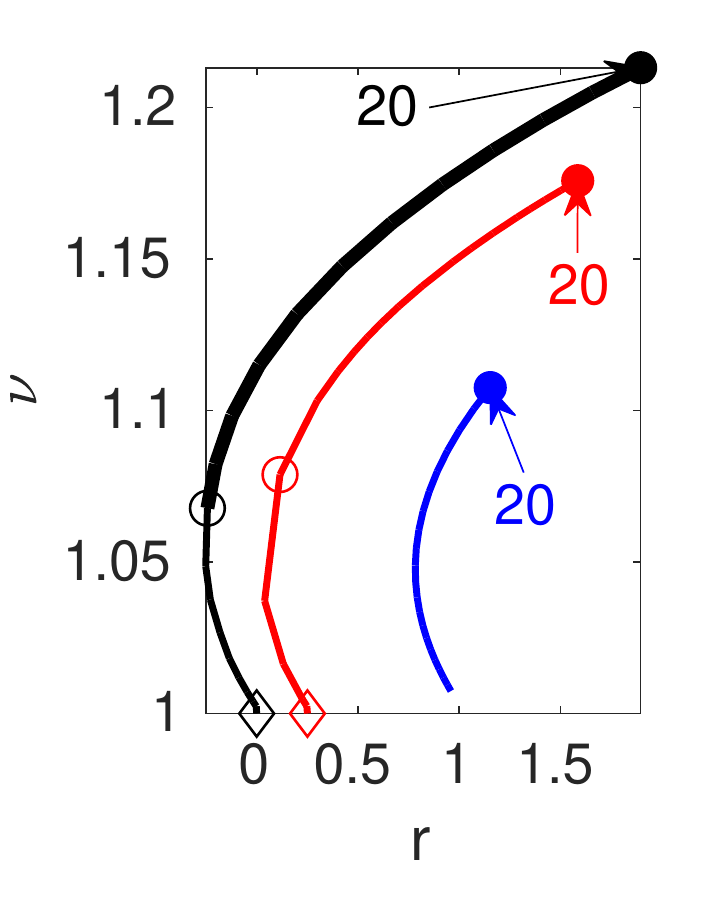}&
\hs{-4mm}\ig[width=0.2\tew]{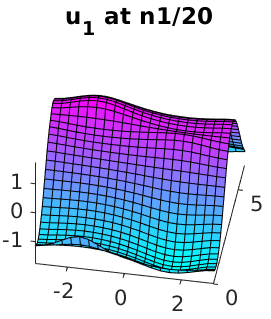}\ig[width=0.2\tew]{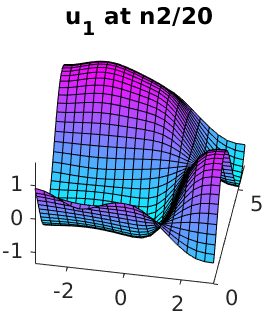}
\ig[width=0.2\tew]{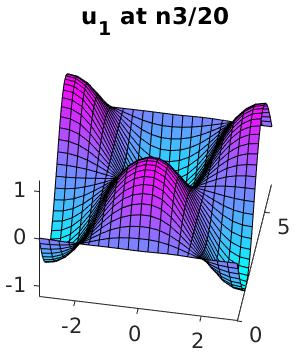}\\
(c) pBCs, BDs, $\max$ and $\nu$ over $r$&(d) sample solutions\\
\hs{-2mm}\ig[width=0.2\tew]{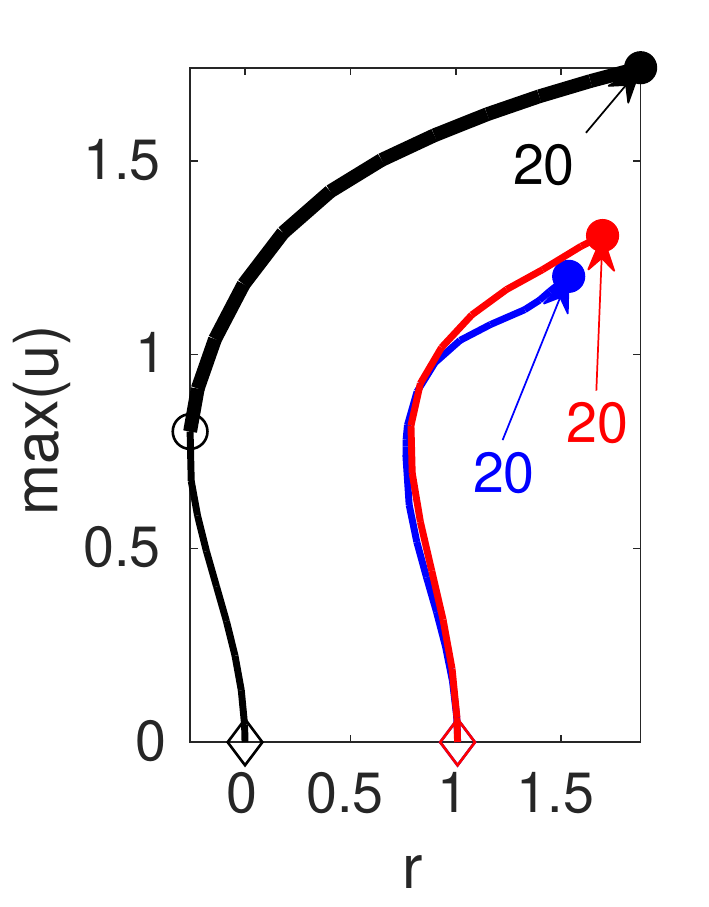}\ig[width=0.2\tew]{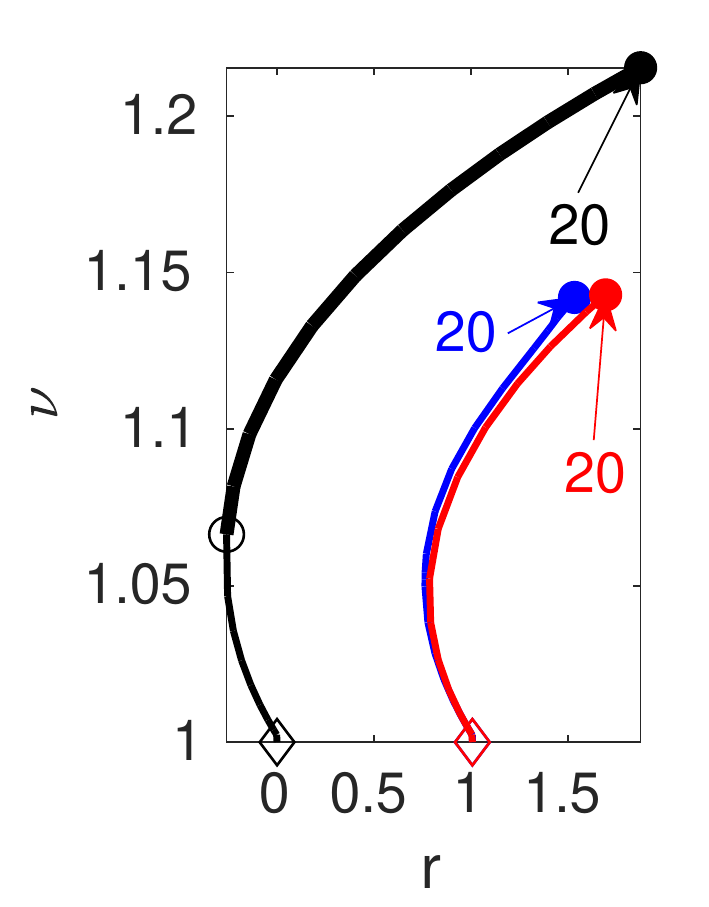}&
\hs{-4mm}\ig[width=0.2\tew]{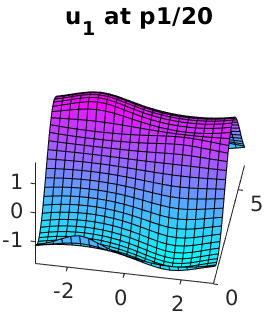}\ig[width=0.2\tew]{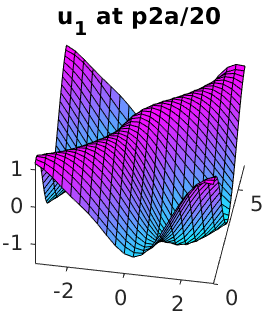}
\ig[width=0.2\tew]{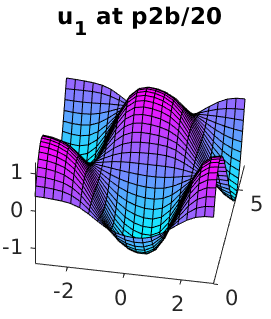}
\end{tabular}}
\ece 

\vs{-5mm}
   \caption{{\small (a,b) Sample outputs from {\tt cmds1.m}, \reff{pcgl} 
with NBCs, $\al=\beta=0.5$, fixed period $T=2\pi$ with $\nu$ as free parameter. 
First (black, n1) 
and second (n2, red) Hopf branches obtained via {\tt hoswibra}, 
and third branch (n3, blue) obtained from 
{\tt poiniguess}; sample solutions in (b). (c,d) Sample outputs from {\tt cmds1.m}, like (a,b) but with with pBCs. Consequently, the 2nd HP is double with 
standing waves and (modulated) TWs bifurcating. 
  \label{pff1}}}
\end{figure}

\appendix 
\section{Some background and formulas}
For details and background on the basic algorithms for Hopf branch point 
(HBP) detection and localization, branch--switching to and continuation 
of Hopf orbits, and Floquet multiplier computations we refer to \cite{hotheo}. 
However, in \S\ref{appa} we first briefly repeat the pertinent formulas that are implemented 
in \pdep, including the augmented systems for Hopf computations 
with constraints. We focus on the arclength parametrization 
setting ({\tt p.sw.para=4}), which is more convenient and robust than 
the 'natural parametrization' ({\tt p.sw.para=3}). 
In \S\ref{appb} we then give the formulas used for branch--switching 
as multipliers go through $\pm 1$. We mix the presentation of 
formulas with their \pdep\ implementation, and in \S\ref{appd} we give an 
overview of the used data structures and functions.

\subsection{Basics}\label{appa} 
First of all, the detection of Hopf bifurcation points (HBPs) requires the 
{\tt p.sw.bifcheck=2} setting \cite[\S2.1]{hotheo}, which is essentially controlled as follows:
{\small
\huga{
\begin{split}
&\text{{\tt p.nc.eigref(1:ne)} contains shifts near which eigenvalues are computed; guesses for these shifts}\\
&\text{ can be obtained via {\tt initeig}. However, except for \S\reff{brusec},  here we use {\tt p.nc.eigref=0}.) }
\end{split}\label{eigref}\\
\label{mus}
\begin{split}
&\text{A bisection for localization of a possible HBP is started 
if $|\re \mu|<{\tt p.nc.mu1}$ for the eigenvalue with}\\
&\text{ the smallest abs.~real part, and a HBP is accepted if $|\re\mu|<{\tt p.nc.mu2}$ at the end of the bisection.}
\end{split}
}
}
The default branch switching {\tt hoswibra} to a Hopf branch generates 
\def\dlam{\del_\lam}\def\dsd{{\rm d}s}
\huga{\label{hoswitr} 
\lam=\lam_H+\del_s\dlam
\quad u(t)=u_0+2\al\del_s\Re(\er^{-\ri \om_H t}\Psi), }
as an initial guess for a periodic solution of \reff{tformd} 
with period near $2\pi/\om$. Here $\Psi$ is the (complex) eigenvector 
associated to $\ri\om_H$, and $\del_\lam,\al$ are computed from 
the normal form 
\huga{\label{honf2}
0=r\bigg[\mu_r'(\lam_H)(\lam-\lam_H)+c_1|r|^2\bigg]. 
} 
of the bifurcation equation on the center manifold, see \cite[\S2.2]{hotheo}. After the coefficients $\del_\lam$ and $\al$ in \reff{hoswitr} are computed 
(in {\tt hogetnf}), 
$\del_s$ is chosen in such a way that the initial step length is {\tt ds} in the norm \reff{xinorm} below. The computation of $\del_\lam$, which currently 
is only implemented for semilinear systems, i.e., in FEM form 
$M\dot u=Ku-Mf(u)$, and even then is often 
not very reliable, can be skipped by calling {\tt hoswibra(\ldots,aux)} with 
{\tt aux.dlam} set to some value (usually {\tt aux.dlam=0} is the best 
choice). 

To compute Hopf orbits, after rescaling $t\mapsto Tt$ with unknown period $T$, the time evolution and periodicity condition for $u$ read 
\huga{\label{tform3} 
M\dot u=-TG(u,\lam), \quad u(\cdot,0)=u(\cdot,1), 
}
and the time-translational phase condition and arclength equation read 
\hual{\label{pca}
\phi&:=\xi_\phi\int_0^1\spr{u(t),\dot u_0(t)}\dd t\stackrel !=0, \\
\label{ala}\psi&:=
\hoxi\sum_{j=1}^m\spr{u(t_j){-}u_0(t_j),u_0'(t_j)}_\Om+
(1{-}\hoxi)\bigl[w_T(T{-}T_0)T_0'+(1{-}w_T)(\lam{-}\lam_0)\lam_0'\bigr]{-}\rds
\stackrel!=0,
}
where $T_0,\lam_0$ and $u_0$ are from the previous step, 
$\hoxi, w_T$ are weights,  $'$ denotes differentiation wrt arclength, 
and $\spr{u,v}_\Om$ stands for $\int_\Om \spr{u(x),v(x)}\dd x$, 
with $\spr{a,b}$ the standard $\R^N$ scalar product.%
\footnote{By default ({\tt p.hopf.y0dsw=2}), $\dot u_0$ in \reff{pca} is evaluated (in {\tt sety0dot})  by 2nd order finite differences; 
alternatively, {\tt p.hopf.y0dsw=1} means first order FD, and 
{\tt p.hopf.y0dsw=0} replaces $\dot u_0$ by $M\dot u_0=-TG(u_0,\lam_0)$.}
Numerically, we use $\spr{u,v}_\Om=\spr{Mu,v}$, where $M$ is the mass matrix 
belonging to the FEM mesh.  The steplength is {\tt ds} in the weighted norm 
\huga{\label{xinorm}
\|(u,T,\lam)\|_\xi=\sqrt{\hoxi\left(\sum_{j=1}^m \|u(t_j)\|_2^2\right)+(1-\hoxi)
\bigl[w_TT^2+
(1-w_T)\lam^2\bigr]}. 
}

In \reff{pca}, $\xi_\phi$ with standard setting $\xi_\phi=10$ 
({\tt p.hopf.pcfac}, see Appendix B)  is another weight, 
which can be helpful to balance the Jacobian $\CA$, see \reff{news}. 
To improve convergence of Newton loops, it sometimes turns out to be useful to set $\xi_\phi$ to somewhat larger values. 
Also, while $\dot u_0$ in \reff{pca} is in principle 
available from $M\dot u_0=-TG(u_0,\lam_0)$ (which is used 
for {\tt p.hopf.y0dsw=0}, it often appears more robust 
to explicitly approximate $\dot u_0$ via finite differences, 
for which we set {\tt p.hopf.y0dsw=2}. Finally, 
for a system with $n_H$ Hopf constraints 
$Q_H(u)=0$, and hence $n_H$ additional free parameters $w\in\R^{n_H}$, 
we also add $w_w\spr{w-w_0,w'}$ to $\psi$, where $w_w$ is a weight for 
the auxiliary parameters $w$. 

Letting $U=(u,T,\lam,w)$, and writing $\CG(U)=0$ for \reff{tform3} (see also 
\reff{gjdef}), in each continuation step we need to solve 
\huga{\label{fsac} 
H(U):=\bpm \CG(U)\\\phi(u)\\\psi(U)\\Q_H(U)\epm\stackrel{!}{=}\bpm 0\\0\\0\\0\epm
\in\R^{mn_u+2+n_H}, }
where $m$ is the number of time slices $t_j$, $j=1,\ldots,m$, $n_u$ the number of PDE unknowns at 
each time slice, and $n_H=:{\tt p.hopf.nqh}$ is the number of Hopf constraints 
(encoded in {\tt p.hopf.qfh}). To solve \reff{fsac} we use Newton's method, i.e., 
\huga{\label{news} 
U^{j+1}{=}U^{j}{-}\CA(U^j)^{-1}H(U^j), \quad 
\CA{=}\bpm \pa_u \CG&\pa_T\CG&\pa_\lam\CG&\pa_a \CG\\
\pa_u\phi&0&0&0\\
\hoxi\tau_u&(1{-}\hoxi)w_T\tau_T&(1{-}\hoxi)(1{-}w_T)\tau_\lam&w_a\tau_a\\
\pa_uQ_H&\pa_T Q_H&\pa_\lam Q&\pa_aQ_H
\epm, 
}
where of course we never form $\CA^{-1}$ but instead use {\tt p.fuha.blss} 
to solve linear systems of type $\CA U=b$. These systems are of bordered 
type, and thus it is often advantageous to use bordered system solvers, 
see \cite{lsstut}. For the case of fixed $T$ POs ({\tt p.hopf.freeT=0}), 
the second column of $\CA$ in \reff{news} is deleted, and consequently 
one additional parameter $a$ must be freed, see \S\ref{cglextsec}. 

For the time discretization we have 
\hual{\label{ydef0}
&u=(u_1,\ldots,u_m)=(u(t_1), u(t_2),\ldots, u(t_m)),\\
&\text{($m$ time slices, stored in {\tt p.hopf.y(1:p.nu,1:m)),}}\notag 
} 
where $u_m=u1$ is redundant but convenient. To assemble $\CG$ in \reff{tform3} we 
use modifications of TOM, yielding, with  $h_j=t_j-t_{j-1}$ and $u_0:=u_{m-1}$, 
\huga{\label{gjdef}
(\CG(u))_j=-h_{j-1}^{-1}M(u_j-u_{j-1})-\frac 1 2 T(G(u_j)+G(u_{j-1})), 
\quad \CG_m(u)=u_m-u_1.
}
The Jacobian is $\pa_u\CG=A_1$, where  we set, 
as it is also used for the Floquet multipliers with free $\ga$ (see \cite[\S2.4]{hotheo} and \S\ref{appb}), 
\renewcommand{\arraystretch}{1.5}
\huga{\label{fl0}
A_\ga=
\bpm M_1&0&0&0&\ldots&-H_1&0\\
-H_2&M_2&0&0&\ldots&0&0\\
0&-H_3&M_3&0&\ldots&0&0\\
\vdots&\ldots&\ddots&\ddots&\ddots&\vdots&\vdots\\
0&\ldots&\ldots&\ddots&\ddots&0&0\\
0&\ldots&\ldots&0&-H_{m-1}&M_{m-1}&0\\
-\ga\,I&0&\ldots&\ldots&\ldots&0&I
\epm,}
\renewcommand{\arraystretch}{1}
where $\ds M_j=-h_{j-1}^{-1}M-\frac 1 2 T \pa_uG(u_j)$, 
$\ds H_j=-h_{j-1}^{-1}M+\frac 1 2 T \pa_uG(u_{j-1})$, and $I$ is 
the ${n_u\times n_u}$ identity matrix.%
\footnote{In \reff{gjdef} and \reff{fl0} we assume that $M$ is regular, 
i.e., \reff{tform3} does not contain algebraic constraints as for instance 
for the 2nd order system formulation of the KS in {\tt kspbc2}; see also Remark  \ref{flremc} for this case.} 
The Jacobians $\pa_uG\in\R^{n_u\times n_u}$ in $M_j, H_j$ are computed 
as for steady state problems, e.g., via {\tt p.fuha.sGjac}, or 
by {\tt numjac} if {\tt p.sw.jac}=0 (but still locally in time). 

For $Q_H(u(\cdot,\cdot),\lam,w,a)=0$, the user 
must provide a function handle {\tt p.hopf.qfh}, similar to 
$0=Q(u,\lam,w)={\tt p.fuha.qf}$ for the steady case. 
Moreover, when switching to a Hopf branch (and if {\tt p.nc.nq} was greater 
$0$ for 
the steady continuation), the user has to drop 
the stationary constraints, i.e., reset {\tt p.nc.nq=0} and 
{\tt p.nc.ilam=p.nc.ilam(1)} to just 
the primary active parameter, while the other active 
parameters $a\in \R^{n_H}$ for $Q_H$ should be set in 
${\tt p.hopf.ilam}\in \N^{n_H}$, which now acts as a pointer to these 
'secondary' active parameters. 
Finally, the user must provide a function handle in {\tt p.hopf.qfhjac} 
to a function that returns $\pa_u Q_H$ from the last line in \reff{news}. 
On the other hand, $\pa_T \CG, \pa_\lam \CG, \pa_a\CG$, and $\pa_TQ_H, \pa_\lam Q_H$ 
and $\pa_a Q_H$ in \reff{news} are cheap from 
numerical differentiation and hence taken care of automatically.

\brem\label{flremc}{\rm 
 If \reff{tform3} contains algebraic constraint components as 
in {\tt kspbc2}, then we modify \reff{gjdef}. For instance, 
if the second component is algebraic, then we (automatically, 
in {\tt tomassemF} and {\tt tomassempbc}) set 
$(\CG(u)_2)_j=-\frac 1 2 TG(u_j)$, and accordingly also modify 
\reff{fl0}. 
} \eex\erem 

\subsection{Floquet multipliers, and bifurcation from periodic orbits}\label{appb} 
\def\gacrit{\ga_{{\rm crit}}}
The Floquet multipliers $\ga$ of a periodic orbit $u_H$ are obtained 
 from 
finding nontrivial solutions $(v,\ga)$ of the variational boundary value 
problem 
\hual{
M\dot v(t)&=-T\pa_u G(u(t))v(t),\label{fl1}\\ 
v(1)&=\ga v(0).
}
By translational invariance of \reff{tform3}, there always is the trivial multiplier $\ga_1=1$. 
Equivalently, the multipliers $\ga$ are the eigenvalues of 
the monodromy matrix $\CM(u_0)=\pa_u \Phi(u_0,T)$, where $\Phi(u_0,t)$ is 
the solution of the initial value problem \reff{tform3} with $u(0)=u_0$ 
from $u_H$. 
Thus, $\CM(u_0)$ depends on $u_0$, but the multipliers $\ga$ do not. 
$\CM(u_0)$ has the eigenvalues $1,\ga_2,\ldots,\ga_{n_u}$, where 
$\ga_2,\ldots,\ga_{n_u}$ are the multipliers of the linearized Poincar\'e map 
$\Pi(\cdot; u_0)$, which maps a point $\ut_0$ in a hyperplane 
$\Sigma$ through $u_0$ and transversal to $u_H$ to its first return to $\Sigma$, see, e.g., \cite[Theorem 1.6]{kuz04}. 
Thus, a necessary condition for the bifurcation 
from a branch $\lam\mapsto u_H(\cdot,\lam)$ of periodic orbits 
is that at some $(u_H(\cdot,\lam_0),\lam_0)$, 
additional to the trivial multiplier $\ga_1=1$ there is a 
second multiplier $\gacrit=\ga_2$ (or a complex conjugate pair $\ga_{2,3}$) 
with $|\ga_2|=1$, which generically leads to the following bifurcations 
(see, e.g., \cite[Chapter 7]{seydel} or \cite{kuz04} for more details): 
\bci
\item[(i)]\label{biftypes} $\ga_2=1$, yields a fold of the periodic orbit, or a transcritical or pitchfork bifurcation of periodic orbits.
\item[(ii)] $\ga_2=-1$, yields a period--doubling bifurcation, i.e., the bifurcation 
of periodic orbits $\ut(\cdot;\lam)$ with approximately double the period, 
$\ut(\tilde T;\lam)=\ut(0;\lam)$, $\tilde{T}(\lam)\approx 2T(\lam)$ for 
$\lam$ near $\lam_0$.  
\item[(iii)] $\ga_{2,3}=\er^{\pm\ri \vt}$ , $\vt\ne 0,\pi$, 
yields a torus (or Naimark--Sacker) bifurcation, i.e., the bifurcation 
of periodic orbits $\ut(\cdot,\lam)$ with two ``periods'' $T(\lam)$ and 
$\tilde T(\lam)$; 
if $T(\lam)/\tilde T(\lam)\not\in \Q$, then $\R\ni t\mapsto \ut(t)$ is dense in 
certain tori.  
\eci 
Numerics for (iii) are difficult even for low dimensional ODEs, but 
there are various algorithms for (i),(ii), and below we explain the simple ones 
so far used in \pdep. 
First we are interested in the computation of the multipliers. 
Using the same discretization for $v$ as for $u$, it follows that $\ga$ 
and $v=(v_1,\ldots,v_m)$ have to satisfy 
\huga{\label{fl2}
v_1=M_1^{-1}H_1v_{m-1},\quad v_2=M_2^{-1}H_2v_1,\quad \ldots,\quad v_{m-1}=M_{m-1}^{-1}H_{m-1}v_{m-2},\quad v_m=\ga v_1, 
}
for some $\ga\in\C$. 
Thus, $\CM(u_{j_0})$ can be obtained 
from certain products involving the $M_j$ and the $H_j$, for instance%
\footnote{In \cite[(2.40)]{hotheo} we used $\CM(u_{m-1})$, but $\CM(u_{1})$ 
seems more convenient for branch--switching}
\huga{\label{fl3} 
\CM(u_{1})=M_1^{-1}H_1M_{m-1}^{-1}H_{m-1}\cdots M_2^{-1}H_2. 
} 
Thus, a simple way to compute the $\ga_j$ is to compute the 
product \reff{fl3} and subsequently  (a number of) 
the eigenvalues of $\CM(u_{1})$. We call this \fla\ (Floquet Algorithm 1, 
implemented in {\tt floq}), and using 
\huga{\label{emudef}
\emu:=|\ga_1-1|} 
as a measure of accuracy we find that this works 
fast and accurately for our dissipative examples. Typically 
$\emu<10^{-10}$, although at larger amplitudes of $u_H$, 
and if there are large multipliers, this may 
go up to $\emu\sim 10^{-8}$, which is the (default) tolerance we require 
for the computation of $u_H$ itself. Thus, in the software 
 we give a warning if $\emu$ exceeds a certain tolerance 
$\fltol$.   
However, for the optimal control example in \S\ref{ocsec}, 
where we naturally have multipliers $\ga_j$ with 
$|\ga_j|>10^{30}$ and larger, 
\fla\ completely fails to compute any meaningful multipliers.

More generally, in for instance \cite{fj91,lust01} it is explained that methods 
based directly on \reff{fl3} 
\bci
\item%
may give considerable numerical errors, in particular if 
there are both, very small and very large multipliers $\ga_j$; 
\item%
discard much useful information, for instance eigenvectors 
of $\CM(u_l)$, $l\ne m-1$, which are useful for branch switching.
\eci
As an alternative, \cite{lust01} suggests to use a periodic Schur decomposition 
\cite{BGD93} to compute the multipliers (and subsequently 
pertinent eigenvectors), and gives examples that in certain cases 
this gives much better accuracy, according to \reff{emudef}. 
See also \cite{kressner01, kress06} for similar ideas and results. 
We thus provide an algorithm \flb\  (Floquet Algorithm 2, 
implemented in {\tt floqps}), 
which, based on {\tt pqzschur} from \cite{kressner01}, 
computes a periodic Schur decomposition of the matrices involved 
in \reff{fl3}, from which we immediately obtain the multipliers, see 
\cite{hotheo} for details. 

\brem\label{flremd}{\rm 
Also for $n_H>0$, the computation of Floquet multipliers is based on \reff{fl0}, 
i.e., ignores the Hopf-constraints $Q_H$. Since these constraints are 
typically used to eliminate neutral directions, ignoring these typically 
leads to Floquet multipliers close to $1$, additional to the trivial 
multiplier $1$ from (time--) translational invariance. This may lead to 
wrong stability assessments of periodic orbits (see \cite[\S4]{symtut} 
for an example), which however usually can identified by suitable 
inspection of the multipliers. 
} \eex\erem 

Here, additional to \cite{hotheo}, we give the (somewhat preliminary, 
see Remark \ref{bsprem}) algorithms for branch 
switching for the case of $\gacrit=\pm 1$. 
First, the (simple) localization 
of a BP is done in {\tt hobifdetec.m} via bisection. 
For the case $\gacrit=1$ with associated eigenvector 
$v_1$ (i.e., $\gacrit\approx 1$ but not equal to the trivial 
multiplier, and hence $v_1\ne \pa_t u_H(0)$) 
we then use \reff{fl3}, i.e., 
\huga{\label{fl2b}
v_2=M_2^{-1}H_2v_1,\quad \ldots,\quad v_{m-1}=M_{m-1}^{-1}H_{m-1}v_{m-2},
}
and additionally $v_m=v_1$, to obtain a tangent predictor 
$V=(v_1,v_2,\ldots,v_{m-1}|v_m)$ for the bifurcating branch. 

For the case $\gacrit=-1$ we double the period $T$, i.e., set 
(recall that $t_1=0$) 
\huga{t_{{\rm new}}=\frac 1 2 (t_1,t_2,t_3,\ldots,t_{m-1},1+t_1,1+t_2,\ldots,
1+t_{m-1}), 
} 
redefine $m=2m$, and use 
$(V,-V|v_1)$ with $V$ from \ref{fl2b} as predictor. The pertinent 
function is {\tt poswibra.m}, which, besides the orbit and branch point, 
and the new directory and the (initial) step length {\tt ds}, can 
take some additional arguments {\tt aux}. For instance, {\tt aux.sw=1} 
{\em forces} {\tt poswibra} to take $\gacrit$ from near $1$, while 
{\tt aux.sw=-1} takes $\gacrit$ from near $-1$. This is sometimes 
necessary because there may be $\ga_j$ close to both, $\pm 1$, 
and the bisection for the localization does not distinguish these 
(or even multipliers elsewhere near the unit circle).

\section{Data structure and function overview}
\label{appd}
The  Hopf setting naturally reuses and extends the stationary \pdep\ setting 
explained in, e.g., \cite{qsrc}. As usual, here we assume that the problem 
is described by the struct {\tt p}, and for convenience list the main 
subfields of {\tt p}  in Table \ref{tab1a}. 

\taskip
{\small
\begin{longtable}{|p{0.07\tew}|p{0.38\tew}|
p{0.07\tew}|p{0.38\tew}|}
\caption{{Main fields in the structure {\tt p} for steady problems, see \cite{qsrc} 
for more details.
}\label{tab1a}}
\endfirsthead\endhead\endfoot\endlastfoot
\hline
field&purpose&field &purpose\\\hline
fuha&{\bf fu}nction {\bf ha}ndles, e.g., fuha.G, \ldots& nc&
{\bf n}umerical {\bf c}ontrols, e.g., nc.tol, \ldots\\
sw&{\bf sw}itches such as sw.bifcheck,\ldots&
sol&values/fields calculated at runtime\\
pdeo&OOPDE data if OOPDE is used&
mesh&mesh data (if the {\tt pdetoolbox} is used)\\
plot&switches and controls for plotting&file&switches etc for file output\\
bel&controls for lssbel (bordered elimination)&ilup&controls for lssAMG 
(ilupack parameters)\\ 
\hline
usrlam&\multicolumn{3}{p{0.85\tew}|}{vector of user set target 
values for the primary parameter, default usrlam=[];} \\
mat&\multicolumn{3}{p{0.75\tew}|}
{problem matrices, in general data that is not saved to file}\\
\hline
\end{longtable}
}\teskip

For the continuation of time-periodic orbits, the field {\tt p.hopf} contains the pertinent data; it is typically created and filled by calling {\tt p=hoswibra(..)}. This inter alia  calls {\tt p=hostanparam(p,aux)}, which 
can be used as a reference for the default values of the Hopf parameters. 
The unconstrained Hopf setting does not need any user setup additional to the 
functions such as {\tt p.fuha.sG, p.fuha.sGjac} already needed for stationary problems. In case of Hopf constraints, 
the user has to provide two function handles in {\tt p.hopf.qfh} and 
{\tt p.hopf.qfhder} to functions which compute $Q_H$ from \reff{fsac} 
and the last row of $\CA$ from \reff{news}, respectively. 
The only changes of the core {\tt p2p} library concern some queries 
whether we consider a Hopf problem, in which case basic 
routines such as {\tt cont} call a Hopf version, i.e., {\tt hocont}. 
Table \ref{hotab1} gives an overview of {\tt p.hopf}, and 
Table \ref{hotab2} lists the main Hopf orbit related functions. 

\taskip
{\small 
\begin{longtable}{|p{0.16\tew}|p{0.84\tew}|}
\caption{Standard (and additional, at bottom) entries in {\tt p.hopf}.\label{hotab1}}
\endfirsthead\endhead\endfoot\endlastfoot
\hline
field&purpose\\
\hline
y&for {\tt p.sw.para=4}: unknowns in the form ($u=(u_1,\ldots,u_m)=(u(t_1), u(t_2),\ldots, u(t_m))$,
\text{ ($m$ time slices, y=$n_u\times m$ matrix);} \\
&for {\tt p.sw.para=3}: $y$ augmented by $\tilde y$ and $T,\lam$ 
($(2n_u{+}2)\times m$ matrix), see \cite{hotheo}.  \\
y0d&for {\tt p.sw.para=4}: $M\dot u_0$ for the phase condition \reff{pca}, 
 ($n_u\times m$ matrix); \\
&for {\tt p.sw.para=3}: $M\dot u_0(0)$  for the phase condition 
\cite[(36)]{hotheo},
($2n_u{+}2$ vector). \\
y0dsw&(for {\tt p.sw.para=4}) controls how $\dot u_0$ in \reff{pca} is computed: 0 for using the PDE \reff{tform3}, 2 for using FD (default).\\
pcfac&weight for the phase condition \reff{pca}, default=10\\
tau&tangent, for {\tt p.sw.para=4}, $(mn_u+2+n_H)\times 1$ vector, see third line in \reff{news}\\
ysec&for {\tt p.sw.para=3}, secant between two solutions $(y_0,T_0,\lam_0)$, $(y_1,T_1,\lam_1)$, $(2n_u{+}2)\times m$ matrix\\
sec&if sec=1, then use secant tau (instead of tangent) predictor for {\tt  p.sw.para=4}\\
t, T, lam&time discretization vector, current period and param.value\\
xi,tw,qw&weights for the arclength \reff{ala}, xi=$\hoxi$, tw=$w_T$, qw=$w_a$; \\
x0i&index for plotting $t\mapsto u(\vx({\tt x0i})$;\\
plot&aux. vars to control hoplot during hocont; see the description of {\tt hoplot}; default plot=[]\\
wn&struct containing the winding number related settings for {\tt initeig}\\
tom&struct containing TOM settings, including the mass matrix $M$\\
jac&switch to control assembly of $\pa_u\CG$. jac=0: numerically (only recommended for testing); jac=1: via {\tt hosjac}.  Note 
that for {\tt p.sw.jac=0} the local matrices $\pa_u G(u(t_j))$ are 
obtained via {\tt numjac}, but this is still much faster than using {\tt p.hopf.jac=0}. \\
flcheck&0 to switch off multiplier-comp.~during cont., 1 to use 
{\tt floq}, 2 to use {\tt floqps}\\
nfloq&\# of multipliers (of largest modulus) to compute (if flcheck=1)\\
fltol& tolerance for multiplier $\ga_1$ 
(give warning if $|\ga_1-1|>${\tt p.hopf.fltol})\\
muv1,muv2&vectors of stable and unstable multipliers, respectively\\
pcheck&if 1, then compute residual in hoswibra (predictor check)\\
bisec&\# of bisection used for BP localizations\\
\hline
&\hs{30mm}Additional entries in case of Hopf constraints\\
\hline 
ilam&pointer to the $n_Q$ additional active parameters in {\tt p.u(p.nu+1:end)}; the pointer to the primary active parameter is still in {\tt p.nc.ilam(1)}. \\
qfh, qfhder&(handles to) functions returning the Hopf constraints $Q_H(U)$ and 
the derivatives (last lines of $H$ in \reff{fsac} and $\CA$ in \reff{news}, 
respectively). \\
\hline
spar,kwnr&index of speed parameter, and spatial wave-nr for TW continuation, 
set in twswibra\\
\hline
\end{longtable}
}

{\small 
\begin{longtable}{|p{0.2\tew}|p{0.8\tew}|} 
\caption{Overview of main functions related to Hopf bifurcations and 
periodic orbits; see {\tt p2phelp} for argument lists and more comments.
\label{hotab2}}
\endfirsthead\endhead\endfoot\endlastfoot
\hline
name&purpose, remarks\\
\hline
hoswibra&branch switching at Hopf bifurcation point, see comments below \\
twswibra&branch switching at Hopf bifurcation point to Traveling Wave branch 
(which is continued as a rel.equilibrium)\\ 
hoswipar&change the active continuation parameter, see also swiparf\\
hoplot&plot the data contained in hopf.y. Space-time plot in 1D; in 2D and 3D: 
snapshots at (roughly) $t=0$, $t=T/4$, $t=T/2$ and $t=3T/4$; see also hoplotf; \\
initeig&find guess for $\om_1$; see also initwn\\
floq&compute {\tt p.hopf.nfloq} multipliers during continuation 
({\tt p.hopf.flcheck}=1)\\
floqps&use periodic Schur to compute (all) multipliers during continuation 
(flcheck=2)\\
floqap, floqpsap&a posteriori versions of floq and floqps, respectively\\
hobra&standard--setting for p.fuha.outfu (data on branch), 
template for adaption to a given problem \\
hostanufu&standard function called after each continuation step\\
plotfloq&plot previously computed multipliers\\
\hline
hpcontini&init Hopf point continuation\\
hpcontexit&exit Hopf point continuation\\
hpjaccheck&check user implementation of \reff{HPC4} against finite differences\\
hploc&use extended system \reff{HPC1} for Hopf point localization\\
hobifdetec&detect bifurcations {\em from} Hopf orbits, and use bisection for localization, based on multipliers \\
poswibra&branch switching {\em from} Hopf orbits\\
hobifpred&compute predictor for branch switching {\em from} Hopf orbits\\
hotintxs&time integrate \reff{tformd} from the data 
contained in p.hopf and u0, with output of $\|u(t)-u_0\|_\infty$, 
and saving $u(t)$ to disk at specified values\\
tintplot*d&plot output of hotintxs; $x{-}t$--plots for *=1, else snapshots 
at specified times \\
hopftref&refine the $t$-mesh in the arclength setting at user 
specified time $t^*$\\
hogradinf&convenience function returning the time $t^*$ where $\|\pa_t u(\cdot,t)\|_{\infty}$ is maximal; may be useful for hopftref. \\
\hline
initwn&init vectors for computation of initial guess for spectral shifts $\om_j$\\
hogetnf&compute initial guesses for dlam, al from the normal form  coefficients of bifurcating Hopf branches, see \cite[(16)]{hotheo}\\
hocont&main continuation routine; called by cont if p.sol.ptype$>$2\\
hostanparam&set standard parameters \\
hostanopt&set standard options for hopf computations\\
hoinistep&perform 2 initial steps and compute secant, used if {\tt p.sw.para=3}\\
honloopext,honloop&the arclength Newton loop, and the Newton loop with fixed $\lam$\\
sety0dot&compute $\dot u_0$ for the phase condition \reff{pca}\\
tomsol&use TOM to compute periodic orbit in p.sw.para=3 setting.\\
tomassemG&use TOM to assemble $\CG$, see \cite[(26)]{hotheo}; see also 
{\tt tomassem, tomassempbc} \\
gethoA&put together the extended Jacobian $\CA$ from \cite[(27)]{hotheo}\\
hopc&the phase condition $\phi$ from\cite[(19)]{hotheo},  
and $\pa_u\phi$.  \\
arc2tom, tom2arc&convert arclength data to tomsol data, e.g., to call tomsol for mesh adaptation. tom2arc to go back. \\
ulamcheckho&check for and compute solutions at user specified values in 
p.usrlam\\
poiniguess&generate initial guess for periodic orbit continuation based on 
{\tt hopf} data structures, but without need of a HBP. Alternative to 
{\tt hoswibra}, see \S\ref{cglextsec}. \\
\hline
hosrhs,hosrhsjac&interfaces to p.fuha.G and p.fuha.Gjac at fixed $t$, internal functions called by tomassempbc, together with hodummybc\\
horhs,hojac&similar to hosrhs, horhsjac, for p.sw.para=3, see also {\tt hobc} 
and {\tt hobcjac}\\
\hline
\end{longtable}
}\teskip

Besides {\tt cont}, the functions 
{\tt initeig}, {\tt hoswibra}, {\tt poswiba}, {\tt hoplot}, {\tt twswibra}, {\tt floqap}, 
{\tt floqplot}, {\tt hotintxs}, {\tt tintplot*d},  and {\tt hopftref}  are 
most likely to be called directly by the user, 
and {\tt hobra} (the branch data) and 
{\tt hostanufu} (called after each continuation step) are likely to be adapted by the user. 
As usual, all functions in Table \ref{hotab2} can be most easily overloaded 
by copying them to the given problem directory and modifying them there.  

In {\tt p=hoswibra(dir,fname,ds,para,varargin)}, 
the auxiliary argument {\tt aux=varargin$\{2\}$} \linebreak ({\tt varargin$\{1\}$} is 
the new directory)  
can for instance have the following 
fields:  
\bci
\item {\tt aux.tl}=30: 
number of (equally spaced) initial mesh-points in $t\in[0,1]$ 
(might be adaptively refined by TOM for {\tt p.sw.para=3}, or via {\tt hopftref} or 
{\tt uhopftref}  for {\tt p.sw.para=4}). 
\item {\tt aux.hodel}=1e-4: used for the finite differences in {\tt hogetnf}. 
\item {\tt aux.al, aux.dlam} (no preset): these can be used to pass a guess 
for $\al$ and $\del_\lam$ in \reff{hoswitr} and thus circumvent {\tt hogetnf}; 
useful for quasilinear problems and for problems with constraints (for 
which {\tt hogetnf} will not work), or more generally when the computation of $\al, \del_\lam$ via {\tt hogetnf} seems to give unreliable results. 
\item {\tt aux.z}: The coefficients $z_{1},\ldots,z_m$ in the ad hoc modification 
\reff{eqswibra} of \reff{hoswitr} used for Hopf points of higher multiplicity 
$m$. 
\eci 

For the other functions listed above we refer to the m-files for 
description of their arguments, and to the demo directories for 
examples of usage and customization.


\renewcommand{\refname}{References}
\renewcommand{\arraystretch}{1.05}\renewcommand{\baselinestretch}{1}
\small
\bibliographystyle{alpha}
\newcommand{\etalchar}[1]{$^{#1}$}

\end{document}